\documentclass{art_ar}

\usepackage{amssymb}

\usepackage{pstricks, pst-plot, pst-node, pst-tree}
\usepackage{array}
\usepackage{fancyhdr}
\usepackage{longtable}
\usepackage{tabularx}
\usepackage{rotating}

\usepackage{makeidx}
\usepackage{amsmath}
\usepackage{theorem}
\usepackage[mathscr]{eucal}
\usepackage{enumerate}
\usepackage[dvips]{epsfig}
\usepackage{float}  

\theoremstyle{plain}
\newtheorem{Def}{Definition}[section]
\newtheorem{Sat}[Def]{Proposition}
\newtheorem{The}[Def]{Theorem}
\newtheorem{Kor}[Def]{Corollary}
\newtheorem{Bem}[Def]{Remark}
\newtheorem{Lem}[Def]{Lemma}

\newtheorem{Bsp}[Def]{Example}

\newcommand{\Id}{\operatorname{Id}}

\begin{document}

\begin{frontmatter}


\title{Rooted Tree Analysis for Order Conditions
of Stochastic Runge-Kutta Methods for the Weak Approximation of
Stochastic Differential Equations}
\author{Andreas R{\"o}{\ss}ler}
\ead{roessler@mathematik.tu-darmstadt.de}
\address{Darmstadt University of Technology, Fachbereich Mathematik, Schlossgartenstr.7,
D-64289 Darmstadt, Germany}





\begin{abstract}
A general class of stochastic Runge-Kutta methods for the weak
approximation of It{\^o} and Stratonovich stochastic differential
equations with a multi-dimensional Wiener process is introduced.
Colored rooted trees are used to derive an expansion of the
solution process and of the approximation process calculated with
the stochastic Runge-Kutta method. A theorem on general order
conditions for the coefficients and the random variables of the
stochastic Runge-Kutta method is proved by rooted tree analysis.
This theorem can be applied for the derivation of stochastic
Runge-Kutta methods converging with an arbitrarily high order.
\end{abstract}
\end{frontmatter}
%
%
%
%
\section{Introduction} \label{Sec:Introduction}
In recent years many numerical methods have been proposed for the
approximation of stochastic differential equations (SDEs), see
e.g.\ \cite{KP99}, \cite{LRS02}, \cite{Mac01}, \cite{MilTret04},
\cite{New91}, \cite{Ruem82} and \cite{Ta90}. Mainly, numerical
methods for strong and for weak approximations can be
distinguished. While strong approximations focus on a good
approximation of the path of a solution, weak approximations are
applied if a good distributional approximation is needed. In
Section~\ref{Sec:RK-Methods-for-general-SDE-systems} of the
present paper, a class of stochastic Runge-Kutta (SRK) methods for
the weak approximation of It{\^o} and Stratonovich SDEs is
introduced. As in the deterministic setting, order conditions for
SRK methods are calculated by comparing the numerical solution
with the exact solution over one step assuming exact initial
values. Therefore, the actual solution of the SDE and the
numerical approximation process have to be expanded by a
stochastic Taylor series. However, even for low orders such
expansions become much more complex than in the deterministic
setting where it is already a lengthy task. In order to handle
this task in an easy way, a rooted tree theory based on three
different kinds of colored nodes is established in
Section~\ref{Sec:Rooted-Tree-Theory}, which is a generalization of
the rooted tree theory due to Butcher~\cite{Butcher87}. Thus,
colored trees are applied in
Section~\ref{Sec:Taylor-Expansion-SDE} and
\ref{Sec:Expansion-SRK-method} to give a representation of the
solution and the approximation process calculated with the SRK
method in order to allow a rooted tree analysis of order
conditions. A similar approach with two different kinds of nodes
has been introduced by Burrage~\&~Burrage~\cite{BuBu96},
\cite{BuBu00a} for a SRK method converging in the strong sense as
well as in Komori~et~al.~\cite{KoMiSu97} for ROW-type schemes for
Stratonovich SDEs. Finally, the main
Theorem~\ref{Ito-St-Theo-conv-cond-tree-main:Wm} presented in
Section~\ref{Sec:general-cond-order-conv-p} immediately yields all
order conditions for the coefficients and the random variables of
the introduced SRK method such that it converges with an
arbitrarily given order in the weak sense. As a result of this
theorem, the lengthy calculation and comparison of Taylor
expansions can be
avoided. \\ \\
Let $(\Omega, \mathcal{F}, P)$ be a probability space with a
filtration $(\mathcal{F}_t)_{t \geq 0}$ and let
$\mathcal{I}=[t_0,T]$ for some $0 \leq t_0 < T < \infty$. We
consider the solution $(X_t)_{t \in \mathcal{I}}$ of either a
$d$-dimensional It{\^o} stochastic differential equation system
\begin{equation} \label{St-lg-sde-ito-1}
    d X_t = a(t,X_t) \, dt + b(t,X_t) \, dW_t
\end{equation}
or a $d$-dimensional Stratonovich stochastic differential equation
system
\begin{equation} \label{St-lg-sde-strato-1}
    d X_t = a(t,X_t) \, dt + b(t,X_t) \, \circ dW_t .
\end{equation}
Let $X_{t_0} = x_0 \in \mathbb{R}^d$ be the
$\mathcal{F}_{t_0}$-measurable initial condition such that for
some $l \in \mathbb{N}$ holds $E(\|X_{t_0}\|^{2l})<\infty$ where
$\| \cdot \|$ denotes the Euclidean norm if not stated otherwise.
Here, $W=((W_t^1, \ldots, W_t^m))_{t \geq 0}$ is an
$m$-dimensional Wiener process w.r.t.\ $(\mathcal{F}_t)_{t \geq
0}$. SDE~(\ref{St-lg-sde-ito-1}) and (\ref{St-lg-sde-strato-1})
can be written in integral form
\begin{equation} \label{Intro-Ito-St-SDE1-integralform-Wm}
    X_t = x_0 + \int_{t_0}^t a(s,X_s) \, ds + \sum_{j=1}^m
    \int_{t_0}^t b^j(s,X_s) \, * dW_s^j
\end{equation}
for $t \in \mathcal{I}$, where the $j$th column of the $d \times
m$-matrix function $b=(b^{i,j})$ is denoted by $b^j$ for $j=1,
\ldots,m$. Here, the second integral w.r.t.\ the Wiener process
has to be interpreted either as an It{\^o} integral in case of SDE
(\ref{St-lg-sde-ito-1}) or as a Stratonovich integral in case of
SDE (\ref{St-lg-sde-strato-1}), which is indicated by the
asterisk. \\ \\
The solution $(X_t)_{t \in \mathcal{I}}$ of a Stratonovich SDE
with drift $a$ and diffusion $b$ is also a solution of an It{\^o}
SDE as in (\ref{St-lg-sde-ito-1}) and therefore also a diffusion
process, however with the modified drift
\begin{equation}
    \tilde{a}^i(t,x) = a^i(t,x) + \frac{1}{2} \sum_{j=1}^d
    \sum_{k=1}^m b^{j,k}(t,x)
    \frac{\partial b^{i,k}}{\partial x^j} (t,x)
\end{equation}
for $i=1, \ldots, d$ and provided that $b$ is sufficiently
differentiable, i.e.\
\begin{equation}
    \begin{split} \label{St-lg-sde-Strato-Ito-equals}
    X_t & = X_{t_0} + \int_{t_0}^t a(s,X_s) \, ds + \int_{t_0}^t b(s,X_s) \,
    \circ dW_s \\
    & = X_{t_0} + \int_{t_0}^t \tilde{a}(s,X_s) \, ds + \int_{t_0}^t b(s,X_s) \, dW_s.
    \end{split}
\end{equation}
The solution of the stochastic differential
equation~(\ref{Intro-Ito-St-SDE1-integralform-Wm}) is sometimes
denoted by $X^{t_0,X_{t_0}}$ in order to emphasize the initial
condition.
We suppose that the drift $a : \mathcal{I} \times \mathbb{R}^d
\rightarrow \mathbb{R}^d$ and the diffusion $b : \mathcal{I}
\times \mathbb{R}^d \rightarrow \mathbb{R}^{d \times m}$ are
measurable functions satisfying a linear growth and a Lipschitz
condition
\begin{alignat}{3}
    & \| a(t,x) \| + \| b(t,x) \| \leq C \left( 1 + \| x \|
    \right) \label{growth1} \\
    & \| a(t,x)-a(t,y) \| + \| b(t,x)-b(t,y) \| \leq C \left\| x-y
    \right\| \label{lip1}
\end{alignat}
for all $x,y \in \mathbb{R}^d$ and all $t \in \mathcal{I}$ with
some constant $C>0$. Then the conditions of the Existence and
Uniqueness Theorem are fulfilled for the It{\^o}
SDE~(\ref{St-lg-sde-ito-1}) (see, e.g., \cite{KS99}). If the
conditions also hold with $a$ replaced by the modified drift
$\tilde{a}$ in the It{\^o} SDE, then the Existence and Uniqueness
Theorem also applies to the Stratonovich
SDE~(\ref{St-lg-sde-strato-1}).
\\ \\
In the following, let $C_P^l(\mathbb{R}^d, \mathbb{R})$ denote the
space of $l$ times continuously differentiable functions $g \in
C^l(\mathbb{R}^d, \mathbb{R})$ for which all partial derivatives
up to order $l$ have polynomial growth. That is, for which there
exist constants $K > 0$ and $r \in \mathbb{N}$ depending on $g$,
such that $| \partial_x^i g(x) | \leq K \, (1+\|x\|^{2r})$ holds
for all $x \in \mathbb{R}^d$ and any partial derivative
$\partial_x^i g$ of order $i \leq l$. \\ \\
Let $\mathcal{I}_h = \{t_0, t_1, \ldots, t_N \}$ be a
discretization of the time interval $\mathcal{I}=[t_0,T]$ such
that
\begin{equation}
    0 \leq t_0 < t_1 < \ldots < t_N = T
\end{equation}
and define $h_n = t_{n+1} - t_n$ for $n=0,1, \ldots, N-1$ with the
maximum step size
\begin{equation*}
    h = \max_{0 \leq n \leq N-1} h_n.
\end{equation*}
In the following, we consider a class of approximation processes
of the type
$    Y^{t,x}(t+h) = A(t,x,h; \xi) $
where $\xi$ is a random variable or in general a vector of random
variables, with moments of sufficiently high order, and $A$ is a
vector valued function of dimension $d$. We write $Y_n =
Y^{t_0,X_{t_0}}(t_n)$ and we construct the sequence
\begin{equation} \label{St-lg-method}
    \begin{split}
    Y_0 &= X_{t_0} \\
    Y_{n+1} &= A(t_n, Y_n, h_n;
    \xi_n), \qquad n=0,1, \ldots,N-1,
    \end{split}
\end{equation}
where $\xi_0$ is independent of $Y_0$, while $\xi_n$ for $n \geq
1$ is independent of $Y_0, \ldots, Y_n$ and $\xi_0, \ldots,
\xi_{n-1}$. Then we can define weak convergence with some order
$p$ of an approximation process.
\begin{Def}
    A time discrete approximation process $Y$ converges weakly
    with order~$p$ to the solution process $X$ of SDE~(\ref{St-lg-sde-ito-1})
    or SDE~(\ref{St-lg-sde-strato-1}) as $h \to 0$
    if for each $f \in C_P^{2(p+1)}(\mathbb{R}^d, \mathbb{R})$
    there exists a constant $C_f$, which does not depend on $h$, and
    a finite $h_0 > 0$ such that
    \begin{equation}
        | E(f(X_{t})) - E(f(Y(t))) | \leq C_f \, h^p
    \end{equation}
    holds for each $h \in \, ]0,h_0[$ and $t \in \mathcal{I}_h$.
\end{Def}
Since we are interested in calculating a global approximation
converging in the weak sense with some desired order $p$, we make
use of the following theorem due to Milstein~(1986)~\cite{Mil86}
which is stated with an appropriate notation.
\begin{The} \label{St-lg-theorem-main}
    Let $X$ be the solution of SDE~(\ref{St-lg-sde-ito-1}) or of
    SDE~(\ref{St-lg-sde-strato-1}).
    Suppose the following conditions hold:
    \begin{enumerate}[(i)]
    \item \label{St-lg-cond1} the coefficients $a^i$ in the case of SDE
    (\ref{St-lg-sde-ito-1}), $\tilde{a}^i$ in the case of SDE
    (\ref{St-lg-sde-strato-1}) and $b^{i,j}$ are continuous,
    satisfy a Lipschitz condition (\ref{lip1}) and belong to
    $C_P^{2(p+1)}(\mathbb{R}^d, \mathbb{R})$ with respect to $x$ for $i=1, \ldots,
    d$, $j=1, \ldots, m$,
    \item \label{St-lg-cond2} for sufficiently large $r$ (specified below)
    the moments $E(\|Y_n\|^{2r})$ do exist and are uniformly bounded with respect
    to $N$ and $n=0,1, \ldots, N$,
    \item \label{St-lg-cond3} assume that for all $f \in C_P^{2(p+1)}(\mathbb{R}^d,
    \mathbb{R})$ the following {\emph{local error estimation}}
    \begin{equation} \label{St-lg-local-error-estimation}
        |E(f(X^{t,x}(t+h))) - E(f(Y^{t,x}(t+h)))| \leq K(x) \, h^{p+1}
    \end{equation}
    is valid for $x \in \mathbb{R}^d$, any $h>0$ with $t, t+h \in \mathcal{I}$ and $K \in
    C_P^0(\mathbb{R}^d, \mathbb{R})$.
    \end{enumerate}
    Then for all $N$ and all $n=0,1, \ldots, N$ the following
    {\emph{global error estimation}}
    \begin{equation} \label{St-lg-global-error-estimation}
        | E(f(X^{t_0,X_{t_0}}(t_n))) - E(f(Y^{t_0,X_{t_0}}(t_n))) | \leq C \, h^p
    \end{equation}
    holds for
    all $f \in
    C_P^{2(p+1)}(\mathbb{R}^d, \mathbb{R})$, where $C$ is a
    constant and where $h$ is the maximum step size of the
    corresponding discretization $\mathcal{I}_h$,
    i.e.\ the method (\ref{St-lg-method}) has order of accuracy $p$
    in the sense of weak approximation.
\end{The}
A proof of Theorem~\ref{St-lg-theorem-main} can be found in
\cite{MilTret04}, \cite{Mil95}, \cite{Mil86} and \cite{Roe03}.
Lemma~\ref{St-lg-Lem-bound-1} gives sufficient conditions such
that condition~(\ref{St-lg-cond2}) of
Theorem~\ref{St-lg-theorem-main} (see also \cite{MilTret04},
\cite{Mil95}, \cite{Mil86}) holds.
\begin{Lem} \label{St-lg-Lem-bound-1}
    Suppose that for $Y_n$ given by (\ref{St-lg-method}) and $h<1$ the conditions
    \begin{alignat}{2}
        & \| E(A(t_n,x,h; \xi_n)-x) \| \leq C_1 (1+\|x\|) \, h, \label{St-lg-Lem-bound-cond1} \\
        & \| A(t_n,x,h; \xi_n)-x \| \leq M(\xi_n) (1+ \|x\|) \,
        h^{1/2} \label{St-lg-Lem-bound-cond2}
    \end{alignat}
    hold where $M(\xi_n)$ has moments of all orders, i.e.\
    $E((M(\xi_n))^i) \leq C_2$, $i \in \mathbb{N}$, with constants
    $C_1$ and $C_2$ independent of $h$. Then for
    every even number $2r$ the expectations $E(\|Y_n\|^{2r})$
    exist and are uniformly bounded with respect to $N$ and $n=1,
    \ldots,N$, if only $E(\|Y_0\|^{2r})$ exists.
\end{Lem}
%
%
%
%
\section{A Class of Stochastic Runge-Kutta Methods}
\label{Sec:RK-Methods-for-general-SDE-systems}
%
%
In the following a class of {\emph{stochastic Runge-Kutta
methods}} is introduced for the approximation of both It{\^o} and
Stratonovich stochastic differential equation systems w.r.t.\ an
$m$-dimensional Wiener process. In order to preserve the most
possible generality, the considered class of stochastic
Runge-Kutta methods is of type (\ref{St-lg-method}) and has the
following structure
\begin{equation}
\begin{split}
    Y_0 &= x_0 \\
    Y_{n+1} &= A(t_n, Y_n, h_n; \theta_{\nu}(h_n) : \nu \in \mathcal{M})
\end{split}
\end{equation}
where $\mathcal{M}$ is an arbitrary finite set of multi-indices
with $\kappa = |\mathcal{M}|$ elements and $\theta_{\nu}(h)$, $\nu
\in \mathcal{M}$, are some suitable random variables. For the weak
approximation of the solution $(X_t)_{t \in \mathcal{I}}$ of the
$d$-dimensional SDE
system~(\ref{Intro-Ito-St-SDE1-integralform-Wm}), considered
either with respect to It{\^o} or Stratonovich calculus, the
general class of $s$-stage stochastic Runge-Kutta methods is given
by
\begin{alignat}{3} \label{St-srk-method-m}
    Y_0 &=& \,\, x_0 \notag \\
    Y_{n+1} &=& \,\, Y_n &+ \sum_{i=1}^s z_i^{(0,0)} \, \,
    a \left( t_n+c_i^{(0,0)} h_n, H_i^{(0,0)} \right) \\
    && &+ \sum_{i=1}^s \sum_{k=1}^m \sum_{\nu \in \mathcal{M}} z_i^{(k,\nu)}
    \, \,
    b^k \left(t_n+c_i^{(k,\nu)} h_n, H_i^{(k,\nu)} \right) \notag
\end{alignat}
for $n=0,1, \ldots, N-1$ with
\begin{alignat}{3}
    {H_i^{(0,0)}} & = & \,\, Y_n &+ \sum_{j=1}^{s} Z_{ij}^{(0,0),(0,0)}
    \, \, a \left( t_n+c_j^{(0,0)} h_n, H_j^{(0,0)} \right) \notag \\
    & & &+ \sum_{j=1}^{s} \sum_{r=1}^m \sum_{\mu \in \mathcal{M}}
    Z_{ij}^{(0,0),(r,\mu)} \, \, b^r \left( t_n+c_j^{(r,\mu)} h_n, H_j^{(r,\mu)} \right) \notag \\
    {H_i^{(k,\nu)}} & = & \,\, Y_n &+ \sum_{j=1}^{s} Z_{ij}^{(k,\nu),(0,0)}
    \, \, a \left( t_n+c_j^{(0,0)} h_n, H_j^{(0,0)} \right) \notag \\
    & & &+ \sum_{j=1}^{s} \sum_{r=1}^m \sum_{\mu \in \mathcal{M}}
    Z_{ij}^{(k,\nu),(r,\mu)} \, \, b^r \left( t_n+c_j^{(r,\mu)} h_n, H_j^{(r,\mu)} \right) \notag
\end{alignat}
for $i=1, \ldots,s$, $k=1, \ldots, m$ and $\nu \in \mathcal{M}$,
where
\begin{alignat*}{5}
    z_i^{(0,0)} & = \alpha_i \, h_n &\qquad \qquad \qquad
    z_i^{(k,\nu)} & = \sum_{\iota \in \mathcal{M}} {\gamma_i^{(\iota)}}^{(k,\nu)} \,
    {\theta_{\iota}}(h_n) \notag \\
    Z_{ij}^{(0,0),(0,0)} & = A_{ij}^{(0,0),(0,0)} \, h_n
    &\qquad
    Z_{ij}^{(0,0),(r,\mu)} & =
    \sum_{\iota \in \mathcal{M}}
    {B_{ij}^{(\iota)}}^{(0,0),(r,\mu)} \, {\theta_{\iota}}(h_n)
    \notag \\
    Z_{ij}^{(k,\nu),(0,0)} & = A_{ij}^{(k,\nu),(0,0)} \, h_n
    &\qquad
    Z_{ij}^{(k,\nu),(r,\mu)} & =
    \sum_{\iota \in \mathcal{M}}
    {B_{ij}^{(\iota)}}^{(k,\nu),(r,\mu)} \, {\theta_{\iota}}(h_n) \notag
\end{alignat*}
for $i,j = 1, \ldots, s$. Here $\alpha_i,
{\gamma_i^{(\iota)}}^{(k,\nu)}, A_{ij}^{(k,\nu),(0,0)},
{B_{ij}^{(\iota)}}^{(k,\nu),(r,\mu)} \in \mathbb{R}$ are the
coefficients of the SRK method and as usual the weights can be
defined by
\begin{equation} \label{SRK-weights-condition}
    c^{(0,0)} = A^{(0,0),(0,0)} e , \qquad \qquad c^{(k,\nu)} =
    A^{(k,\nu),(0,0)} e ,
\end{equation}
with $e=(1, \ldots, 1)^T$. If
$A_{ij}^{(k,\nu),(0,0)}={B_{ij}^{(\iota)}}^{(k,\nu),(r,\mu)}=0$
for $j \geq i$ and $0 \leq k \leq m$, $\iota, \nu \in \mathcal{M}
\cup \{0\}$, then (\ref{St-srk-method-m}) is called an explicit
SRK method, otherwise it is called implicit. The class of SRK
methods introduced above can be characterized by an extended
Butcher array
\begin{equation}
{\setlength{\extrarowheight}{6pt}
\begin{tabular}{c|c|c|c|c}
    $c^{(0,0)}$ & $A^{(0,0),(0,0)}$ & ${B^{(\iota_1)}}^{(0,0),(r,\mu)}$ & $\ldots \ldots \ldots \ldots$ &
    ${B^{(\iota_{\kappa})}}^{(0,0),(r,\mu)}$ \\
    \cline{1-5}
    $c^{(k,\nu)}$ & $A^{(k,\nu),(0,0)}$ & ${B^{(\iota_1)}}^{(k,\nu),(r,\mu)}$ & $\ldots \ldots \ldots \ldots$ &
    ${B^{(\iota_{\kappa})}}^{(k,\nu),(r,\mu)}$ \\
    \hline
    & $\alpha^T$ & ${{\gamma^{(\iota_1)}}^{(k,\nu)}}^T$ & $\ldots \ldots \ldots \ldots$
    & ${{\gamma^{(\iota_{\kappa})}}^{(k,\nu)}}^T$
\end{tabular}
}
\end{equation}
for $k,r=1, \ldots, m$ and $\iota_i, \nu, \mu \in \mathcal{M}$ for
$1 \leq i \leq \kappa$. We assume that the random variables
$\theta_{\nu}(h_n)$ satisfy the moment condition
\begin{equation} \label{St-SRK-moment-condition-1:Wm}
    \begin{split}
    E \left( \theta_{\nu_1}^{p_1}(h_n) \cdot \ldots \cdot
    \theta_{\nu_{\kappa}}^{p_{\kappa}}(h_n) \right)
    = O \left( h_n^{(p_1 + \ldots + p_{\kappa} ) / 2} \right)
    \end{split}
\end{equation}
for all $p_i \in \mathbb{N}_0$ and $\nu_i \in \mathcal{M}$, $1
\leq i \leq \kappa$. The moment condition ensures a contribution
of each random variable having an order of magnitude
$O(\sqrt{h})$. This condition is in accordance with the order of
magnitude of the increments of the Wiener process. Further, the
moment condition is necessary for the estimates of the reminder
terms of the Taylor expansion of the SRK approximation presented
in Section~\ref{Sec:general-cond-order-conv-p}.
\\ \\
Some SRK schemes which belong to the introduced general class of
SRK methods can be found in \cite{Roe03a}, \cite{Roe04a} and
\cite{Roe03}. Further, many Runge-Kutta type schemes proposed in
recent literature like in \cite{KP99}, \cite{KoMiSu97},
\cite{Mac01} or \cite{ToVA02} are covered. Usually, the set
$\mathcal{M}$ may consist of some multi-indices $(j_1,\ldots,j_l)$
with $0 \leq j_i \leq m$ for $i=1, \ldots, l$ and the random
variables may be chosen as multiple It{\^o} or Stratonovich
integrals of type $I_{(j_1,\ldots,j_l)}/h^q$ or $J_{(j_1, \ldots,
j_l)}/h^q$, depending on the calculus that is
used. \\ \\
For example, the SRK scheme RI1WM due to
R{\"o}{\ss}ler~\cite{Roe03} for the It{\^o}
SDE~(\ref{St-lg-sde-ito-1}) in the case of $d \geq 1$ and $m \geq
1$ with $\mathcal{M} = \{\{j_1\}, \{j_1,j_2\} : 1 \leq j_1, j_2
\leq m\}$ is defined by (\ref{St-srk-method-m}) with
\begin{alignat*}{5}
    z_i^{(0,0)} & = \alpha_i \cdot h_n &\qquad \qquad
    z_i^{(k,l)} & = {\gamma_i^{(k)}}^{(k,l)} \, \hat{I}_{(k)}
    + {\gamma_i^{(k,l)}}^{(k,l)} \, \tfrac{\hat{I}_{(k,l)}}{\sqrt{h_n}} \\
    Z_{ij}^{(0,0),(0,0)} & =
    A_{ij}^{(0,0),(0,0)} \cdot h_n &\qquad \qquad
    Z_{ij}^{(0,0),(r,s)} & = {B_{ij}^{(r)}}^{(0,0),(r,s)} \, \hat{I}_{(r)} \\
    Z_{ij}^{(k,l),(0,0)} & =
    A_{ij}^{(k,l),(0,0)} \cdot h_n &\qquad \qquad
    Z_{ij}^{(k,l),(r,s)} & = {B_{ij}^{(0)}}^{(k,l),(r,s)} \, \sqrt{h_n}
\end{alignat*}
for $1 \leq k,l,r,s \leq m$. Further, we define
${B_{ij}^{(r)}}^{(0,0),(r,s)}=0$ in the case of $r \neq s$,
${B_{ij}^{(0)}}^{(k,l),(r,s)}=0$ in the case of $l \neq r$ or $l
\neq s$ and $A_{ij}^{(k,l),(0,0)}=0$ in the case of $k \neq l$ for
$i,j = 1, \ldots, s$. Here, $\hat{I}_{(k)}$, $1 \leq k \leq m$,
are independent random variables defined by $P(\hat{I}_{(k)} = \pm
\sqrt{3 h_n} ) = \tfrac{1}{6}$ and $P(\hat{I}_{(k)} = 0) =
\tfrac{2}{3}$. The $\hat{I}_{(k,l)}$, $1 \leq k,l \leq m$, are
defined by $\hat{I}_{(k,l)} = \tfrac{1}{2} (\hat{I}_{(k)} \,
\hat{I}_{(l)} + V_{k,l} )$ with independent random variables
$V_{k,l}$ such that $P(V_{k,l} = \pm h_n)=\tfrac{1}{2}$ for $l=1,
\ldots, k-1$, $V_{k,k}=-h_n$ and $V_{l,k}=-V_{k,l}$ for $l=k+1,
\ldots,m$ and $k=1, \ldots, m$.
Thus, we can characterize the SRK method~(\ref{St-srk-method-m})
by the following Butcher array for $1 \leq k,l \leq m$ with $k
\neq l$:
\renewcommand{\arraystretch}{1.8}
\begin{equation*}
\begin{tabular}{c|c|c|c}
    $c^{(0,0)}$ & ${A}^{(0,0),(0,0)}$ & ${B^{(k)}}^{(0,0),(k,k)}$ & \\
    \cline{1-4}
    $c^{(k,k)}$ & ${A}^{(k,k),(0,0)}$ & ${B^{(0)}}^{(k,k),(k,k)}$ & ${B^{(0)}}^{(k,l),(l,l)}$ \\
    \hline
    & $\alpha^T$ & ${{\gamma^{(k)}}^{(k,k)}}^T$ & ${{\gamma^{(k,k)}}^{(k,k)}}^T$ \\
    \cline{2-4}
    & & ${{\gamma^{(k)}}^{(k,l)}}^T$ & ${{\gamma^{(k,l)}}^{(k,l)}}^T$
\end{tabular}
\end{equation*}
%
%
%
\begin{table}[htbp]
\begin{center}
\renewcommand{\arraystretch}{1.3}
\begin{tabular}{r|ccc|ccc|ccc}
    $0$ & & & & & & \\
    $\frac{2}{3}$ & $\frac{2}{3}$ & &  & $1$ & &  & \\
    $\frac{2}{3}$ & $-\frac{1}{3}$ & $1$ & & $0$ & $0$ &  & & \\
    \cline{1-10}
    $0$ & & & & & & \\
    $1$ & $1$ & &  & $1$ & &  & $1$ \\
    $1$ & $1$ & $0$ & & $-1$ & $0$ &  & $-1$ & $0$ \\
    \hline
    & $\frac{1}{4}$ & $\frac{1}{2}$ & $\frac{1}{4}$ &
    $\frac{1}{2}$ & $\frac{1}{4}$ & $\frac{1}{4}$ & $0$ & $\frac{1}{2}$ & $-\frac{1}{2}$\\
    \cline{2-10}
    & & & & $-\frac{1}{2}$ & $\frac{1}{4}$ & $\frac{1}{4}$ & $0$ &
    $\frac{1}{2}$ & $-\frac{1}{2}$
\end{tabular}
\caption{SRK scheme RI1WM of order $p=2.0$ for It{\^o} SDEs.}
\label{SRK-scheme-RI1WM}
\end{center}
\end{table}
The coefficients of the order $2.0$ SRK scheme RI1WM are given in
Table~\ref{SRK-scheme-RI1WM}. For detailed calculations of the
order conditions and the corresponding coefficients we refer to
\cite{Roe03}. \\ \\
As an example for a SRK scheme due to R\"o{\ss}ler applicable to
the Stratonovich SDE~(\ref{St-lg-sde-strato-1}) with $d \geq 1$
and $m \geq 1$ fulfilling a commutativity condition (see
\cite{Roe04a}, \cite{Roe03} for details) we choose now
$\mathcal{M} = \{ j : 1 \leq j \leq m\}$ and
\begin{alignat*}{5}
    z_i^{(0,0)} & = \alpha_i \cdot h_n &\qquad \qquad
    z_i^{(k,k)} & = {\gamma_i^{(k)}}^{(k,k)} \, \hat{I}_{(k)}(h_n) \\
    Z_{ij}^{(0,0),(0,0)} & =
    A_{ij}^{(0,0),(0,0)} \cdot h_n &\qquad \qquad
    Z_{ij}^{(0,0),(k,k)} & = {B_{ij}^{(k)}}^{(0,0),(k,k)} \, \hat{I}_{(k)}(h_n)
    \\
    Z_{ij}^{(k,k),(0,0)} & =
    A_{ij}^{(k,k),(0,0)} \cdot h_n &\qquad \qquad
    Z_{ij}^{(k,k),(l,l)} & = {B_{ij}^{(l)}}^{(k,k),(l,l)} \, \hat{I}_{(l)}(h_n)
\end{alignat*}
for $k,l \in \mathcal{M}$ and $i,j = 1, \ldots, s$. The
coefficients of such a method can be represented by the Butcher
array taking for $k \neq l$ the form
\renewcommand{\arraystretch}{1.8}
\begin{equation*}
\begin{tabular}{c|c|c|c}
    $c^{(0,0)}$ & ${A}^{(0,0),(0,0)}$ & ${B^{(k)}}^{(0,0),(k,k)}$ & \\
    \cline{1-4}
    $c^{(k,k)}$ & ${A}^{(k,k),(0,0)}$ & ${B^{(k)}}^{(k,k),(k,k)}$ &${B^{(l)}}^{(k,k),(l,l)}$ \\
    \hline
    & $\alpha^T$ & ${{\gamma^{(k)}}^{(k,k)}}^T$ &
\end{tabular}
\end{equation*}
For detailed calculations of the order conditions we refer to
\cite{Roe04a} and \cite{Roe03}. The coefficients of the order 2.0
SRK scheme RS1WM are presented in Table~\ref{SRK-scheme-RS1WM}.
%
%
\begin{table}[htbp]
\begin{center}
\renewcommand{\arraystretch}{1.3}
\begin{tabular}{r|cccc|cccc|cccc}
    $0$ & & & &  & & & &  &&&&\\
    $0$ & $0$ & & &  & $0$ & & &  & \\
    $1$ & $1$ & $0$ & & & $-\frac{3}{4}$ & $\frac{3}{4}$ & &  &&&& \\
    $0$ & $0$ & $0$ & $0$ & & $1$ & $0$ & $0$ &   &&&& \\
    \cline{1-13}
    $0$ & & & &  & & & &  &&&& \\
    $0$ & $0$ & & &  & $\frac{2}{3}$ & & &   & $0$ &&& \\
    $1$ & $1$ & $0$ & & & $\frac{1}{12}$ & $\frac{1}{4}$ & &   & $\frac{1}{4}$ & $\frac{3}{4}$ && \\
    $1$ & $1$ & $0$ & $0$ & & $-\frac{5}{4}$ & $\frac{1}{4}$ & $2$ &   & $\frac{1}{4}$ & $\frac{3}{4}$ & $0$ & \\
    \hline
    & $0$ & $0$ & $\frac{1}{2}$ & $\frac{1}{2}$ &
    $\frac{1}{8}$ & $\frac{3}{8}$ & $\frac{3}{8}$ & $\frac{1}{8}$
    &
\end{tabular}
\caption{SRK scheme RS1WM of order $p=2.0$ for Stratonovich SDEs.}
\label{SRK-scheme-RS1WM}
\end{center}
\end{table}
%
%
%
%
%
\section{Stochastic Rooted Tree Theory}
\label{Sec:Rooted-Tree-Theory}
\begin{figure}[H]
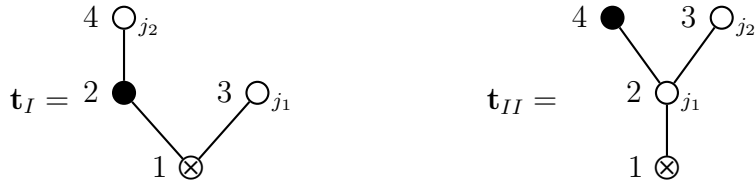

\begin{center}
\begin{tabular}{ccc}
    $\textbf{t}_I = \begin{array}{c}
    \text{\pstree[treemode=U, dotstyle=otimes, dotsize=3.2mm, levelsep=0.1cm, radius=1.6mm, treefit=loose]
    {\Tn}{
    \pstree[treemode=U, dotstyle=otimes, dotsize=3.2mm, levelsep=1cm, radius=1.6mm, treefit=loose]
    {\Tdot~[tnpos=l]{1 }}
    {\pstree{\TC*~[tnpos=l]{2}}{\TC~[tnpos=l]{4}~[tnpos=r]{$\!\!_{j_2}$}}
    \TC~[tnpos=l]{3}~[tnpos=r]{$\!\!_{j_1}$}}
    }}
    \end{array}$
    & \qquad \qquad \qquad &
    $\textbf{t}_{II} = \begin{array}{c}
    \text{\pstree[treemode=U, dotstyle=otimes, dotsize=3.2mm, levelsep=0.1cm, radius=1.6mm, treefit=loose]
    {\Tn}{
    \pstree[treemode=U, dotstyle=otimes, dotsize=3.2mm, levelsep=1cm, radius=1.6mm, treefit=loose]
    {\Tdot~[tnpos=l]{1 }} {\pstree{\TC~[tnpos=l]{2}~[tnpos=r]{$\!\!_{j_1}$}}{\TC*~[tnpos=l]{4}
    \TC~[tnpos=l]{3}~[tnpos=r]{$\!\!_{j_2}$}}}
    }}
    \end{array}$
\end{tabular}
\caption{Two elements of $LTS$ with $j_1, j_2 \in \{1, \ldots,
m\}$.} \label{St-S-tree-examples-tI+tII:Wm}
\end{center}
\end{figure}

\begin{figure}[htb]
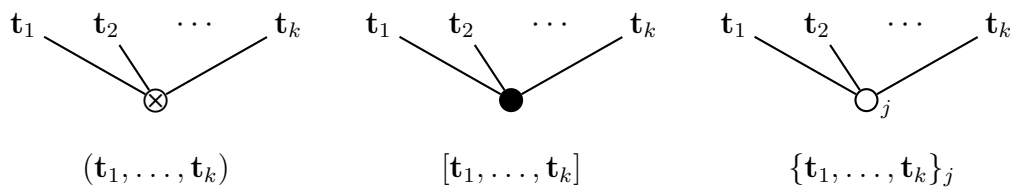

\begin{center}
\begin{tabular}{ccccc}
    \pstree[treemode=U, dotstyle=otimes, dotsize=3.2mm, levelsep=0.1cm, radius=1.6mm, treefit=loose]
    {\Tn}{
    \pstree[treemode=U, dotstyle=otimes, dotsize=3.2mm, levelsep=1cm, radius=1.6mm, treefit=loose, nodesepB=1mm]
    {\Tdot} {\Tr{$\textbf{t}_1$} \Tr{$\textbf{t}_2$} \Tr[edge=none]{$\cdots$} \Tr{$\textbf{t}_k$}}
    }
    & &
    \pstree[treemode=U, dotstyle=otimes, dotsize=3.2mm, levelsep=0.1cm, radius=1.6mm, treefit=loose]
    {\Tn}{
    \pstree[treemode=U, dotstyle=otimes, dotsize=3.2mm, levelsep=1cm, radius=1.6mm, treefit=loose, nodesepB=1mm]
    {\TC*} {\Tr{$\textbf{t}_1$} \Tr{$\textbf{t}_2$} \Tr[edge=none]{$\cdots$} \Tr{$\textbf{t}_k$}}
    }
    & &
    \pstree[treemode=U, dotstyle=otimes, dotsize=3.2mm, levelsep=0.1cm, radius=1.6mm, treefit=loose]
    {\Tn}{
    \pstree[treemode=U, dotstyle=otimes, dotsize=3.2mm, levelsep=1cm, radius=1.6mm, treefit=loose, nodesepB=1mm]
    {\TC~[tnpos=r]{$\!\!_{j}$}} {\Tr{$\textbf{t}_1$} \Tr{$\textbf{t}_2$} \Tr[edge=none]{$\cdots$} \Tr{$\textbf{t}_k$}}
    }
    \\
    $(\textbf{t}_1, \ldots, \textbf{t}_k)$ & &
    $[\textbf{t}_1, \ldots, \textbf{t}_k]$ & & $\,\, \{\textbf{t}_1, \ldots, \textbf{t}_k\}_j$
\end{tabular}
\caption{Writing a colored S-tree with brackets.}
\label{St-tree-bracket-together}
\end{center}
\end{figure}

\begin{figure}[H]
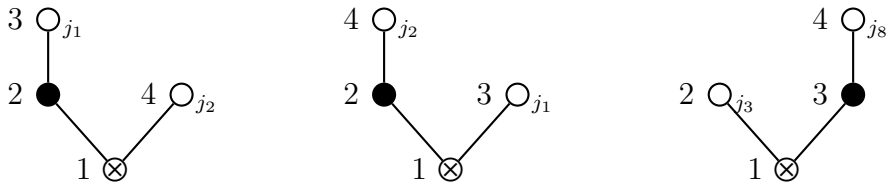

\begin{center}
    \begin{tabular}{ccccc}
    \pstree[treemode=U, dotstyle=otimes, dotsize=3.2mm, levelsep=0.1cm, radius=1.6mm, treefit=loose]
    {\Tn}{
    \pstree[treemode=U, dotstyle=otimes, dotsize=3.2mm, levelsep=1cm, radius=1.6mm, treefit=loose]
    {\Tdot~[tnpos=l]{1 }} {\pstree{\TC*~[tnpos=l]{2}}{\TC~[tnpos=l]{3}~[tnpos=r]{$\!\!_{j_1}$}}
    \TC~[tnpos=l]{4}~[tnpos=r]{$\!\!_{j_2}$}}
    }
    & \qquad \qquad &
    \pstree[treemode=U, dotstyle=otimes, dotsize=3.2mm, levelsep=0.1cm, radius=1.6mm, treefit=loose]
    {\Tn}{
    \pstree[treemode=U, dotstyle=otimes, dotsize=3.2mm, levelsep=1cm, radius=1.6mm, treefit=loose]
    {\Tdot~[tnpos=l]{1 }} {\pstree{\TC*~[tnpos=l]{2}}{\TC~[tnpos=l]{4}~[tnpos=r]{$\!\!_{j_2}$}}
    \TC~[tnpos=l]{3}~[tnpos=r]{$\!\!_{j_1}$}}
    }
    & \qquad \qquad &
    \pstree[treemode=U, dotstyle=otimes, dotsize=3.2mm, levelsep=0.1cm, radius=1.6mm, treefit=loose]
    {\Tn}{
    \pstree[treemode=U, dotstyle=otimes, dotsize=3.2mm, levelsep=1cm, radius=1.6mm, treefit=loose]
    {\Tdot~[tnpos=l]{1 }} {\TC~[tnpos=l]{2}~[tnpos=r]{$\!\!_{j_3}$} \pstree{\TC*~[tnpos=l]{3}}{
    \TC~[tnpos=l]{4}~[tnpos=r]{$\!\!_{j_8}$}}}
    }
    \end{tabular}
\caption{Trees of the same equivalence class.}
\label{St-equal-trees:Wm}
\end{center}
\end{figure}

%
The SDE system~(\ref{Intro-Ito-St-SDE1-integralform-Wm}) can be
represented by an autonomous SDE system
\begin{equation} \label{Ito-St-SDE1-autonom-Wm}
    X_t = x_0 + \int_{t_0}^t a(X_s) \, ds + \sum_{j=1}^m
    \int_{t_0}^t b^j(X_s) \, * dW_s^j
\end{equation}
with one additional equation representing time. Hence, it is
sufficient to treat autonomous SDE systems in the following.
First of all, we give a definition of colored graphs which will be
suitable in the rooted tree theory for SDEs w.r.t.\ a
multi-dimensional Wiener process (see \cite{Roe04b}).
\begin{Def} \label{Def:rooted-S-trees:Wm}
    Let $l$ be a positive integer.
    \begin{enumerate}
        \item A monotonically labelled {\emph{S-tree (stochastic
        tree)}} $\textbf{t}$
        with $l=l(\textbf{t})$ nodes is a pair of maps $\textbf{t}=(\textbf{t}',\textbf{t}'')$
        with
        \begin{equation*}
            \begin{split}
                \textbf{t}' & : \{2, \ldots, l\} \to \{1, \ldots,
                l-1\} \\
                \textbf{t}'' & : \{1, \ldots, l\} \to \mathcal{A} \\
            \end{split}
        \end{equation*}
        so that $\textbf{t}'(i) < i$ for $i=2, \ldots, l$. Unless otherwise noted, we choose
        the set $\mathcal{A} = \{ \gamma, \tau,
        \sigma_{j_k}, k \in \mathbb{N} \}$ where $j_k$ is a variable index with $j_k \in \{1, \ldots, m\}$.
        \item $LTS$ denotes the set of all monotonically labelled
        S-trees w.r.t.\ $\mathcal{A}$.
        Here two trees $\textbf{t}=(\textbf{t}',\textbf{t}'')$ and $\textbf{u}=(\textbf{u}',\textbf{u}'')$
        just differing by their colors $\textbf{t}''$ and $\textbf{u}''$ are
        considered to be identical if there exists a bijective
        map $\pi : \mathcal{A} \to \mathcal{A}$
        with $\pi(\gamma)=\gamma$ and $\pi(\tau)=\tau$ so that
        $\textbf{t}''(i) = \pi(\textbf{u}''(i))$ holds for $i=1, \ldots,
        l$.
    \end{enumerate}
\end{Def}
So $\textbf{t}'$ defines a father son relation between the nodes,
i.e.\ $\textbf{t}'(i)$ is the father of the son $i$. Furthermore
the color $\textbf{t}''(i)$, which consists of one element of the
set $\mathcal{A}$, is added to the node $i$ for $i=1, \ldots,
l(\textbf{t})$. Here, $\tau = $~\pstree[treemode=U,
dotstyle=otimes, dotsize=3.2mm, levelsep=0.1cm, radius=1.6mm,
treefit=loose]
    {\Tn}{
    \pstree[treemode=U, dotstyle=otimes, dotsize=3.2mm, levelsep=0cm, radius=1.6mm, treefit=loose]
    {\TC*~[tnpos=r]{}} {}
    }
is a deterministic node and $\sigma_{j_k} = $~\pstree[treemode=U,
dotstyle=otimes, dotsize=3.2mm, levelsep=0.1cm, radius=1.6mm,
treefit=loose]
    {\Tn}{
    \pstree[treemode=U, dotstyle=otimes, dotsize=3.2mm, levelsep=0cm, radius=1.6mm, treefit=loose]
    {\TC~[tnpos=r]{$\!\!{_{j_k}}$}} {}
    }
is a stochastic node with a variable index $j_k \in \{1, \ldots,
m\}$. In the case of $\gamma \in \mathcal{A}$ the node of type
$\gamma = \,\,$~\pstree[treemode=U, dotstyle=otimes,
dotsize=3.2mm, levelsep=0.1cm, radius=1.6mm, treefit=loose]
    {\Tn}{
    \pstree[treemode=U, dotstyle=otimes, dotsize=3.2mm, levelsep=0cm, radius=1.6mm, treefit=loose]
    {\Tdot~[tnpos=r]{ }} {}
    }
is denoted as the {\emph{root}} and always sketched as the lowest
node of the graph. However, in the case of $\mathcal{A}=\{ \tau,
\sigma_{j_k}, k \in \mathbb{N} \}$, the nodes $\tau$ and
$\sigma_{j_k}$ may also serve as the root of the tree. The
variable index $j_k$ is associated with the $j_k$th component of
the corresponding $m$-dimensional Wiener process of the considered
SDE. In case of a one-dimensional Wiener process one can omit the
variable indices since we have $j_k = 1$ for all $k \in
\mathbb{N}$ (see also \cite{Roe03}). As an example
Figure~\ref{St-S-tree-examples-tI+tII:Wm} presents two elements of
$LTS$.
\noindent For the labelled S-tree $\textbf{t}_I$ in
Figure~\ref{St-S-tree-examples-tI+tII:Wm} we have
$\textbf{t}_I'(2)=\textbf{t}_I'(3)=1$ and $\textbf{t}_I'(4)=2$.
The color of the nodes is given by $\textbf{t}_I''(1)=\gamma$,
$\textbf{t}_I''(2)=\tau$, $\textbf{t}_I''(3)=\sigma_{j_1}$ and
$\textbf{t}_I''(4)=\sigma_{j_2}$.
\begin{Def} \label{Def:order-S-tree-W1}
    Let $\textbf{t}=(\textbf{t}',\textbf{t}'') \in LTS$.
    We denote by $d(\textbf{t}) = \sharp \{ i : \textbf{t}''(i) = \tau \}$
    the number of deterministic nodes, by
    $s(\textbf{t}) = \sharp \{ i : \textbf{t}''(i) = \sigma_{j_k}, k \in
    \mathbb{N} \}$ the number of stochastic nodes and by
    $n(\textbf{t}) = \sharp \{ i : \textbf{t}''(i) = \textbf{t}''(i+1) = \sigma_{j_k},
    k \in \mathbb{N} \}$ the number of pairs of stochastic nodes with the same variable index.
    The {\emph{order}} $\rho(\textbf{t})$ of the tree $\textbf{t}$ is defined as
    $\rho(\textbf{t}) = d(\textbf{t}) + \tfrac{1}{2} s(\textbf{t})$ with $\rho(\gamma) = 0$.
\end{Def}
The order of the trees $\textbf{t}_I$ and $\textbf{t}_{II}$
presented in Figure~\ref{St-S-tree-examples-tI+tII:Wm} can be
calculated as $\rho(\textbf{t}_I)=\rho(\textbf{t}_{II})=2$. Every
labelled S-tree can be written as a combination of three different
brackets defined as follows.
\begin{Def}
    If $\textbf{t}_1, \ldots, \textbf{t}_k$ are colored trees then we denote by
    \begin{equation*}
        (\textbf{t}_1, \ldots, \textbf{t}_k), \,\,\,\,\,
        [\textbf{t}_1, \ldots, \textbf{t}_k] \,\,\,\,\, \text{ and }
        \,\,\,\,\, \{\textbf{t}_1, \ldots, \textbf{t}_k \}_j
    \end{equation*}
    the tree in which $\textbf{t}_1, \ldots, \textbf{t}_k$ are each joined by a
    single branch to $\,\,$
    \pstree[treemode=U, dotstyle=otimes, dotsize=3.2mm, levelsep=0.1cm, radius=1.6mm, treefit=loose]
    {\Tn}{
    \pstree[treemode=U, dotstyle=otimes, dotsize=3.2mm, levelsep=0cm, radius=1.6mm, treefit=loose]
    {\Tdot} {}
    }
    $\,$,
    \pstree[treemode=U, dotstyle=otimes, dotsize=3.2mm, levelsep=0.1cm, radius=1.6mm, treefit=loose]
    {\Tn}{
    \pstree[treemode=U, dotstyle=otimes, dotsize=3.2mm, levelsep=0cm, radius=1.6mm, treefit=loose]
    {\TC*} {}
    }
    $\,\,$and
    \pstree[treemode=U, dotstyle=otimes, dotsize=3.2mm, levelsep=0.1cm, radius=1.6mm, treefit=loose]
    {\Tn}{
    \pstree[treemode=U, dotstyle=otimes, dotsize=3.2mm, levelsep=0cm, radius=1.6mm, treefit=loose]
    {\TC~[tnpos=r]{$\!\!_j$}} {}
    },
    respectively (see Figure~\ref{St-tree-bracket-together}).
\end{Def}
Therefore proceeding recursively, for the two examples
$\textbf{t}_I$ and $\textbf{t}_{II}$ in
Figure~\ref{St-S-tree-examples-tI+tII:Wm} we obtain
    $\textbf{t}_I = ([
    \text{\pstree[treemode=U, dotstyle=otimes, dotsize=3.2mm, levelsep=0.1cm, radius=1.6mm, treefit=loose]
    {\Tn}{
    \pstree[treemode=U, dotstyle=otimes, dotsize=3.2mm, levelsep=0cm, radius=1.6mm, treefit=loose]
    {\TC~[tnpos=r]{$\!\!_{j_2}$}} {}
    }} ],
    \text{\pstree[treemode=U, dotstyle=otimes, dotsize=3.2mm, levelsep=0.1cm, radius=1.6mm, treefit=loose]
    {\Tn}{
    \pstree[treemode=U, dotstyle=otimes, dotsize=3.2mm, levelsep=0cm, radius=1.6mm, treefit=loose]
    {\TC~[tnpos=r]{$\!\!_{j_1}$}} {}
    }} ) = ([\sigma_{j_2}] , \sigma_{j_1})$
    and
    $\textbf{t}_{II} = (\{
    \text{\pstree[treemode=U, dotstyle=otimes, dotsize=3.2mm, levelsep=0.1cm, radius=1.6mm, treefit=loose]
    {\Tn}{
    \pstree[treemode=U, dotstyle=otimes, dotsize=3.2mm, levelsep=0cm, radius=1.6mm, treefit=loose]
    {\TC*} {}
    }},
    \text{\pstree[treemode=U, dotstyle=otimes, dotsize=3.2mm, levelsep=0.1cm, radius=1.6mm, treefit=loose]
    {\Tn}{
    \pstree[treemode=U, dotstyle=otimes, dotsize=3.2mm, levelsep=0cm, radius=1.6mm, treefit=loose]
    {\TC~[tnpos=r]{$\!\!_{j_2}$}} {}
    }}
    \}_{j_1}) = (\{ \tau, \sigma_{j_2} \}_{j_1})$. \\ \\
%
%
Due to the fact that we are interested in calculating weak
approximations, it will turn out that we can concentrate our
considerations to one representative tree of each equivalence
class.
\begin{Def} \label{St-tree-equivalence:Wm}
    Let $\textbf{t}=(\textbf{t}',\textbf{t}'')$ and $\textbf{u}=(\textbf{u}',\textbf{u}'')$
    be elements of $LTS$. Then
    the trees $\textbf{t}$ and $\textbf{u}$ are equivalent, i.e.\ $\textbf{t} \sim \textbf{u}$, if the
    following hold:
    \begin{enumerate}[(i)]
        \item $l(\textbf{t})=l(\textbf{u})$
        \item There exist two bijective maps
            \begin{equation*}
                \begin{split}
                \psi &: \{1, \ldots, l(\textbf{t})\} \rightarrow \{1,
                \ldots, l(\textbf{t})\} \quad \text{ with } \quad
                \psi(1)=1, \\
                \pi &: \mathcal{A} \rightarrow \mathcal{A} \quad \text{ with } \quad
                \pi(\gamma)=\gamma \quad \text{ and } \quad
                \pi(\tau)=\tau,
                \end{split}
            \end{equation*}
            so that the following diagram commutes
            \begin{equation*}
                \begin{psmatrix}[colsep=3cm, rowsep=.5cm]
                    \{2,\ldots,l(\textbf{t})\} & \{1,\ldots,l(\textbf{t})\} & \\
                    & & \mathcal{A} \\
                    \{2,\ldots,l(\textbf{t})\} & \{1,\ldots,l(\textbf{t})\} &
                    \psset{arrows=->, nodesep=3pt}
                    \everypsbox{\scriptstyle}
                    \ncline{1,1}{1,2}\taput{\textbf{t}'}
                    \ncline{3,1}{3,2}\taput{\textbf{u}'}
                    \ncline{1,1}{3,1}\trput{\psi}
                    \ncline{1,2}{3,2}\trput{\psi}
                    \ncline{1,2}{2,3}\taput{\textbf{t}''}
                    \ncline{3,2}{2,3}\tbput{\pi(\textbf{u}'')}
                \end{psmatrix}
            \end{equation*}
    \end{enumerate}
    The set of all equivalence classes under the relation $\sim$
    is denoted by $TS = LTS / \sim$. We denote by $\alpha(\textbf{t})$ the cardinality of
    $\textbf{t}$, i.e.\ the number of possibilities of monotonically
    labelling the nodes of $\textbf{t}$ with numbers $1, \ldots, l(\textbf{t})$.
\end{Def}
Thus, a monotonically labelled S-tree $\textbf{u}$ is equivalent
to $\textbf{t}$, if each label $i$ is replaced by $\psi(i)$ and if
each stochastic node $\sigma_{j_k}$ with variable index $j_k$ is
replaced by an other stochastic node $\pi(\sigma_{j_k})$.
Thus, all trees in Figure~\ref{St-equal-trees:Wm} belong to the
same equivalence class as $\textbf{t}_I$ in the example above,
since the indices $j_1$ and $j_2$ are just renamed either by $j_2$
and $j_1$ or $j_8$ and $j_3$, respectively. Finally the graphs
differ only in the labelling of their number indices.
%
%
For every rooted tree $\textbf{t} \in LTS$, there exists a
corresponding {\emph{elementary differential}} which is a direct
generalization of the differential in the deterministic case (see,
e.g., \cite{Butcher87}). For $j \in \{1, \ldots, m\}$, the
elementary differential is defined recursively by
\begin{equation*}
    F(\gamma)(x) = f(x), \qquad
    F(\tau)(x) = a(x), \qquad
    F(\sigma_j)(x) = b^j(x),
\end{equation*}
for single nodes and by
\begin{equation} \label{St-elementary-differential-F:Wm}
    F(\textbf{t})(x) =
    \begin{cases}
    f^{(k)}(x) \cdot (F(\textbf{t}_1)(x), \ldots, F(\textbf{t}_k)(x)) &
    \text{for } \textbf{t}=(\textbf{t}_1, \ldots, \textbf{t}_k) \\
    a^{(k)}(x) \cdot (F(\textbf{t}_1)(x), \ldots,
    F(\textbf{t}_k)(x)) & \text{for } \textbf{t}=[\textbf{t}_1, \ldots, \textbf{t}_k] \\
    {b^j}^{(k)}(x) \cdot (F(\textbf{t}_1)(x),
    \ldots, F(\textbf{t}_k)(x)) & \text{for } \textbf{t}=\{\textbf{t}_1, \ldots,
    \textbf{t}_k\}_j
    \end{cases}
\end{equation}
for a tree $\textbf{t}$ with more than one node. Here $f^{(k)}$,
$a^{(k)}$ and ${b^j}^{(k)}$ define a symmetric $k$-linear
differential operator, and one can choose the sequence of labelled
S-trees $\textbf{t}_1, \ldots, \textbf{t}_k$ in an arbitrary
order. For example, the $I$th component of $a^{(k)} \cdot
(F(\textbf{t}_1), \ldots, F(\textbf{t}_k))$ can be written as
\begin{align*}
    ( a^{(k)} \cdot (F(\textbf{t}_1), \ldots, F(\textbf{t}_k)) )^I
    &= \sum_{J_1, \ldots, J_k=1}^d \frac{\partial^k
    a^I}{\partial x^{J_1} \ldots \partial x^{J_k}} \,
    (F^{J_1}(\textbf{t}_1), \ldots, F^{J_k}(\textbf{t}_k))
\end{align*}
where the components of vectors are denoted by superscript
indices, which are chosen as capitals.
As a result of this we get for $\textbf{t}_I$ and
$\textbf{t}_{II}$ the elementary differentials
\begin{equation*}
    \begin{split}
    &F(\textbf{t}_I) = f'' (a'(b^{j_2}), b^{j_1}) = \sum_{J_1,J_2=1}^d
    \frac{\partial^2 f}{\partial x^{J_1} \partial x^{J_2}}
    \big( \sum_{K_1=1}^d \frac{\partial a^{J_1}}{\partial x^{K_1}}
    \, b^{K_1,j_2} \cdot b^{J_2,j_1} \big) \\
    &F(\textbf{t}_{II}) = f' ({b^{j_1}}'' (a, b^{j_2})) = \sum_{J_1=1}^d
    \frac{\partial f}{\partial x^{J_1}} \big( \sum_{K_1, K_2 =1}^d
    \frac{\partial^2 b^{J_1,j_1}}{\partial x^{K_1} \partial
    x^{K_2}} \, a^{K_1} \cdot b^{K_2,j_2} \big)
    \end{split}
\end{equation*}
It has to be pointed out that the elementary differentials for the
trees presented in Figure~\ref{St-equal-trees:Wm} coincide with
$F(\textbf{t}_I)$ if the variable indices $j_i$ are simply renamed
by a suitable bijective mapping $\pi$.
%
%
%
%
%
\section{Taylor Expansion for It{\^o} and Stratonovich SDEs}
\label{Sec:Taylor-Expansion-SDE}
%
For the expansion of the expectation of some functional applied to
the solution $(X_t)_{t \in \mathcal{I}}$ of the $d$-dimensional
SDE~(\ref{Ito-St-SDE1-autonom-Wm}) considered either w.r.t.\
It{\^o} or Stratonovich calculus, some subsets $LTS(I)$ and
$LTS(S)$ of $LTS$ have to be introduced, respectively.
\begin{Def} \label{Ito-Stratonovich-subset-LTS(*)}
    For $* \in \{I,S \}$ let $LTS(*)$ denote the set of trees
    $\textbf{t} 
    \in LTS$ having a
    root $\gamma = \,\,$~\pstree[treemode=U, dotstyle=otimes, dotsize=3.2mm, levelsep=0.1cm, radius=1.6mm, treefit=loose]
    {\Tn}{
    \pstree[treemode=U, dotstyle=otimes, dotsize=3.2mm, levelsep=0cm, radius=1.6mm, treefit=loose]
    {\Tdot} {}
    }$\,\,\,$
    and which can be constructed by a finite number
    of steps of the form
    \begin{enumerate}[a)]
        \item adding a deterministic node $\tau = $
        \pstree[treemode=U, dotstyle=otimes, dotsize=3.2mm, levelsep=0.1cm, radius=1.6mm, treefit=loose]
        {\Tn}{
        \pstree[treemode=U, dotstyle=otimes, dotsize=3.2mm, levelsep=0cm, radius=1.6mm, treefit=loose]
        {\TC*} {}
        }, or
        \item adding two stochastic nodes $\sigma_{j_k} = $
        \pstree[treemode=U, dotstyle=otimes, dotsize=3.2mm, levelsep=0.1cm, radius=1.6mm, treefit=loose]
        {\Tn}{
        \pstree[treemode=U, dotstyle=otimes, dotsize=3.2mm, levelsep=0cm, radius=1.6mm, treefit=loose]
        {\TC~[tnpos=r]{$\!\!_{j_k}$}} {}
        }, where both nodes get the same new variable index $j_k$ for some $k \in \mathbb{N}$.
        Additionally, in the case of $*=I$
        neither of the two nodes is allowed to be the father of the other.
    \end{enumerate}
    The nodes have to be labelled in the same order as they have been
    added by the construction of the tree. Further $TS(*) = LTS(*)/ \sim$
    denotes the equivalence class under the relation of
    Definition~\ref{St-tree-equivalence:Wm} restricted to $LTS(*)$
    and $\alpha_*(\textbf{t})$ denotes the cardinality of $\textbf{t}$ in $LTS(*)$
    for $* \in \{I,S\}$, respectively.
\end{Def}
Since the number of stochastic nodes is always even with
$n(\textbf{t}) = s(\textbf{t}) / 2$, the order $\rho(\textbf{t})$
has to be an integer and $\textbf{t}$ owns the variable indices
$j_1, \ldots, j_{n(\textbf{t})}$. As the construction of the trees
in $LTS(I)$ is more restrictive than of the ones in $LTS(S)$, it
holds $LTS(I) \subset LTS(S)$.
\begin{figure}[H]
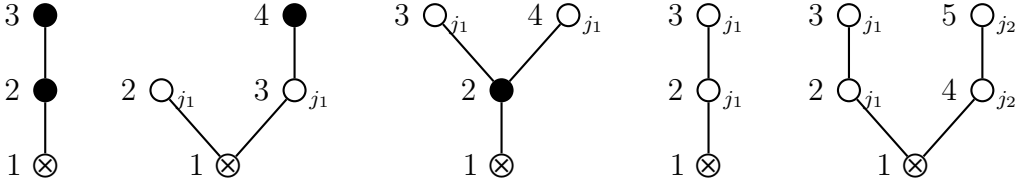

\begin{center}
    \begin{tabular}{ccccccccccc}
    \pstree[treemode=U, dotstyle=otimes, dotsize=3.2mm, levelsep=0.1cm, radius=1.6mm, treefit=loose]
    {\Tn}{
    \pstree[treemode=U, dotstyle=otimes, dotsize=3.2mm, levelsep=1cm, radius=1.6mm, treefit=loose]
    {\Tdot~[tnpos=l]{1 }} {\pstree{\TC*~[tnpos=l]{2}}{\TC*~[tnpos=l]{3}}}
    }
    &
    &
    \pstree[treemode=U, dotstyle=otimes, dotsize=3.2mm, levelsep=0.1cm, radius=1.6mm, treefit=loose]
    {\Tn}{
    \pstree[treemode=U, dotstyle=otimes, dotsize=3.2mm, levelsep=1cm, radius=1.6mm, treefit=loose]
    {\Tdot~[tnpos=l]{1 }} {\TC~[tnpos=l]{2}~[tnpos=r]{$\!\!_{j_1}$}
    \pstree{\TC~[tnpos=l]{3}~[tnpos=r]{$\!\!_{j_1}$}}{\TC*~[tnpos=l]{4}}}
    }
    &
    &
    \pstree[treemode=U, dotstyle=otimes, dotsize=3.2mm, levelsep=0.1cm, radius=1.6mm, treefit=loose]
    {\Tn}{
    \pstree[treemode=U, dotstyle=otimes, dotsize=3.2mm, levelsep=1cm, radius=1.6mm, treefit=loose]
    {\Tdot~[tnpos=l]{1 }} {\pstree{\TC*~[tnpos=l]{2}}{\TC~[tnpos=l]{3}~[tnpos=r]{$\!\!_{j_1}$}
    \TC~[tnpos=l]{4}~[tnpos=r]{$\!\!_{j_1}$}}}
    }
    &
    &
    \pstree[treemode=U, dotstyle=otimes, dotsize=3.2mm, levelsep=0.1cm, radius=1.6mm, treefit=loose]
    {\Tn}{
    \pstree[treemode=U, dotstyle=otimes, dotsize=3.2mm, levelsep=1cm, radius=1.6mm, treefit=loose]
    {\Tdot~[tnpos=l]{1 }} {\pstree{\TC~[tnpos=l]{2}~[tnpos=r]{$\!\!_{j_1}$}}{\TC~[tnpos=l]{3}~[tnpos=r]{$\!\!_{j_1}$}}}
    }
    &
    &
    \pstree[treemode=U, dotstyle=otimes, dotsize=3.2mm, levelsep=0.1cm, radius=1.6mm, treefit=loose]
    {\Tn}{
    \pstree[treemode=U, dotstyle=otimes, dotsize=3.2mm, levelsep=1cm, radius=1.6mm, treefit=loose]
    {\Tdot~[tnpos=l]{1 }} {\pstree{\TC~[tnpos=l]{2}~[tnpos=r]{$\!\!_{j_1}$}}{\TC~[tnpos=l]{3}~[tnpos=r]{$\!\!_{j_1}$}}
    \pstree{\TC~[tnpos=l]{4}~[tnpos=r]{$\!\!_{j_2}$}}{\TC~[tnpos=l]{5}~[tnpos=r]{$\!\!_{j_2}$}}}
    }
    \end{tabular}
\caption{Some trees which belong to $LTS(I)$ or $LTS(S)$.}
\label{Trees-in-LTS(*):Wm}
\end{center}
\end{figure}
All trees of Figure~\ref{Trees-in-LTS(*):Wm} belong to $LTS(S)$,
however only the first three trees belong to $LTS(I)$. For the
last tree, there is a similar tree $(\{ \sigma_{j_2} \}_{j_1}, \{
\sigma_{j_2} \}_{j_1})$ which belongs to $LTS(I)$. The only
difference is the sequence of the construction, i.e.\ the correct
father-son relationship for the stochastic nodes. Clearly, a tree
like $(\{ \tau \}_{j_1} )$ or $(\{ [ \sigma_{j_1} ] \}_{j_1} )$
neither belongs to $LTS(I)$ nor to $LTS(S)$. \\ \\
The following result gives an expansion for the solution process
of an It{\^o} and a Stratonovich SDE, respectively, by the use of
colored rooted trees.
\begin{The} \label{St-tree-expansion-exact-sol:Wm}
    Let $(X_t)_{t \in \mathcal{I}}$ be the solution of the
    stochastic differential equation
    system~(\ref{Ito-St-SDE1-autonom-Wm})
    with initial value $X_{t_0} = x_0 \in \mathbb{R}^d$.
    Then for $p \in \mathbb{N}_0$ and $f, a^i, \tilde{a}^i, b^{i,j}
    \in C_P^{2(p+1)}(\mathbb{R}^d, \mathbb{R})$
    for $i=1, \ldots, d$, $j=1, \ldots, m$ and for $t \in [t_0,T]$
    the following truncated expansion holds:
    \begin{equation}
    \begin{split} \label{Strato-tree-expansion-exact-sol-formula1:Wm}
        E^{t_0,x_0}(f(X_{t})) = &\sum_{\substack{\textbf{t} \in LTS(*) \\
        \rho(\textbf{t}) \leq p}} \,\,
        \sum_{j_1, \ldots, j_{s(\textbf{t})/2}=1}^m
        \frac{F(\textbf{t})(x_0) }{2^{s(\textbf{t})/2} \, \rho(\textbf{t})!} \,
        (t-t_0)^{\rho(\textbf{t})} + O( (t-t_0)^{p+1} ) \\
        = &\sum_{\substack{\textbf{t} \in TS(*) \\ \rho(\textbf{t}) \leq p}} \,\,
        \sum_{j_1, \ldots, j_{s(\textbf{t})/2}=1}^m
        \frac{\alpha_*(\textbf{t}) \, F(\textbf{t})(x_0) }{2^{s(\textbf{t})/2}
        \, \rho(\textbf{t})!} \, (t-t_0)^{\rho(\textbf{t})} + O( (t-t_0)^{p+1}
        )
    \end{split}
    \end{equation}
    Here, $* = I$ for the It{\^o} version of SDE~(\ref{Ito-St-SDE1-autonom-Wm}),
    and $* = S$ for the Stratonovich version of SDE~(\ref{Ito-St-SDE1-autonom-Wm}).
\end{The}
{\bf{Proof.}} For a proof we refer to Theorem~3.2 and Theorem~4.2
together with Proposition~5.1 in~\cite{Roe04b}. \hfill $\Box$
%
%
%
\section{Taylor Expansion for the Stochastic Runge-Kutta method}
\label{Sec:Expansion-SRK-method}
%
In order to derive conditions such that the stochastic Runge-Kutta
method (\ref{St-srk-method-m}) converges in the weak sense with
some specified order, a Taylor expansion of the numerical solution
based on colored rooted trees has to be developed.
We follow the approach of Butcher~\cite{Butcher87} in a similar
way as in Burrage and Burrage~\cite{BuBu96}, \cite{BuBu00a},
Hairer~\cite{Hai81}, Hairer, N{\o}rsett
and Wanner~\cite{HNW93} and R{\"o}{\ss}ler~\cite{Roe03}. \\ \\
For notational convenience, for the set of multi-indices
$\mathcal{M}$ we here put $\overline{\mathcal{M}} = \mathcal{M}
\cup \{0\}$ and we set $\theta_0(h) = h$ and denote by
$\theta(h)=(\theta_0(h), \theta_{\nu_1}(h), \ldots,
\theta_{\nu_{\kappa}}(h))^T$, $\nu_i \in \mathcal{M}$, the
corresponding $\kappa+1$-dimensional vector of random
variables\footnote{Then $Y(t)=A(t_0, Y(t_0), \theta(t-t_0))$ and
$Y(t_0)=A(t_0, Y(t_0),0, \ldots, 0)$.} with
$\kappa=|\mathcal{M}|$. Further, it is assumed that
$\theta_{\nu}(0) = 0$ for all $\nu \in \overline{\mathcal{M}}$.
Due to condition~(\ref{SRK-weights-condition}), it is sufficient
to consider autonomous SRK methods~(\ref{St-srk-method-m}) in the
following. We denote $t_n$ by $t_0$ and for a given $t=t_0 + h$
the approximations $Y_n$ and $Y_{n+1}$ are denoted by $Y(t_0)$ and
$Y(t)$ in (\ref{St-srk-method-m}), respectively. Further, the
values $H_i^{(k,\nu)}$ are denoted by $H_i^{(k,\nu)}(t)$ in order
to stress the dependency on $t$ of the random variables
$\theta_0(t-t_0), \theta_{\nu_1}(t-t_0), \ldots,
\theta_{\nu_{\kappa}}(t-t_0)$ appearing in $H_i^{(k,\nu)}$.
For the Taylor expansion of the SRK method $Y(t)=A(t_0,Y(t_0),
\theta(t-t_0))$ as a function of $\theta_0, \theta_{\nu_1},
\ldots, \theta_{\nu_{\kappa}}$, the differential operator
$\mathcal{D}^k$ for $k \in \mathbb{N}$ is introduced as
\begin{equation} \label{Taylor-Operator-Dk}
    \mathcal{D}^k = \sum_{\nu_1, \ldots, \nu_k \in \overline{\mathcal{M}}}
    \Delta {\theta}_{\nu_1} \cdot \Delta {\theta}_{\nu_2} \cdot
    \ldots \cdot \Delta {\theta}_{\nu_k} \cdot
    \frac{\partial^k}{\partial
    {\theta}_{\nu_1} \partial {\theta}_{\nu_2} \ldots
    \partial {\theta}_{\nu_k}}
\end{equation}
with $\Delta {\theta}_{\nu} = {\theta}_{\nu}(h) -
{\theta}_{\nu}(0)$ and we denote by $\mathcal{D}^0 \equiv \Id$.
Under the assumption that $f$, $a$ and $b^j$, $1 \leq j \leq m$,
are sufficiently differentiable, we apply the Theorem of Taylor
and get for $n \in \mathbb{N}$
\begin{equation} \label{SRK-w1-Taylor-D}
    f(Y(t)) = \sum_{k=0}^{n} \frac{\mathcal{D}^k
    f(Y(t_0))}{k!} + \mathcal{R}_n(t,t_0)
\end{equation}
with a remainder term $\mathcal{R}_n$ which can be written in
Lagrange form as
\begin{equation}
    \mathcal{R}_n(t,t_0) = \frac{\mathcal{D}^{n+1} f(Y(t_0 + \xi \,
    h))}{(n+1)!}
\end{equation}
with some $\xi \in \, ]0,1[$ and $h=t-t_0$. \\ \\
The next step is the computation of $\mathcal{D}^k f(Y(t_0))$ for
$k \in \mathbb{N}_0$, i.e.\ the $k$th derivative of the numerical
solution $f(Y(t))$. Therefore, generalized versions of the
{\emph{Leibniz formula}} and of {\emph{ Fa{\`a} di Bruno's
formula}} (see, e.g., \cite{HNW93}) are helpful. \\ \\
To begin with, a multi-dimensional version of the Leibniz formula
fitted to the expansion of the SRK method is given. Let $q \in
\mathbb{N}$, $k \in \{0,1, \ldots, m\}$ and $\nu \in
\overline{\mathcal{M}}$. Then the formula
\begin{equation} \label{Lem-Leibniz-formula-2-eqn1}
    \begin{split}
    \sum_{\nu_1, \ldots, \nu_q \in \overline{\mathcal{M}}}
    &\frac{\partial^q H_i^{(k, \nu)}(t)^J}{\partial \theta_{\nu_1}
    \ldots \partial \theta_{\nu_q}}
    = \, \, q \cdot  \sum_{j=1}^s
    A_{ij}^{(k,\nu),(0,0)} \sum_{\nu_1, \ldots, \nu_{q-1} \in \overline{\mathcal{M}}}
    \frac{\partial^{q-1}
    a(H_j^{(0,0)}(t))^J}{\partial \theta_{\nu_1}
    \ldots \partial \theta_{\nu_{q-1}}} \\
    + & \,\, q \cdot \sum_{j=1}^s \sum_{r=1}^m
    \sum_{\mu, \nu_q \in \mathcal{M}} {B_{ij}^{(\nu_q)}}^{(k,\nu),(r,\mu)}
    \sum_{\nu_1, \ldots, \nu_{q-1} \in \overline{\mathcal{M}}} \frac{\partial^{q-1}
    b^r (H_j^{(r,\mu)}(t))^J}{\partial \theta_{\nu_1} \ldots \partial
    \theta_{\nu_{q-1}}} \\
    + & \sum_{j=1}^s A_{ij}^{(k,\nu),(0,0)} \, \theta_0(t-t_0)
    \sum_{\nu_1, \ldots, \nu_q \in \overline{\mathcal{M}}} \frac{\partial^q a(H_j^{(0,0)}(t))^J}{\partial
    \theta_{\nu_1} \ldots \theta_{\nu_q}} \\
    + & \sum_{j=1}^s \sum_{r=1}^m
    \sum_{\mu, \iota \in \mathcal{M}}
    {B_{ij}^{(\iota)}}^{(k,\nu),(r,\mu)} \, \theta_{\iota}(t-t_0)
    \sum_{\nu_1, \ldots, \nu_q \in \overline{\mathcal{M}}} \frac{\partial^q
    b^r(H_j^{(r,\mu)}(t))^J}{\partial \theta_{\nu_1} \ldots \theta_{\nu_q}}
    \end{split}
\end{equation}
can be easily calculated (see also \cite{Roe03}, Lemma~2.5.3).
In order to state a generalized version of Fa{\`a} di Bruno's
formula~\cite{HNW93}, we introduce a special set of trees
corresponding to the derivatives of the composition of two
functions.
\begin{figure}[htbp]
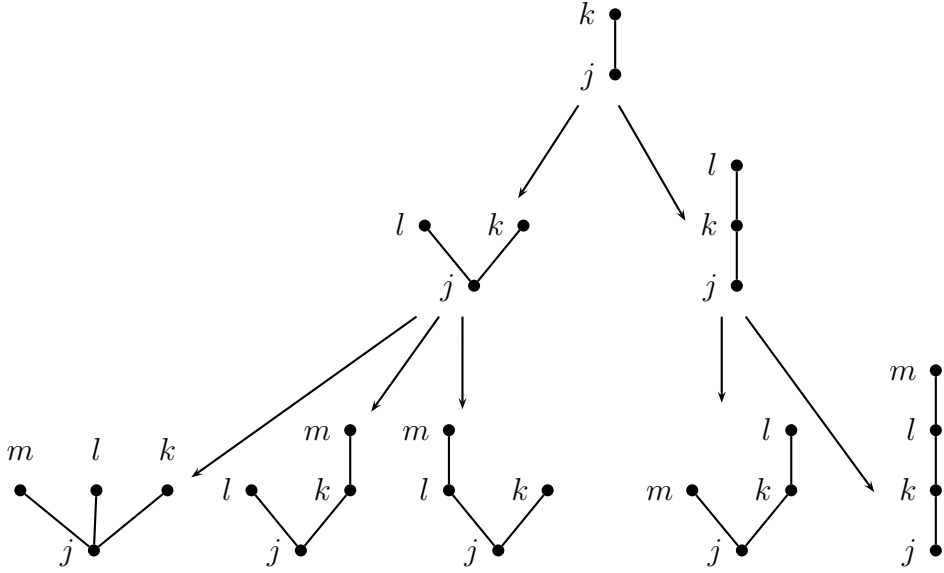

\begin{center}
    \begin{psmatrix}[colsep=0.2cm, rowsep=0.8cm]
    & & &
    \pstree[treemode=U, dotstyle=otimes, dotsize=1.6mm, levelsep=0.1cm, radius=0.8mm, treefit=tight]
    {\Tn}{
    \pstree[treemode=U, dotstyle=otimes, dotsize=1.6mm, levelsep=0.8cm, radius=0.8mm, treefit=tight]
    {\TC*~[tnpos=l]{$j$}} {\TC*~[tnpos=l]{$k$}}
    } \\
    & &
    \pstree[treemode=U, dotstyle=otimes, dotsize=1.6mm, levelsep=0.1cm, radius=0.8mm, treefit=tight]
    {\Tn}{
    \pstree[treemode=U, dotstyle=otimes, dotsize=1.6mm, levelsep=0.8cm, radius=0.8mm, treefit=tight]
    {\TC*~[tnpos=l]{$j$}} {\TC*~[tnpos=l]{$l$} \TC*~[tnpos=l]{$k$}}
    }
    & \qquad \qquad &
    \pstree[treemode=U, dotstyle=otimes, dotsize=1.6mm, levelsep=0.1cm, radius=0.8mm, treefit=tight]
    {\Tn}{
    \pstree[treemode=U, dotstyle=otimes, dotsize=1.6mm, levelsep=0.8cm, radius=0.8mm, treefit=tight]
    {\TC*~[tnpos=l]{$j$}} {\pstree{\TC*~[tnpos=l]{$k$}}{\TC*~[tnpos=l]{$l$}}}
    } \\
    \pstree[treemode=U, dotstyle=otimes, dotsize=1.6mm, levelsep=0.1cm, radius=0.8mm, treefit=tight]
    {\Tn}{
    \pstree[treemode=U, dotstyle=otimes, dotsize=1.6mm, levelsep=0.8cm, radius=0.8mm, treefit=tight]
    {\TC*~[tnpos=l]{$j$}} {\TC*~[tnpos=a]{$m$} \TC*~[tnpos=a]{$l$} \TC*~[tnpos=a]{$k$}}
    }
    & \quad
    \pstree[treemode=U, dotstyle=otimes, dotsize=1.6mm, levelsep=0.1cm, radius=0.8mm, treefit=tight]
    {\Tn}{
    \pstree[treemode=U, dotstyle=otimes, dotsize=1.6mm, levelsep=0.8cm, radius=0.8mm, treefit=tight]
    {\TC*~[tnpos=l]{$j$}} {\TC*~[tnpos=l]{$l$} \pstree{\TC*~[tnpos=l]{$k$}}{\TC*~[tnpos=l]{$m$}}}
    }
    & \quad
    \pstree[treemode=U, dotstyle=otimes, dotsize=1.6mm, levelsep=0.1cm, radius=0.8mm, treefit=tight]
    {\Tn}{
    \pstree[treemode=U, dotstyle=otimes, dotsize=1.6mm, levelsep=0.8cm, radius=0.8mm, treefit=tight]
    {\TC*~[tnpos=l]{$j$}} {\pstree{\TC*~[tnpos=l]{$l$}}{\TC*~[tnpos=l]{$m$}} \TC*~[tnpos=l]{$k$}}
    }
    & &
    \pstree[treemode=U, dotstyle=otimes, dotsize=1.6mm, levelsep=0.1cm, radius=0.8mm, treefit=tight]
    {\Tn}{
    \pstree[treemode=U, dotstyle=otimes, dotsize=1.6mm, levelsep=0.8cm, radius=0.8mm, treefit=tight]
    {\TC*~[tnpos=l]{$j$}} {\TC*~[tnpos=l]{$m$} \pstree{\TC*~[tnpos=l]{$k$}}{\TC*~[tnpos=l]{$l$}}}
    }
    & \qquad \qquad &
    \pstree[treemode=U, dotstyle=otimes, dotsize=1.6mm, levelsep=0.1cm, radius=0.8mm, treefit=tight]
    {\Tn}{
    \pstree[treemode=U, dotstyle=otimes, dotsize=1.6mm, levelsep=0.8cm, radius=0.8mm, treefit=tight]
    {\TC*~[tnpos=l]{$j$}} {\pstree{\TC*~[tnpos=l]{$k$}}{\pstree{\TC*~[tnpos=l]{$l$}}{\TC*~[tnpos=l]{$m$}}}}
    }
    \psset{arrows=->, nodesep=.2cm}
    \ncline{1,4}{2,3}
    \ncline{1,4}{2,5}
    \ncline{2,3}{3,1}
    \ncline{2,3}{3,2}
    \ncline{2,3}{3,3}
    \ncline{2,5}{3,5}
    \ncline{2,5}{3,7}
    \end{psmatrix}
\caption{Some special trees representing the derivatives of
$g(h)^J$.} \label{SRK-expansion-special-trees-Faa-di-Bruno}
\end{center}
\end{figure}
For example, if we consider $g \circ h$, we get for the $J$th
component of the third derivative
\begin{equation}
    \begin{split}
    & \frac{\partial^3 g(h)^J}{\partial x^K \partial x^L \partial
    x^M} = \sum_{K_1,K_2,K_3} g^J_{K_1 K_2 K_3}(h)
    \left( \frac{\partial h^{K_1}}{\partial x^M} \cdot \frac{\partial
    h^{K_2}}{\partial x^L} \cdot \frac{\partial h^{K_3}}{\partial x^K}
    \right) \\
    & \quad + \sum_{K_1,K_2} g_{K_1 K_2}^J(h)
    \left( \frac{\partial h^{K_1}}{\partial x^L}
    \cdot \frac{\partial^2 h^{K_2}}{\partial x^K \partial x^M} \right) +
    \sum_{K_1,K_2} g_{K_1 K_2}^J(h)
    \left( \frac{\partial^2 h^{K_1}}{\partial x^L \partial x^M}
    \cdot \frac{\partial h^{K_2}}{\partial x^K} \right) \\
    & \quad + \sum_{K_1,K_2} g_{K_1 K_2}^J(h)
    \left( \frac{\partial h^{K_1}}{\partial x^M}
    \cdot \frac{\partial^2 h^{K_2}}{\partial x^K \partial x^L} \right) +
    \sum_{K_1} g_{K_1}^J(h)
    \left( \frac{\partial^3 h^{K_1}}{\partial x^K \partial x^L \partial x^M}
    \right) .
    \end{split}
\end{equation}
The corresponding special trees are presented in the last line of
Figure~\ref{SRK-expansion-special-trees-Faa-di-Bruno}. Here the
number $m$ of indices $K_1, \ldots, K_m$ depends on the number of
ramifications of the root. Each time $g(h)^J$ is differentiated,
one has to
\begin{enumerate}[(i)]
    \item differentiate the first factor $g_{K_1 \ldots}^J$, i.e.,
    add a new branch to the root $j$, \\
    \item increase the number of derivatives of each of the $h$
    functions by 1, which is presented by lengthening the corresponding
    branch.
\end{enumerate}
So each time we differentiate, we have to add a new label.
{\emph{All}} trees which are obtained in this way are those
{\emph{special trees}} which have no ramifications except at the
root. \\ \\
In order to take into account colored stochastic trees with their
meaning for the expansion of the SRK method, special trees having
either a root of type $\gamma$, $\tau$ or $\sigma_{j}$ have to be
considered in the following. This is due to the analysis of the
composed functions $f(Y(t))$, $a(H_i^{(0,0)}(t))$ and
$b^j(H_i^{(j,\nu)}(t))$.
\begin{Def}
    The set of {\emph{special labelled
    trees}} with $q$ nodes having no ramifications except at the
    root is denoted by $SLTS_q$. For $\textbf{u} \in SLTS_q$ we denote by
    $m=m(\textbf{u})$ the number of ramifications of the root of $\textbf{u}$.
    Further we denote by $SLTS_q^{(M)} \subset SLTS_q$ with
    $M \subset \mathcal{A}= \{ \gamma, \tau, \sigma_{j_k} : k \in \mathbb{N} \}$
    the set of special labelled trees in
    $SLTS_q$ having a root of type $\pi$ with $\pi \in M$.
\end{Def}
Now a formula similar to Fa{\`a} di Bruno's formula fitted to the
stochastic setting can be stated.
\begin{Lem} \label{Lem-Faa-di-Bruno-1}
    For $q \in \mathbb{N}$,  $\pi \in \mathcal{A}$
    and functions $g: \mathbb{R}^{d}
    \rightarrow \mathbb{R}^{r}$ and $h: \mathbb{R}^{\kappa+1}
    \rightarrow \mathbb{R}^{d}$,
    the multi-dimensional chain rule
    \begin{equation} \label{Lem-Faa-di-Bruno-1-eqn1}
        \begin{split}
        \sum_{\nu_1, \ldots, \nu_q \in \overline{\mathcal{M}}}
        & \frac{\partial^q g(h)^J}{\partial \theta_{\nu_1} \ldots
        \partial \theta_{\nu_q}} = \sum_{\textbf{u} \in SLTS_{q+1}^{(\pi)}}
        \sum_{K_1, \ldots, K_{m(\textbf{u})} = 1}^d g_{K_1 \ldots K_{m(\textbf{u})}}^J (h) \,
        \times \\
        \times & \left( \Big( \sum_{\nu_1, \ldots, \nu_{\delta_1} \in \overline{\mathcal{M}}}
        \frac{\partial^{\delta_1} h^{K_1}}{\partial \theta_{\nu_1}
        \ldots \partial \theta_{\nu_{\delta_1}}} \Big) \cdot \ldots
        \cdot \Big( \sum_{\nu_1, \ldots, \nu_{\delta_{m(\textbf{u})}} \in \overline{\mathcal{M}}}
        \frac{\partial^{\delta_{m(\textbf{u})}} h^{K_{m(\textbf{u})}}}{\partial \theta_{\nu_1}
        \ldots \partial \theta_{\nu_{\delta_{m(\textbf{u})}}}} \Big) \right)
        \end{split}
    \end{equation}
    holds. Here $m=m(\textbf{u})$ denotes the number of ramifications of
    the root of the special tree $\textbf{u}=(\textbf{u}_1, \ldots,
    \textbf{u}_m)_{\pi}$ with a root of type $\pi \in \mathcal{A}$
    and $\delta_i=l(\textbf{u}_i)$ describes
    the number of nodes of the subtree $\textbf{u}_i$ for $1 \leq i \leq m$
    with $\delta_1 + \ldots + \delta_{m(\textbf{u})}=q$ for all $\textbf{u} \in
    SLTS_{q+1}^{(\pi)}$.
\end{Lem}
{\bf Proof.} We prove Lemma~\ref{Lem-Faa-di-Bruno-1} by induction
on $q$. For $q=1$ and $\pi \in \mathcal{A}$ we have
\begin{equation*}
    \sum_{\nu_1 \in \overline{\mathcal{M}}} \frac{\partial g(h)^J}{\partial
    \theta_{\nu_1}} = \sum_{\nu_1 \in \overline{\mathcal{M}}} \sum_{K_1=1}^d
    \frac{\partial g(h)^J}{\partial x^{K_1}} \cdot \frac{\partial
    h^{K_1}}{\partial \theta_{\nu_1}} = \sum_{\textbf{u} \in
    SLTS_2^{(\pi)}} \sum_{K_1=1}^d g_{K_1}^J(h) \left(
    \sum_{\nu_1 \in \overline{\mathcal{M}}} \frac{\partial h^{K_1}}{\partial
    \theta_{\nu_1}} \right)
\end{equation*}
with the set $SLTS_2^{(\pi)} = \{ \, ( \tau )_{\pi} \, \}$,
$m(\textbf{u})=1$ and $\delta_1=1$. Assuming now that the
hypothesis~(\ref{Lem-Faa-di-Bruno-1-eqn1}) holds for $q$, we prove
it for $q+1$. Therefore we write shortly
\begin{equation} \label{Bezeichnung-Faa-di-Bruno-1-short}
    (h^{K})^{(\delta)} = \sum_{\nu_1, \ldots,
    \nu_{\delta} \in \overline{\mathcal{M}}} \frac{\partial^{\delta} h^K}{\partial
    \theta_{\nu_1} \ldots \partial \theta_{\nu_{\delta}}}
\end{equation}
and we thus get
\begin{equation*}
    \begin{split}
    & \sum_{\nu_1, \ldots, \nu_{q+1} \in \overline{\mathcal{M}}}
    \frac{\partial^{q+1}
    g(h)^J}{\partial \theta_{\nu_1} \ldots \partial
    \theta_{\nu_{q+1}}}
    = \sum_{\nu_{q+1} \in \overline{\mathcal{M}}}
    \frac{\partial}{\partial \theta_{\nu_{q+1}}} \Big(
    \sum_{\nu_1, \ldots, \nu_q \in \overline{\mathcal{M}}} \frac{\partial^q
    g(h)^J}{\partial \theta_{\nu_1} \ldots \partial
    \theta_{\nu_q}} \Big) \\
    &= \sum_{\nu_{q+1} \in \overline{\mathcal{M}}} \frac{\partial}{\partial
    \theta_{\nu_{q+1}}} \Big( \sum_{\textbf{u} \in SLTS_{q+1}^{(\pi)}}
    \sum_{K_1, \ldots, K_m = 1}^d g_{K_1 \ldots K_m}^J(h) \cdot
    (h^{K_1})^{(\delta_1)} \cdot \ldots \cdot
    (h^{K_m})^{(\delta_m)} \Big) \\
    &= \sum_{\textbf{u} \in SLTS_{q+1}^{(\pi)}} \sum_{K_1, \ldots,
    K_m, K = 1}^d g_{K_1 \ldots K_m K}^J(h) \cdot
    (h^K)^{(1)} \cdot (h^{K_1})^{(\delta_1)} \cdot \ldots \cdot
    (h^{K_m})^{(\delta_m)} \\
    &+ \sum_{\textbf{u} \in SLTS_{q+1}^{(\pi)}} \sum_{K_1, \ldots,
    K_m = 1}^d g_{K_1 \ldots K_m}^J(h) \cdot
    (h^{K_1})^{(\delta_1+1)} \cdot (h^{K_2})^{(\delta_2)} \cdot
    \ldots \cdot (h^{K_m})^{(\delta_m)} \\
    &+ \, \ldots \\
    &+ \sum_{\textbf{u} \in SLTS_{q+1}^{(\pi)}} \sum_{K_1, \ldots,
    K_m = 1}^d g_{K_1 \ldots K_m}^J(h) \cdot
    (h^{K_1})^{(\delta_1)} \cdot \ldots \cdot
    (h^{K_{m-1}})^{(\delta_{m-1})} \cdot (h^{K_m})^{(\delta_m+1)}
    \\
    &= \sum_{\textbf{u} \in SLTS_{q+2}^{(\pi)}} \sum_{K_1, \ldots,
    K_m = 1}^d g_{K_1 \ldots K_m}^J(h) \cdot
    (h^{K_1})^{(\delta_1)} \cdot \ldots \cdot
    (h^{K_m})^{(\delta_m)} .
    \end{split}
\end{equation*}
\hfill $\Box$ \\ \\
As in the deterministic setting, the density $\gamma(\textbf{t})$
of a tree is a measure of its non-bushiness and can be similarly
defined for stochastic colored trees.
\begin{Def}
    For $\textbf{t}=(\textbf{t}',\textbf{t}'') \in LTS$
    let $\gamma(\textbf{t})$ be defined
    recursively by
    \begin{alignat*}{2}
        \gamma(\textbf{t}) &= 1 &\qquad &\text{if } \, l(\textbf{t})=1, \\
        \gamma(\textbf{t}) &= \prod_{i=1}^m \gamma(\textbf{t}_i) &\qquad
        &\text{if } \, \textbf{t}=(\textbf{t}_1, \ldots, \textbf{t}_m), \\
        \gamma(\textbf{t}) &= l(\textbf{t}) \prod_{i=1}^m \gamma(\textbf{t}_i) &\qquad
        &\text{if } \, \textbf{t}=[\textbf{t}_1, \ldots, \textbf{t}_m] \,
        \text{ or } \, \textbf{t}=\{\textbf{t}_1, \ldots, \textbf{t}_m\}_{j}.
    \end{alignat*}
\end{Def}
\begin{figure}[htbp]
\begin{center}
    \def\dedge{\ncline[linestyle=dotted]}
    \begin{psmatrix}[colsep=1.0cm, rowsep=0.6cm]
    \pstree[treemode=U, dotstyle=otimes, dotsize=1.6mm, levelsep=0.0cm, radius=1.2mm, treefit=loose]
    {\Tn}{
    \pstree[treemode=U, dotstyle=otimes, dotsize=1.6mm, levelsep=0.6cm, radius=1.2mm, treefit=loose]
    {\TC*} {
    \pstree
    {\TC}{\TC*}
        \pstree
        {\TC*}{
            \pstree
            {\TC}{\TC \TC*} \TC* } }
    }
    &
    \pstree[treemode=U, dotstyle=otimes, dotsize=1.6mm, levelsep=0.0cm, radius=1.2mm, treefit=loose]
    {\Tn}{
    \pstree[treemode=U, dotstyle=otimes, dotsize=1.6mm, levelsep=0.6cm, radius=1.2mm, treefit=loose]
    {\TC*} {
    \pstree
    {\TC[edge=\dedge]}{\TC*}
        \pstree
        {\TC*[edge=\dedge]}{
            \pstree
            {\TC}{\TC \TC*} \TC* } }
    }
    &
    \pstree[treemode=U, dotstyle=otimes, dotsize=1.6mm, levelsep=0.0cm, radius=1.2mm, treefit=loose]
    {\Tn}{
    \pstree[treemode=U, dotstyle=otimes, dotsize=1.6mm, levelsep=0.6cm, radius=1.2mm, treefit=loose]
    {\TC*} {
    \pstree
    {\TC[edge=\dedge]}{\TC*[edge=\dedge]}
        \pstree
        {\TC*[edge=\dedge]}{
            \pstree
            {\TC[edge=\dedge]}{\TC \TC*} \TC*[edge=\dedge] } }
    }
    \\
    $\gamma(\textbf{t}) = \, 8 \qquad \qquad \qquad$ &
    $\cdot \, 2 \, \cdot \, 5 \qquad$ &
    $\quad \cdot \, 3 \, = 240$
    \end{psmatrix}
\caption{Example for the definition of $\gamma(\textbf{t})$ for a
tree $\textbf{t} \in LTS$.} \label{SRK-expansion-example-gamma(t)}
\end{center}
\end{figure}
In order to have a more suitable notation for the proof of the
main theorem of this section, i.e.\ the theorem about the
expansion of the approximation calculated with the stochastic
Runge-Kutta method by rooted trees, we introduce the following
notation:
\begin{Def} \label{St-Def-SRK-expansion-1}
    Let $\textbf{t}=(\textbf{t}',\textbf{t}'') \in LTS$ be a tree with
    $l=l(\textbf{t})>1$ nodes which are
    denoted by $i_1 < i_2 < \ldots < i_l$, consisting of $s=s(\textbf{t})
    \leq l$ stochastic nodes $\sigma_{j_1}, \sigma_{j_2}, \ldots,
    \sigma_{j_s}$. Then
    we denote for $i \in \{i_2, \ldots, i_l\}$ by
    \begin{equation}
        {Z}_{\textbf{t}'(i),i} = \begin{cases}
        z_i^{(0,0)} \quad & \text{ if } \, \textbf{t}''(i) = \tau \, \text{ and }
        \, \textbf{t}''(\textbf{t}'(i)) = \gamma \\
        \sum_{j_k=1}^m \sum_{\nu_k \in \mathcal{M}}
        z_i^{(j_k,\nu_k)} \quad & \text{ if } \, \textbf{t}''(i) = \sigma_{j_k} \, \text{ and }
        \, \textbf{t}''(\textbf{t}'(i)) = \gamma \\
        Z_{\textbf{t}'(i),i}^{(0,0),(0,0)} \quad & \text{ if } \, \textbf{t}''(i) = \tau
        \, \text{ and } \, \textbf{t}''(\textbf{t}'(i)) = \tau \\
        Z_{\textbf{t}'(i),i}^{(j_k,\nu_k),(0,0)} \quad & \text{ if } \, \textbf{t}''(i) = \tau
        \, \text{ and } \, \textbf{t}''(\textbf{t}'(i)) = \sigma_{j_k} \\
        \sum_{j_k=1}^m \sum_{\nu_k \in \mathcal{M}}
        Z_{\textbf{t}'(i),i}^{(0,0),(j_k,\nu_k)} \quad & \text{ if } \, \textbf{t}''(i) =
        \sigma_{j_k} \, \text{ and } \, \textbf{t}''(\textbf{t}'(i)) = \tau \\
        \sum_{j_r=1}^m \sum_{\nu_r \in \mathcal{M}}
        Z_{\textbf{t}'(i),i}^{(j_k,\nu_k),(j_r,\nu_r)} \quad & \text{ if } \, \textbf{t}''(i) =
        \sigma_{j_r} \, \text{ and } \, \textbf{t}''(\textbf{t}'(i)) = \sigma_{j_k}
        \end{cases} .
    \end{equation}
    Further, we denote by
    \begin{equation}
        \Phi_{i_1}(\textbf{t}) = \sum_{i_2, \ldots, i_l=1}^s
        Z_{\textbf{t}'(i_2),i_2} \cdot \ldots \cdot
        Z_{\textbf{t}'(i_l),i_l}
    \end{equation}
    the corresponding coefficient function and define
    $\Phi_{i_1}(\textbf{t})=1$ if $l(\textbf{t})=1$.
\end{Def}
We will now state a proposition which allows a representation of
the derivatives of the stochastic Runge-Kutta method w.r.t.\
rooted trees.
\begin{Sat} \label{Sat-Ableitungen-H(t)+Y(t)-1}
    Let $q \in \mathbb{N}$, $k \in \{0,1, \ldots, m\}$, $\nu \in \overline{\mathcal{M}}$
    and $\mathcal{A}=\{\tau, \sigma_{j_r} : r \in \mathbb{N}\}$.
    We denote by
    \begin{equation}
        \begin{split}
        z_{i_1} &= \begin{cases} z_{i_1}^{(0,0)} & \text{if } \,
        \textbf{t}''(i_1) = \tau \\
        \sum_{j_r=1}^m \sum_{\nu_r \in \mathcal{M}}
        z_{i_1}^{(j_r,\nu_r)}
        & \text{if } \, \textbf{t}''(i_1) = \sigma_{j_r}
        \end{cases} \\ 
        Z_{i,i_1}^{(k,\nu)} &= \begin{cases} Z_{i,i_1}^{(k,\nu),(0,0)}
        & \text{if } \, \textbf{t}''(i_1) = \tau \\
        \sum_{j_r=1}^m \sum_{\mu_r \in \mathcal{M}}
        Z_{i,i_1}^{(k,\nu),(j_r,\nu_r)}
        & \text{if } \, \textbf{t}''(i_1) = \sigma_{j_r}
        \end{cases} .
        \end{split}
    \end{equation}
    Then the derivatives of the $J$th component of $H_i^{(k,\nu)}(t_0)$ satisfy
    \begin{equation} \label{Sat-Ableitungen-H(t)+Y(t)-eqn1}
        \begin{split}
        \mathcal{D}^q \, H_i^{(k,\nu)}(t_0)^J
        = \sum_{\substack{\textbf{t} \in LTS \\ l(\textbf{t}) = q}}
        \gamma(\textbf{t}) \sum_{i_1=1}^s
        Z_{i,i_1}^{(k,\nu)} \cdot
        \Phi_{i_1}(\textbf{t}) \cdot F(\textbf{t})(Y(t_0))^J .
        \end{split}
    \end{equation}
    The $J$th component of the numerical solution $Y(t_0)$
    satisfies
    \begin{equation} \label{Sat-Ableitungen-H(t)+Y(t)-eqn2}
        \begin{split}
        \mathcal{D}^q \, Y(t_0)^J
        = \sum_{\substack{\textbf{t} \in LTS \\ l(\textbf{t}) = q}}
        \gamma(\textbf{t}) \sum_{i_1=1}^s z_{i_1} \cdot
        \Phi_{i_1}(\textbf{t}) \cdot F(\textbf{t})(Y(t_0))^J .
        \end{split}
    \end{equation}
\end{Sat}
{\bf Proof.} Because of the similarity of $Y(t)$ and
$H_i^{(k,\nu)}(t)$, it is sufficient to prove the first
equation~(\ref{Sat-Ableitungen-H(t)+Y(t)-eqn1}) only. Then the
second equation~(\ref{Sat-Ableitungen-H(t)+Y(t)-eqn2}) follows
substituting $Z_{ij}^{(k,\nu)(0,0)}$ and $Z_{ij}^{(k,\nu)(r,\mu)}$
by $z_{j}^{(0,0)}$ and $z_{j}^{(r,\mu)}$, respectively, in
$H_i^{(k,\nu)}(t)$, which equals it to $Y(t)$, and by the
definition of $\mathcal{D}^q$ in (\ref{Taylor-Operator-Dk}). \\ \\
We prove equation~(\ref{Sat-Ableitungen-H(t)+Y(t)-eqn1}) by
induction on $q$. For $q=1$ and $\mathcal{A}=\{\tau, \sigma_{j_k}
: k \in \mathbb{N} \}$ there are two trees $\textbf{t}_1=\tau$ and
$\textbf{t}_2=\sigma_{j_1}$ with
$l(\textbf{t}_1)=l(\textbf{t}_2)=1$ in $LTS$ and
\begin{equation}
    \begin{split}
    \sum_{\nu_1 \in \overline{\mathcal{M}}} \Delta \theta_{\nu_1} \cdot \frac{\partial
    H_i^{(k,\nu)}(t_0)^J}{\partial \theta_{\nu_1}} = &
    \sum_{i_1=1}^s Z_{i,i_1}^{(k,\nu),(0,0)} \cdot a(Y(t_0))^J \\
    & +
    \sum_{i_1=1}^s \sum_{j_1=1}^m \sum_{\nu_1 \in \mathcal{M}}
    Z_{i,i_1}^{(k,\nu),(j_1,\nu_1)} \cdot b^{j_1}(Y(t_0))^J \\
    = & \sum_{\substack{\textbf{t} \in LTS \\ l(\textbf{t})=1}} \gamma(\textbf{t})
    \sum_{i_1=1}^s Z_{i,i_1}^{(k,\nu)} \cdot F(\textbf{t})(Y(t_0))^J .
    \end{split}
\end{equation}
For a better understanding, we also consider the case $q=2$. Here
we have to consider the trees $\textbf{t}_3=[\tau]$, $\textbf{t}_4
= [\sigma_{j_1}]$, $\textbf{t}_5 = \{ \tau \}_{j_1}$ and
$\textbf{t}_6=\{ \sigma_{j_2} \}_{j_1}$ with $l(\textbf{t})=2$
nodes in $LTS$. Then we get
\begin{equation}
    \begin{split}
    & \sum_{\nu_1, \nu_2 \in \overline{\mathcal{M}}}
    \Delta \theta_{\nu_1} \cdot \Delta \theta_{\nu_2}
    \cdot \frac{\partial H_i^{(k,\nu)}(t_0)^J}{\partial \theta_{\nu_1}
    \partial \theta_{\nu_2}} \\
    & = \,
    2 \sum_{i_1,i_2=1}^s Z_{i,i_1}^{(k,\nu),(0,0)} \,
    Z_{i_1,i_2}^{(0,0),(0,0)}
    \, \sum_{K_1=1}^d \frac{\partial a(Y(t_0))^J}{\partial x^{K_1}} \, a(Y(t_0))^{K_1} \\
    & + 2 \sum_{i_1,i_2=1}^s \sum_{j_1=1}^m \sum_{\nu_1 \in \mathcal{M}}
    Z_{i,i_1}^{(k,\nu),(0,0)} \,
    Z_{i_1,i_2}^{(0,0),(j_1,\nu_1)}
    \, \sum_{K_1=1}^d \frac{\partial a(Y(t_0))^J}{\partial x^{K_1}} \, b^{j_1}(Y(t_0))^{K_1} \\
    & +
    2 \sum_{i_1,i_2=1}^s \sum_{j_1=1}^m \sum_{\nu_1 \in \mathcal{M}}
    Z_{i,i_1}^{(k,\nu),(j_1,\nu_1)} \, Z_{i_1,i_2}^{(j_1,\nu_1),(0,0)} \,
    \sum_{K_1=1}^d \frac{\partial b^{j_1} (Y(t_0))^J}{\partial x^{K_1}} \, a(Y(t_0))^{K_1} \\
    & +
    2 \sum_{i_1,i_2=1}^s \sum_{j_1,j_2=1}^m \sum_{\nu_1,\nu_2 \in \mathcal{M}}
    Z_{i,i_1}^{(k,\nu),(j_1,\nu_1)} Z_{i_1,i_2}^{(j_1,\nu_1),(j_2,\nu_2)}
    \sum_{K_1=1}^d \frac{\partial b^{j_1}(Y(t_0))^J}{\partial x^{K_1}} b^{j_2}(Y(t_0))^{K_1} \\
    & = \sum_{\substack{\textbf{t} \in LTS \\ l(\textbf{t})=2}} \gamma(\textbf{t})
    \sum_{i_1=1}^s Z_{i,i_1}^{(k,\nu)} \cdot \Phi_{i_1}(\textbf{t}) \cdot F(\textbf{t})(Y(t_0))^J .
    \end{split}
\end{equation}
Now, we assume that
equation~(\ref{Sat-Ableitungen-H(t)+Y(t)-eqn1}) holds for some
$q-1$ and prove the case $q$. The first step is the application of
formula~(\ref{Lem-Leibniz-formula-2-eqn1}) in order to obtain
\begin{equation} \label{Sat-Ableitungen-H(t)+Y(t)-1-proof-eqn0}
    \begin{split}
    \mathcal{D}^q & H_i^{(k,\nu)}(t_0)^J =
    \sum_{\nu_1, \ldots, \nu_q \in \overline{\mathcal{M}}}
    \Delta \theta_{\nu_1} \cdot
    \ldots \cdot \Delta \theta_{\nu_q} \cdot \frac{\partial^q
    H_i^{(k,\nu)}(t_0)^J}{\partial \theta_{\nu_1} \ldots \partial
    \theta_{\nu_q}} \\
    & = \, q \sum_{i_1=1}^s Z_{i,i_1}^{(k,\nu),(0,0)}
    \sum_{\nu_1, \ldots, \nu_{q-1} \in \overline{\mathcal{M}}} \Delta \theta_{\nu_1}
    \cdot \ldots \cdot \Delta \theta_{\nu_{q-1}}
    \frac{\partial^{q-1} a(H_{i_1}^{(0,0)}(t_0))^J}{\partial
    \theta_{\nu_1} \ldots \theta_{\nu_{q-1}}} \\
    & + q \sum_{i_1=1}^s \sum_{j_1=1}^m \sum_{\nu_1 \in \mathcal{M}}
    Z_{i,i_1}^{(k,\nu)(j_1,\nu_1)} \times \\
    & \times
    \sum_{\nu_1,
    \ldots, \nu_{q-1} \in \overline{\mathcal{M}}} \Delta \theta_{\nu_1} \cdot \ldots
    \cdot \Delta \theta_{\nu_{q-1}} \frac{\partial^{q-1}
    b^{j_1}(H_{i_1}^{(j_1,\nu_1)}(t_0))^J}{\partial \theta_{\nu_1} \ldots
    \partial \theta_{\nu_{q-1}}} .
    \end{split}
\end{equation}
As the second step, we make use of Lemma~\ref{Lem-Faa-di-Bruno-1}
twice. First, equation~(\ref{Lem-Faa-di-Bruno-1-eqn1}) is applied
to trees $\textbf{u} \in SLTS_q^{(\tau)}$ (i.e., trees having a
root of type $\tau$) and second, to trees $\textbf{u} \in
SLTS_q^{(\sigma_{j_1})}$ (i.e., trees having a root of type
$\sigma_{j_1}$). Thus with $\delta_1+ \ldots +
\delta_{m(\textbf{u})} = q-1$ we obtain
\begin{equation} \label{Sat-Ableitungen-H(t)+Y(t)-1-proof-eqn1}
    \begin{split}
    \sum_{\nu_1, \ldots, \nu_{q-1} \in \overline{\mathcal{M}}} \Delta \theta_{\nu_1} &
    \cdot \ldots \cdot \Delta \theta_{\nu_{q-1}}
    \frac{\partial^{q-1}
    a(H_{i_1}^{(0,0)}(t_0))^J}{\partial
    \theta_{\nu_1} \ldots \partial \theta_{\nu_{q-1}}} \\
    = & \sum_{\textbf{u} \in SLTS_q^{(\tau)}} \sum_{K_1, \ldots,
    K_m=1}^d a_{K_1 \ldots K_m}^J(H_{i_1}^{(0,0)}(t_0))
    \, \times\\
    & \times \left( \left( \sum_{\nu_1, \ldots, \nu_{\delta_1} \in
    \overline{\mathcal{M}}}
    \Delta \theta_{\nu_1} \cdot \ldots \cdot \Delta
    \theta_{\nu_{\delta_1}} \frac{\partial^{\delta_1}
    H_{i_1}^{(0,0)}(t_0)^{K_1}}{\partial \theta_{\nu_1} \ldots
    \partial \theta_{\nu_{\delta_1}}} \right) \cdot \ldots \right. \\
    & \left. \ldots \cdot
    \left( \sum_{\nu_1, \ldots, \nu_{\delta_m} \in \overline{\mathcal{M}}}
    \Delta \theta_{\nu_1} \cdot \ldots \cdot \Delta
    \theta_{\nu_{\delta_m}} \frac{\partial^{\delta_m}
    H_{i_1}^{(0,0)}(t_0)^{K_m}}{\partial \theta_{\nu_1} \ldots
    \partial \theta_{\nu_{\delta_m}}}
    \right)
    \right)
    \end{split}
\end{equation}
and analogously
\begin{equation} \label{Sat-Ableitungen-H(t)+Y(t)-1-proof-eqn2}
    \begin{split}
    \sum_{\nu_1, \ldots, \nu_{q-1} \in \overline{\mathcal{M}}}
    \Delta \theta_{\nu_1} &
    \cdot \ldots \cdot \Delta \theta_{\nu_{q-1}}
    \frac{\partial^{q-1}
    b^{j_1}(H_{i_1}^{(j_1,\nu_1)}(t_0))^J}{\partial
    \theta_{\nu_1} \ldots \partial \theta_{\nu_{q-1}}} \\
    = & \sum_{\textbf{u} \in SLTS_q^{(\sigma_{j_1})}} \sum_{K_1, \ldots,
    K_m=1}^d {b^{j_1}}_{K_1 \ldots K_m}^J(H_{i_1}^{(j_1,\nu_1)}(t_0))
    \, \times\\
    & \times \left( \left( \sum_{\nu_1, \ldots, \nu_{\delta_1} \in \overline{\mathcal{M}}}
    \Delta \theta_{\nu_1} \cdot \ldots \cdot \Delta
    \theta_{\nu_{\delta_1}} \frac{\partial^{\delta_1}
    H_{i_1}^{(j_1,\nu_1)}(t_0)^{K_1}}{\partial \theta_{\nu_1} \ldots
    \partial \theta_{\nu_{\delta_1}}} \right) \cdot \ldots \right. \\
    & \left. \ldots \cdot
    \left( \sum_{\nu_1, \ldots, \nu_{\delta_m} \in \overline{\mathcal{M}}}
    \Delta \theta_{\nu_1} \cdot \ldots \cdot \Delta
    \theta_{\nu_{\delta_m}} \frac{\partial^{\delta_m}
    H_{i_1}^{(j_1,\nu_1)}(t_0)^{K_m}}{\partial \theta_{\nu_1} \ldots
    \partial \theta_{\nu_{\delta_m}}}
    \right)
    \right) .
    \end{split}
\end{equation}
Finally, we replace the derivatives of $H_{i_1}^{(0,0)}$ and
$H_{i_1}^{(j_1,\nu_1)}$, which appear
in~(\ref{Sat-Ableitungen-H(t)+Y(t)-1-proof-eqn1})
and~(\ref{Sat-Ableitungen-H(t)+Y(t)-1-proof-eqn2}) with $\delta_i
\leq q-1$, $1 \leq i \leq m=m(\textbf{u})$, by the induction
hypothesis (\ref{Sat-Ableitungen-H(t)+Y(t)-eqn1}) and rearrange
the sums. Then we get
for~(\ref{Sat-Ableitungen-H(t)+Y(t)-1-proof-eqn0}):
\begin{equation} \label{Sat-Ableitungen-H(t)+Y(t)-1-proof-eqn3}
    \begin{split}
    & \sum_{\nu_1, \ldots, \nu_q \in \overline{\mathcal{M}}}
    \Delta \theta_{\nu_1} \cdot
    \ldots \cdot \Delta \theta_{\nu_q} \frac{\partial^q
    H_i^{(k,\nu)}(t_0)^J}{\partial \theta_{\nu_1} \ldots \partial
    \theta_{\nu_q}} \\
    & = \, q \sum_{\textbf{u} \in SLTS_q^{(\tau)}}
    \sum_{\substack{\textbf{t}_1 \in LTS
    \\ l(\textbf{t}_1) = \delta_1}} \ldots
    \sum_{\substack{\textbf{t}_m \in LTS \\
    l(\textbf{t}_m) = \delta_m}} \gamma(\textbf{t}_1) \cdot \ldots \cdot
    \gamma(\textbf{t}_m) \, \times \\
    & \times \sum_{i_1=1}^s Z_{i,i_1}^{(k,\nu)(0,0)} \left( \sum_{k_1=1}^s
    Z_{i_1,k_1}^{(0,0)} \Phi_{k_1}(\textbf{t}_1) \cdot \ldots \cdot \sum_{k_m=1}^s
    Z_{i_1,k_m}^{(0,0)} \Phi_{k_m}(\textbf{t}_m) \right) \times \\
    & \times \sum_{K_1, \ldots, K_m=1}^d a_{K_1 \ldots
    K_m}^J(H_{i_1}^{(0,0)}(t_0)) \cdot \left( F(\textbf{t}_1)(Y(t_0))^{K_1} \cdot
    \ldots \cdot F(\textbf{t}_m)(Y(t_0))^{K_m} \right) \\
    & + q
    \sum_{\textbf{u} \in SLTS_q^{(\sigma_{j_1})}}
    \sum_{\substack{\textbf{t}_1 \in
    LTS \\ l(\textbf{t}_1) = \delta_1}} \ldots
    \sum_{\substack{\textbf{t}_m \in LTS
    \\ l(\textbf{t}_m) = \delta_m}} \gamma(\textbf{t}_1) \cdot \ldots \cdot
    \gamma(\textbf{t}_m) \, \times \\
    & \times \sum_{j_1=1}^m \sum_{\nu_1 \in \mathcal{M}} \sum_{i_1=1}^s
    Z_{i,i_1}^{(k,\nu)(j_1,\nu_1)} \left( \sum_{k_1=1}^s
    Z_{i_1,k_1}^{(j_1,\nu_1)} \Phi_{k_1}(\textbf{t}_1) \cdot
    \ldots \cdot \sum_{k_m=1}^s
    Z_{i_1,k_m}^{(j_1,\nu_1)} \Phi_{k_m}(\textbf{t}_m) \right) \times \\
    & \times \sum_{K_1, \ldots, K_m=1}^d {b^{j_1}}_{K_1 \ldots
    K_m}^J(H_{i_1}^{(j_1,\nu_1)}(t_0)) \cdot \left(
    F(\textbf{t}_1)(Y(t_0))^{K_1} \cdot
    \ldots \cdot F(\textbf{t}_m)(Y(t_0))^{K_m} \right)
    \end{split}
\end{equation}
where $i_1$ denotes the root of $\textbf{u}$ and $k_1, \ldots,
k_m$ denote the roots of the trees $\textbf{t}_1, \ldots,
\textbf{t}_m$, respectively. \\ \\
The main difficulty is now to understand that to each tuple of
trees
\begin{equation*}
    (\textbf{u}, \textbf{t}_1, \ldots, \textbf{t}_m) \quad \text{with}
    \quad \textbf{u} \in
    SLTS_q^{(\pi)}, \quad \textbf{t}_i \in LTS, \quad l(\textbf{t}_i) = \delta_i
\end{equation*}
with $\pi \in \mathcal{A}$ and $\sum_{i=1}^m \delta_i = q-1$,
there corresponds exactly one labelled tree
$\textbf{t}=(\textbf{t}',\textbf{t}'') \in LTS$ with
$l(\textbf{t})=q$ such that the root $i_1$ of $\textbf{t}$ is of
type $\pi$ and such that
\begin{equation} \label{Sat-Ableitungen-H(t)+Y(t)-1-proof-eqn4}
    \gamma(\textbf{t}) = q \cdot \gamma(\textbf{t}_1) \cdot \ldots \cdot
    \gamma(\textbf{t}_m)
\end{equation}
and for $\pi=\tau$
\begin{equation} \label{Sat-Ableitungen-H(t)+Y(t)-1-proof-eqn5}
    \begin{split}
    F(\textbf{t})(Y(t_0))^J &= \sum_{K_1, \ldots, K_m=1}^d a_{K_1 \ldots K_m}^J
    (Y(t_0)) \, F(\textbf{t}_1)(Y(t_0))^{K_1} \ldots
    F(\textbf{t}_m)(Y(t_0))^{K_m} \\
    \Phi_{i_1}(\textbf{t}) &= \sum_{k_1, \ldots, k_m=1}^s
    Z_{i_1,k_1}^{(0,0)} \ldots Z_{i_1,k_m}^{(0,0)}
    \Phi_{k_1}(\textbf{t}_1) \ldots \Phi_{k_m}(\textbf{t}_m)
    \end{split}
\end{equation}
or for $\pi = \sigma_{j_1}$
\begin{equation} \label{Sat-Ableitungen-H(t)+Y(t)-1-proof-eqn6}
    \begin{split}
    F(\textbf{t})(Y(t_0))^J &= \sum_{K_1, \ldots, K_m=1}^d {b^{j_1}}_{K_1 \ldots
    K_m}^J (Y(t_0)) \, F(\textbf{t}_1)(Y(t_0))^{K_1} \ldots
    F(\textbf{t}_m)(Y(t_0))^{K_m} \\
    \Phi_{i_1}(\textbf{t}) &= \sum_{k_1, \ldots, k_m=1}^s
    Z_{i_1,k_1}^{(j_1,\nu_1)} \ldots Z_{i_1,k_m}^{(j_1,\nu_1)}
    \Phi_{k_1}(\textbf{t}_1) \ldots \Phi_{k_m}(\textbf{t}_m)
    \end{split}
\end{equation}
holds, respectively. This labelled tree $\textbf{t}$ is obtained
if the branches of $\textbf{u}$ are replaced by the trees
$\textbf{t}_1, \ldots, \textbf{t}_m$ and the corresponding labels
are taken over in a natural way, i.e., in the same order (see
Figure~\ref{SRK-expansion-proof-tree-bijection}).
\begin{figure}[htbp]
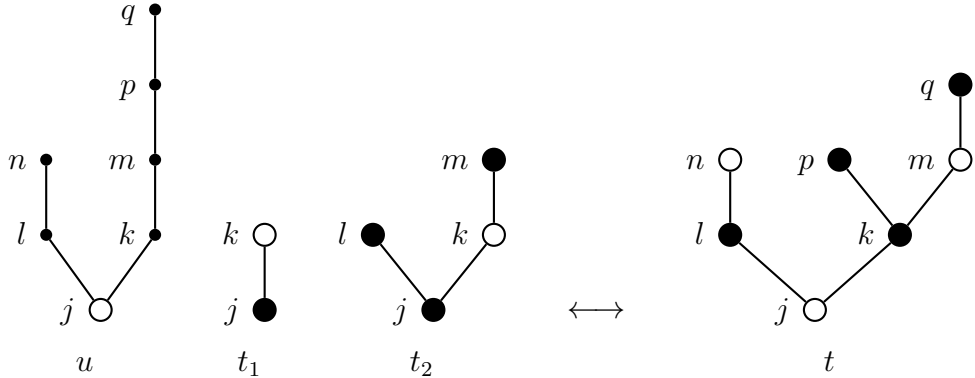

\begin{center}
    \begin{psmatrix}[colsep=0.2cm, rowsep=0.2cm]
    \pstree[treemode=U, dotstyle=otimes, dotsize=3.2mm, levelsep=0.1cm, radius=1.6mm, treefit=loose]
    {\Tn}{
    \pstree[treemode=U, dotstyle=otimes, dotsize=1.6mm, levelsep=1cm, radius=0.8mm, treefit=loose]
    {\TC~[tnpos=l]{$j$}}{
        \pstree{\TC*~[tnpos=l]{$l$}}{\TC*~[tnpos=l]{$n$}}
        \pstree{\TC*~[tnpos=l]{$k$}}{
            \pstree{\TC*~[tnpos=l]{$m$}}{
                \pstree{\TC*~[tnpos=l]{$p$}}{\TC*~[tnpos=l]{$q$}}}}}
    }
    & $\quad$ &
    \pstree[treemode=U, dotstyle=otimes, dotsize=3.2mm, levelsep=0.1cm, radius=1.6mm, treefit=loose]
    {\Tn}{
    \pstree[treemode=U, dotstyle=otimes, dotsize=3.2mm, levelsep=1cm, radius=1.6mm, treefit=loose]
    {\TC*~[tnpos=l]{$j$}} {\TC~[tnpos=l]{$k$}}
    }
    & $\quad$ &
    \pstree[treemode=U, dotstyle=otimes, dotsize=3.2mm, levelsep=0.1cm, radius=1.6mm, treefit=loose]
    {\Tn}{
    \pstree[treemode=U, dotstyle=otimes, dotsize=3.2mm, levelsep=1cm, radius=1.6mm, treefit=loose]
    {\TC*~[tnpos=l]{$j$}}{
        \TC*~[tnpos=l]{$l$} \pstree{\TC~[tnpos=l]{$k$}}{\TC*~[tnpos=l]{$m$}}}
    }
    & $\quad$ & $\longleftrightarrow$ & $\quad$ &
    \pstree[treemode=U, dotstyle=otimes, dotsize=3.2mm, levelsep=0.1cm, radius=1.6mm, treefit=loose]
    {\Tn}{
    \pstree[treemode=U, dotstyle=otimes, dotsize=3.2mm, levelsep=1cm, radius=1.6mm, treefit=loose]
    {\TC~[tnpos=l]{$j$}}{
        \pstree{\TC*~[tnpos=l]{$l$}}{\TC~[tnpos=l]{$n$}}
            \pstree{\TC*~[tnpos=l]{$k$}}{
            \TC*~[tnpos=l]{$p$}
            \pstree{\TC~[tnpos=l]{$m$}}{\TC*~[tnpos=l]{$q$}}}}
    }
    \\
    $u$ & & $t_1$ & & $t_2$ & & & & $t$
    \end{psmatrix}
\caption{Example for the bijection of $(\textbf{u}, \textbf{t}_1,
\ldots, \textbf{t}_m) \leftrightarrow \textbf{t}$ with
$\pi=\sigma$.} \label{SRK-expansion-proof-tree-bijection}
\end{center}
\end{figure}
\\ \\
In this way, for $\pi=\tau$ and $\pi=\sigma_{j_1}$ {\emph{all}}
trees $\textbf{t}=(\textbf{t}',\textbf{t}'') \in LTS$ with
$l(\textbf{t})=q$ appear exactly {\emph{once}}. Thus
(\ref{Sat-Ableitungen-H(t)+Y(t)-1-proof-eqn3}) becomes
(\ref{Sat-Ableitungen-H(t)+Y(t)-eqn1}) after inserting
(\ref{Sat-Ableitungen-H(t)+Y(t)-1-proof-eqn4}),
(\ref{Sat-Ableitungen-H(t)+Y(t)-1-proof-eqn5}) and
(\ref{Sat-Ableitungen-H(t)+Y(t)-1-proof-eqn6}), respectively.
\hfill $\Box$ \\ \\
Since the Taylor expansion contains the coefficients of the SRK
method, we define a coefficient function $\Phi_S$ which assigns to
every tree an {\emph{elementary weight}}. So for every $\textbf{t}
\in TS$ or $\textbf{t} \in LTS$ the function $\Phi_S$ is defined
recursively by
\begin{equation}
    \Phi_S(\textbf{t}) =
    \begin{cases}
    \displaystyle
    \prod_{i=1}^{\lambda} \Phi_S(\textbf{t}_i) & \text{ if } \textbf{t}=(\textbf{t}_1, \ldots,
        \textbf{t}_{\lambda}) \\
    \displaystyle
    {z^{(0,0)}}^T \prod_{i=1}^{\lambda} \Psi^{(0,0)}(\textbf{t}_i) & \text{ if }
        \textbf{t}=[\textbf{t}_1, \ldots, \textbf{t}_{\lambda}] \\
    \displaystyle
        \sum_{\nu \in \mathcal{M}} {z^{(k,\nu)}}^T \prod_{i=1}^{\lambda}
        \Psi^{(k,\nu)}(\textbf{t}_i) &
        \text{ if } \textbf{t}=\{\textbf{t}_1, \ldots, \textbf{t}_{\lambda}\}_k
    \end{cases}
\end{equation}
\quad \\
\noindent where
$\Psi^{(0,0)}(\emptyset)=\Psi^{(k,\nu)}(\emptyset)=e$ with $\gamma
= (\emptyset)$, $\tau = [\emptyset]$, $\sigma_k = \{\emptyset\}_k$
and
\begin{equation}
    \Psi^{(k,\nu)}(\textbf{t}) =
    \begin{cases}
        \displaystyle
        Z^{(k,\nu),(0,0)} \prod_{i=1}^{\lambda} \Psi^{(0,0)}(\textbf{t}_i) & \text{ if }
        \textbf{t}=[\textbf{t}_1, \ldots, \textbf{t}_{\lambda}] \\
    \displaystyle
        \sum_{\mu \in \mathcal{M}} Z^{(k,\nu),(r,\mu)} \prod_{i=1}^{\lambda}
        \Psi^{(r,\mu)}(\textbf{t}_i) & \text{ if }
        \textbf{t}=\{\textbf{t}_1, \ldots, \textbf{t}_{\lambda}\}_r
    \end{cases} .
\end{equation}
Here $e=(1, \ldots, 1)^T$ and the product of vectors in the
definition of $\Psi^{(0,0)}$ and $\Psi^{(k,\nu)}$ is defined by
component-wise multiplication, i.e.\ with $(a_1, \ldots, a_n) *
(b_1, \ldots, b_n) = (a_1 b_1, \ldots, a_n b_n)$.
Now we get immediately the following representation of the
stochastic Runge-Kutta approximation w.r.t.\ rooted trees.
\begin{Kor}
    Assume that the drift $a$ and the diffusion $b^j$, $1 \leq j \leq m$, are
    sufficiently differentiable. Then, the one-step approximation
    $Y(t)=Y(t_0+h)$ with $h \in \, ]0,\infty[$, given by the stochastic Runge-Kutta
    method~(\ref{St-srk-method-m}), can be re\-pre\-sent\-ed as
    \begin{equation}
        \begin{split}
        Y(t)^J & = Y(t_0)^J + \sum_{\substack{\textbf{t} \in LTS \\
        l(\textbf{t}) \leq n}} \frac{\gamma(\textbf{t})
        \sum_{i_1=1}^s z_{i_1} \, \Phi_{i_1}(\textbf{t}) \,
        F(\textbf{t})(Y(t_0))^J}{l(\textbf{t})!}
        + \mathcal{R}_n(t,t_0) \\
        & = Y(t_0)^J + \sum_{\substack{\textbf{t} \in TS \\ l(\textbf{t}) \leq n}}
        \frac{\alpha(\textbf{t}) \, \gamma(\textbf{t})
        \sum_{i_1=1}^s z_{i_1} \, \Phi_{i_1}(\textbf{t}) \,
        F(\textbf{t})(Y(t_0))^J}{l(\textbf{t})!}
        + \mathcal{R}_n(t,t_0)
        \end{split}
    \end{equation}
    for $n \in \mathbb{N}$ and with $\alpha(\textbf{t})$ denoting
    the cardinality of the tree $\textbf{t} \in
    LTS$ with $\mathcal{A}=\{\tau,\sigma_{j_k} : k \in \mathbb{N}\}$.
    Using the coefficient function $\Phi_S$, we get analogously
    \begin{equation}
        \begin{split}
        Y(t)^J & = Y(t_0)^J + \sum_{\substack{\textbf{t} \in LTS \\
        l(\textbf{t}) \leq n}} \sum_{j_1, \ldots,
        j_{s(\textbf{t})}=1}^m
        \frac{\gamma(\textbf{t}) \, \Phi_S(\textbf{t}) \,
        F(\textbf{t})(Y(t_0))^J}{l(\textbf{t})!}
        + \mathcal{R}_n(t,t_0) \\
        & = Y(t_0)^J + \sum_{\substack{\textbf{t} \in TS \\
        l(\textbf{t}) \leq n}} \sum_{j_1, \ldots,
        j_{s(\textbf{t})}=1}^m
        \frac{\alpha(\textbf{t}) \, \gamma(\textbf{t}) \,
        \Phi_S(\textbf{t}) \, F(\textbf{t})(Y(t_0))^J}{l(\textbf{t})!}
        + \mathcal{R}_n(t,t_0) .
        \end{split}
    \end{equation}
\end{Kor}
{\bf Proof.} This follows directly from the Theorem of Taylor
(see~(\ref{SRK-w1-Taylor-D})) and
Proposition~\ref{Sat-Ableitungen-H(t)+Y(t)-1}. \hfill $\Box$ \\ \\
As a final step, we extend this representation of the
approximation $Y(t)$ to our primary problem of a representation
for $f(Y(t))$. Therefore we consider a suitable subset
$LTS(\Delta)$ of $LTS$ w.r.t.\ the set $\mathcal{A}=\{\gamma,
\tau, \sigma_{j_k} : k \in \mathbb{N} \}$, where $\gamma$
represents the function $f$.
\begin{Def}
    Let $LTS(\Delta)$ denote the set of trees $\textbf{t}=(\textbf{t}',\textbf{t}'') \in LTS$
    w.r.t.\ $\mathcal{A}=\{\gamma, \tau, \sigma_{j_k} : k \in \mathbb{N}\}$
    such that
    \begin{enumerate}[a)]
        \item the root is of type $\textbf{t}''(1) = \gamma$  and all other nodes
        are either deterministic or stochastic nodes, i.e.\
        $\textbf{t}''(i) \in \{\tau, \sigma_{j_k} : k \in \mathbb{N}\}$
        for $2 \leq i \leq l(\textbf{t})$,
        \item all stochastic nodes own a different variable index $j_k$, $1 \leq k \leq s(\textbf{t})$,
        i.e.\ for two different stochastic nodes $i \neq l$ holds $\textbf{t}''(i) \neq
        \textbf{t}''(l)$.
    \end{enumerate}
    Further $TS(\Delta) = LTS(\Delta)/ \sim$
    denotes the equivalence class under the relation of
    Definition~\ref{St-tree-equivalence:Wm} restricted to $LTS(\Delta)$
    and $\alpha_{\Delta}(\textbf{t})$ denotes the cardinality of $\textbf{t}$ in $LTS$.
\end{Def}
Here it has to be pointed out that $LTS(I) \subset LTS(S) \subset
LTS(\Delta)$ since the rules of construction for the trees
$\textbf{t}$ in $LTS(I)$ and in $LTS(S)$ are more restrictive than
for the trees $\textbf{t} \in LTS(\Delta)$. However in contrast to
$LTS(I)$ and $LTS(S)$, a tree $\textbf{t} \in LTS(\Delta)$ has
$s(\textbf{t})$ different variable indices $j_1, \ldots, j_{s(t)}$
while a tree $\textbf{u}$ in $LTS(I)$ or $LTS(S)$ has only
$n(\textbf{u})=s(\textbf{u})/2$ different variable indices. For
example, the tree $(\{[\sigma_{j_2}]\}_{j_1})$ is an element of
$LTS(\Delta)$ while it is neither an element of $LTS(I)$ nor of
$LTS(S)$.
With the definition of the set $LTS(\Delta)$, we can now formulate
our main result for the expansion of the stochastic Runge-Kutta
method. It provides an expansion of $f(Y(t))$ which is required
for the calculation of order conditions for the SRK method.
\begin{The} \label{St-Theorem-Expansion-SRK-main-1:Wm}
    For the one-step approximation $Y(t)=Y(t_0+h)$, $h \in \, ]0, \infty[ \,$,
    defined by the stochastic Runge-Kutta
    method~(\ref{St-srk-method-m}), a function $f : \mathbb{R}^d
    \rightarrow \mathbb{R}$ and for $n \in \mathbb{N}$ the expansion
    \begin{equation} \label{St-Theorem-Expansion-SRK-main-eqn1:Wm}
        \begin{split}
        f(Y(t)) &= \sum_{\substack{\textbf{t} \in LTS(\Delta) \\
        l(\textbf{t})-1 \leq n}}
        \sum_{j_1, \ldots, j_{s(\textbf{t})}=1}^m
        \frac{\gamma(\textbf{t}) \cdot \Phi_S(\textbf{t}) \cdot F(\textbf{t})(Y(t_0))}{(l(\textbf{t})-1)!}
        + \mathcal{R}_n(t,t_0) \\
        &= \sum_{\substack{\textbf{t} \in TS(\Delta) \\ l(\textbf{t})-1 \leq n}}
        \sum_{j_1, \ldots, j_{s(\textbf{t})}=1}^m
        \frac{\alpha_{\Delta}(\textbf{t}) \cdot \gamma(\textbf{t}) \cdot \Phi_S(\textbf{t})
        \cdot F(\textbf{t})(Y(t_0))}{(l(\textbf{t})-1)!}
        + \mathcal{R}_n(t,t_0)
        \end{split}
    \end{equation}
    holds provided all necessary derivatives of $f$,
    $a$ and $b^{j}$, $1 \leq j \leq m$, exist.
\end{The}
{\bf Proof.} Let $\mathcal{A}=\{\gamma, \tau, \sigma_{j_k} : k \in
\mathbb{N}\}$. We apply Lemma~\ref{Lem-Faa-di-Bruno-1} with
$\pi=\gamma$ and conclude that
\begin{equation}
    \begin{split}
    \mathcal{D}^q & f(Y(t_0)) = \sum_{\nu_1, \ldots, \nu_q \in \overline{\mathcal{M}}}
    \Delta \theta_{\nu_1} \cdot \ldots \cdot \Delta \theta_{\nu_q}
    \cdot \frac{\partial^q f(Y(t_0))}{\partial \theta_{\nu_1}
    \ldots \partial \theta_{\nu_q}} \\
    = & \sum_{\textbf{u} \in SLTS_{q+1}^{(\gamma)}} \sum_{K_1, \ldots,
    K_m=1}^d f_{K_1 \ldots K_m}(Y(t_0)) \cdot \big(
    \mathcal{D}^{\delta_1} Y(t_0)^{K_1} \ldots
    \mathcal{D}^{\delta_m} Y(t_0)^{K_m} \big)
    \end{split}
\end{equation}
where $m=m(\textbf{u})$ and $\delta_1 + \ldots + \delta_m = q$.
Now Proposition~\ref{Sat-Ableitungen-H(t)+Y(t)-1} yields
\begin{equation}
    \begin{split}
    \mathcal{D}^q f(Y(t_0)) = & \sum_{\textbf{u} \in SLTS_{q+1}^{(\gamma)}}
    \sum_{K_1, \ldots, K_m=1}^d f_{K_1 \ldots K_m} (Y(t_0)) \, \times \\
    & \times \Big(
    \big( \sum_{\substack{\textbf{t}_1 \in LTS \\ l(\textbf{t}_1) = \delta_1}}
    \gamma(\textbf{t}_1) \sum_{k_1=1}^s z_{k_1} \Phi_{k_1}(\textbf{t}_1) \cdot
    F(\textbf{t}_1)(Y(t_0))^{K_1} \big) \cdot \ldots \times \\
    & \times \ldots \cdot
    \big( \sum_{\substack{\textbf{t}_m \in LTS \\ l(\textbf{t}_m) = \delta_m}}
    \gamma(\textbf{t}_m) \sum_{k_m=1}^s z_{k_m} \Phi_{k_m}(\textbf{t}_m) \cdot
    F(\textbf{t}_m)(Y(t_0))^{K_m} \big) \Big)
    \end{split}
\end{equation}
where $\textbf{t}_1, \ldots, \textbf{t}_m \in LTS$ are considered
w.r.t.\ $\mathcal{A}=\{\tau, \sigma_{j_k} : k \in \mathbb{N}\}$
and $k_1, \ldots, k_m$ denote the roots of the trees
$\textbf{t}_1, \ldots, \textbf{t}_m$, respectively. Now nearly the
same considerations as in the proof of
Proposition~\ref{Sat-Ableitungen-H(t)+Y(t)-1} apply: To each tuple
of trees $(\textbf{u}, \textbf{t}_1, \ldots, \textbf{t}_m)$ with
$\textbf{u} \in SLTS_{q+1}^{(\gamma)}$, $\textbf{t}_i \in LTS$,
$l(\textbf{t}_i) = \delta_i$ and with $\sum_{i=1}^m \delta_i = q$,
there corresponds exactly one labelled tree
$\textbf{t}=(\textbf{t}',\textbf{t}'') \in LTS(\Delta)$ with
$l(\textbf{t})=q+1$ such that the root $i_1$ of $\textbf{t}$ is of
type $\gamma$ and
\begin{equation}
    \begin{split}
    \gamma(\textbf{t}) & = \gamma(\textbf{t}_1) \cdot \ldots \cdot
    \gamma(\textbf{t}_m) \\
    F(\textbf{t})(Y(t_0)) & = \sum_{K_1, \ldots, K_m=1}^d f_{K_1 \ldots
    K_m} (Y(t_0)) \cdot F(\textbf{t}_1)(Y(t_0))^{K_1} \ldots
    F(\textbf{t}_m)(Y(t_0))^{K_m} \\
    \tilde{\Phi}(\textbf{t}) :&= \prod_{k \in {\textbf{t}'}^{-1}(i_1)} \sum_{k=1}^s
    z_k \Phi_k(\textbf{t}_k) = \sum_{k_1, \ldots, k_m=1}^s z_{k_1}
    \Phi_{k_1}(\textbf{t}_1) \cdot \ldots \cdot z_{k_m} \Phi_{k_m}(\textbf{t}_m)
    \end{split}
\end{equation}
where $\textbf{t}_k$ denotes the subtree of $\textbf{t}$ having
the node $k$ as a root. \\ \\
The labelled tree $\textbf{t}$ is obtained if the branches of
$\textbf{u}$ are replaced by the trees $\textbf{t}_1, \ldots,
\textbf{t}_m$ and the corresponding labels are taken over in a
natural way, i.e.\ in the same order (see
Figure~\ref{SRK-expansion-proof-tree-bijection}).
In this way {\emph{all}} trees
$\textbf{t}=(\textbf{t}',\textbf{t}'') \in LTS(\Delta)$ with
$l(\textbf{t})=q+1$ appear exactly {\emph{once}}. Applying the
usual tensor notation and substituting $\tilde{\Phi}(\textbf{t})$
by $\Phi_S(\textbf{t})$, we get
\begin{equation}
    \begin{split}
    \mathcal{D}^q f(Y(t_0)) & = \sum_{\substack{\textbf{t} \in LTS(\Delta)
    \\ l(\textbf{t})=q+1}} \gamma(\textbf{t}) \cdot \tilde{\Phi}(\textbf{t})
    \cdot F(\textbf{t})(Y(t_0)) \\
    & = \sum_{\substack{\textbf{t} \in LTS(\Delta) \\
    l(\textbf{t})=q+1}} \sum_{j_1, \ldots, j_{s(\textbf{t})} =
    1}^m \gamma(\textbf{t}) \cdot \Phi_S(\textbf{t}) \cdot
    F(\textbf{t})(Y(t_0)) .
    \end{split}
\end{equation}
With $\Phi_S(\gamma)=1$, $F(\gamma)(Y(t_0))=f(Y(t_0))$ and the
Theorem of Taylor~(\ref{SRK-w1-Taylor-D}) we finally arrive at
(\ref{St-Theorem-Expansion-SRK-main-eqn1:Wm}) which completes the
proof. \hfill $\Box$
\section{Order Conditions for Stochastic Runge-Kutta Methods}
\label{Sec:general-cond-order-conv-p}
In this section, conditions such that the stochastic Runge-Kutta
method~(\ref{St-srk-method-m}) converges in the weak sense with
order $p$ to the solution of the stochastic differential
equation~(\ref{Ito-St-SDE1-autonom-Wm}) are considered. Therefore,
we give a suitable representation of the approximation due to the
SRK method.
\begin{Sat} \label{St-SRK-expansion-expectation-1:Wm}
Let $Y(t)=Y(t_0+h)$ with $h \in \, ]0,h_0[ \,$ and $Y(t_0)=x_0$
denote the one-step approximation defined by the stochastic
Runge-Kutta method~(\ref{St-srk-method-m}). Assume that for the
random variables holds $\theta_{\iota}(h) = \sqrt{h} \cdot
\vartheta_{\iota}$ for $\iota \in \mathcal{M}$ with a bounded
random variable $\vartheta_{\iota}$. Then for $f : \mathbb{R}^d
\rightarrow \mathbb{R}$ and $p \in \mathbb{N}$ the expansion
\begin{equation} \label{St-SRK-expansion-expectation-eqn1:Wm}
    {E}^{t_0,x_0} \left( f \left( Y(t) \right) \right)
    = \sum_{\substack{\textbf{t} \in
    TS(\Delta) \\ \rho(\textbf{t}) \leq p + \tfrac{1}{2}}}
    \sum_{j_1, \ldots, j_{s(\textbf{t})}=1}^m
    \frac{ \alpha_{\Delta}(\textbf{t}) \, \gamma(\textbf{t})
    \, F(\textbf{t})(x_0) \, {E}\left(\Phi_S(\textbf{t}) \right)}
    {(l(\textbf{t})-1)!}
    + O \left( h^{p+1} \right)
\end{equation}
holds for sufficient small $h_0>0$, provided $f, a^i, b^{i,j} \in
C_P^{2(p+1)}(\mathbb{R}^d, \mathbb{R})$ for all $i=1, \ldots, d$
and $j=1, \ldots,m$.
\end{Sat}
{\bf Proof.} Apply
Theorem~\ref{St-Theorem-Expansion-SRK-main-1:Wm} with $n =
2(p+\tfrac{1}{2})$ and simply take the expectation of
equation~(\ref{St-Theorem-Expansion-SRK-main-eqn1:Wm}). By the
definition of $\Phi_S$ and due to
(\ref{St-SRK-moment-condition-1:Wm}), for all $\textbf{t} \in
TS(\Delta)$ the expectation becomes
\[
    E(\Phi_S(\textbf{t})) = O( h^{d(\textbf{t})+ \tfrac{1}{2}
    s(\textbf{t})} ) = O(h^{\rho(\textbf{t})}).
\]
Now, for all trees $\textbf{t} \in TS(\Delta)$ appearing in the
sum of equation~(\ref{St-Theorem-Expansion-SRK-main-eqn1:Wm}) and
which do not appear in the sum of
(\ref{St-SRK-expansion-expectation-eqn1:Wm}), i.e.\ trees with
$l(\textbf{t}) \leq 2p+2$ and $\rho(\textbf{t}) \geq p + 1$, we
have ${E}(\Phi_S(\textbf{t})) = O(h^{p+1})$. As a result of this,
we finally have to prove that ${E}^{t_0,x_0}
(\mathcal{R}_{2p+1}(t,t_0)) = O(h^{p+1})$ holds.
In the following, let $h < 1$. The autonomous version of the SRK
method~(\ref{St-srk-method-m}) can be written as
\begin{equation} \label{SRK-method-implizit-Kronecker-Darstellung}
    \begin{split}
    H^{(k,\nu)} &= \left( e \otimes I \right) Y_n + \sum_{r=0}^m
    \sum_{\mu, \iota \in \overline{\mathcal{M}}} \theta_{\iota}(h)
    \left( {B^{(\iota)}}^{(k,\nu),(r,\mu)} \otimes I \right) G_r
    \left( H^{(r,\mu)} \right) \\
    Y_{n+1} &= Y_n + \sum_{k=0}^m \sum_{\nu, \iota \in
    \overline{\mathcal{M}}} \theta_{\iota}(h) \left(
    {\gamma^{(\iota)}}^{(k,\nu)} \otimes I \right) G_k \left(
    H^{(k,\nu)} \right) .
    \end{split}
\end{equation}
Here, denote $\theta_{0}(h)=h$ and ${\gamma^{(0)}}^{(0,0)} =
\alpha$, ${\gamma^{(0)}}^{(k,\nu)}=0$ for $k \neq 0$ or $\nu \neq
0$, ${B^{(0)}}^{(k,\nu),(0,0)} = A^{(k,\nu),(0,0)}$ and
${B^{(0)}}^{(k,\nu),(r,\mu)} = 0$  for $r \neq 0$ or $\mu \neq 0$.
Further we denote $b^0 = a$, $G_k ( H^{(k,\nu)} ) = ( b^k (
H_1^{(k,\nu)} )^T, \ldots, b^k ( H_s^{(k,\nu)} )^T )^T \in
\mathbb{R}^{d \cdot s}$, $H^{(k,\nu)} = ( {H_1^{(k,\nu)}}^T,
\ldots, {H_s^{(k,\nu)}}^T )^T \in \mathbb{R}^{d \cdot s}$, $I \in
\mathbb{R}^{d \times d}$ and $e = (1, \ldots, 1)^T \in
\mathbb{R}^s$. In the following the norm $\| G_k ( H^{(k,\nu)}) \|
= \max_{1 \leq i \leq s} \| b^k ( H_i^{(k,\nu)} ) \|$ is used.
Then, with the linear growth condition $\| G_k ( H^{(k,\nu)} ) \|
\leq C_1 ( 1 + \| H^{(k,\nu)} \| )$ and
\begin{equation*}
    C_2 = \max_{\iota,k,\nu,r,\mu} \left\{ \left\|
    {B^{(\iota)}}^{(k,\nu),(r,\mu)} \otimes I \right\|, \left\|
    {\gamma^{(\iota)}}^{(k,\nu)} \otimes I \right\|, \left\|
    e \otimes I \right\| \right\}
\end{equation*}
the following inequality holds:
\begin{equation}
    \begin{split}
    \max_{(k,\nu)} \left\| H^{(k,\nu)} \right\| &\leq C_2 \| Y_n \| +
    \sum_{r=0}^m \sum_{\mu, \iota \in \overline{\mathcal{M}}}
    \left| \theta_{\iota}(h) \right| C_2 \, C_1 \left(
    1 + \left\| H^{(r,\mu)} \right\| \right) \\
    &\leq C_2 \| Y_n \| + (m+1) |\overline{\mathcal{M}}|^2 \, \max_{\iota \in
    \overline{\mathcal{M}}} \left| \theta_{\iota}(h) \right|
    \, C_1 C_2 \left( 1 + \max_{(k,\nu)} \left\| H^{(k,\nu)} \right\|
    \right) .
    \end{split}
\end{equation}
Let $C_3 = (m+1) \, |\overline{\mathcal{M}}|^2 \, C_1 C_2$. Then
for $\max_{\iota \in \overline{\mathcal{M}}} |\theta_{\iota}(h)|
\leq \tfrac{1}{2 \, C_3}$ holds
\begin{equation} \label{Abschaetzung_H_1}
    \begin{split}
    \max_{(k,\nu)} \left\| H^{(k,\nu)} \right\| &\leq \left( C_2 \| Y_n \|
    + C_3 \, \max_{\iota \in \overline{\mathcal{M}}} \left|
    \theta_{\iota}(h) \right| \right) \frac{1}{1- C_3 \,
    \max_{\iota \in \overline{\mathcal{M}}} \left|
    \theta_{\iota}(h) \right|} \\
    &\leq 2 \, C_2 \| Y_n \| + 2 \, C_3 \, \max_{\iota \in \overline{\mathcal{M}}} \left|
    \theta_{\iota}(h) \right| \\
    &\leq C_4 \, \left( 1 + \| Y_n \| \right) .
    \end{split}
\end{equation}
Next, consider the $q$th derivative. By
(\ref{Lem-Leibniz-formula-2-eqn1}) and similar considerations, we
obtain with the application of Lemma~\ref{Lem-Faa-di-Bruno-1}
using the notation~(\ref{Bezeichnung-Faa-di-Bruno-1-short}) that
\begin{equation}
    \begin{split}
    & \left\| \sum_{\nu_1, \ldots, \nu_q \in
    \overline{\mathcal{M}}} \frac{\partial^q
    {H^{(k,\nu)}}^J}{\partial \theta_{\nu_1} \ldots \partial
    \theta_{\nu_q}} \right\| \\
    & \leq q \sum_{r=0}^m
    \sum_{\mu,\nu_q \in \overline{\mathcal{M}}} \left\|
    {B^{(\nu_q)}}^{(k,\nu),(r,\mu)} \otimes I \right\| \left\|
    \sum_{\nu_1, \ldots, \nu_{q-1} \in \overline{\mathcal{M}}}
    \frac{\partial^{q-1} G_r \left( H^{(r,\mu)}
    \right)^J}{\partial \theta_{\nu_1} \ldots \partial
    \theta_{\nu_{q-1}}} \right\| \\
    & + \sum_{r=0}^m \sum_{\mu, \iota \in \overline{\mathcal{M}}}
    |\theta_{\iota}(h)| \left\| {B^{(\iota)}}^{(k,\nu),(r,\mu)} \otimes
    I \right\| \left\| \sum_{\nu_1, \ldots, \nu_q \in
    \overline{\mathcal{M}}} \frac{\partial^q
    G_r \left( H^{(r,\mu)} \right)^J}{\partial \theta_{\nu_1} \ldots \partial
    \theta_{\nu_q}} \right\| \\
    & \leq q \, |\overline{\mathcal{M}}|
    C_2 \sum_{r=0}^m \sum_{\mu \in \overline{\mathcal{M}}} \sum_{\textbf{u} \in
    SLTS_q^{(\sigma_r)}} \sum_{K_1, \ldots,
    K_{m(\textbf{u})}=1}^{d \cdot s} \left\| {G_r^J}_{K_1 \ldots
    K_{m(\textbf{u})}} \left( H^{(r,\mu)} \right) \right\| \times
    \\
    & \times \left\| \left( {H^{(r,\mu)}}^{K_1} \right)^{(\delta_1)}
    \right\| \cdot \ldots \cdot \left\| \left( {H^{(r,\mu)}}^{K_{m(\textbf{u})}}
    \right)^{(\delta_{m(\textbf{u})})} \right\| \\
    & + |\overline{\mathcal{M}}| C_2 \, \max_{\iota \in
    \overline{\mathcal{M}}} |\theta_{\iota}(h)| \sum_{r=0}^m
    \sum_{\mu \in \overline{\mathcal{M}}} \sum_{\substack{ \textbf{u} \in
    SLTS_{q+1}^{(\sigma_r)} \\ m(\textbf{u}) > 1}} \sum_{K_1,
    \ldots, K_{m(\textbf{u})}=1}^{d \cdot s} \left\| {G_r^J}_{K_1 \ldots
    K_{m(\textbf{u})}} \left( H^{(r,\mu)} \right) \right\| \times
    \\
    & \times \left\| \left( {H^{(r,\mu)}}^{K_1} \right)^{(\delta_1)}
    \right\| \cdot \ldots \cdot \left\| \left( {H^{(r,\mu)}}^{K_{m(\textbf{u})}}
    \right)^{(\delta_{m(\textbf{u})})} \right\| \\
    & + |\overline{\mathcal{M}}| C_2 \, \max_{\iota \in
    \overline{\mathcal{M}}} |\theta_{\iota}(h)| \sum_{r=0}^m
    \sum_{\mu \in \overline{\mathcal{M}}} \sum_{K_1=1}^{d \cdot s}
    \left\| {G_r^J}_{K_1} \left( H^{(r,\mu)} \right) \right\|
    \cdot \left\| \left( {H^{(r,\mu)}}^{K_1} \right)^{(q)}
    \right\| .
    \end{split}
\end{equation}
Due to the Lipschitz condition and the polynomial growth
condition, we have $\| {G_r^J}_{K_1} ( H^{(r,\mu)} ) \| \leq L$
and $\| {G_r^J}_{K_1 \ldots K_{m(\textbf{u})}} ( H^{(r,\mu)} ) \|
\leq C_5 ( 1 + ( \max_{(k,\nu)} \| H^{(k,\nu)} \| )^{2l} )$ which
is bounded by some constant $C_6$ only depending on $\| Y_n \|$
due to (\ref{Abschaetzung_H_1}). Therefore, with $C_{7} = C_2 \,
|\overline{\mathcal{M}}|^2 \, (m+1)$ follows
\begin{equation}
    \begin{split}
    \max_{J,(k,\nu)} & \left\| \left( {H^{(k,\nu)}}^J \right)^{(q)}
    \right\| \leq q \, C_7 \sum_{\textbf{u} \in
    SLTS_{q}^{(\sigma)}} (d \cdot s)^{m(\textbf{u})} \, C_6
    \prod_{i=1}^{m(\textbf{u})} \max_{J,(k,\nu)} \left\| \left( {H^{(k,\nu)}}^J \right)^{(\delta_i)}
    \right\| \\
    & + C_7 \, \max_{\iota \in \overline{\mathcal{M}}} |\theta_{\iota}(h)| \sum_{\substack{ \textbf{u} \in
    SLTS_{q+1}^{(\sigma)} \\ m(\textbf{u}) > 1 }} (d \cdot s)^{m(\textbf{u})} \,
    C_6 \prod_{i=1}^{m(\textbf{u})} \max_{J,(k,\nu)} \left\| \left( {H^{(k,\nu)}}^J \right)^{(\delta_i)}
    \right\| \\
    & + C_7 \, \max_{\iota \in \overline{\mathcal{M}}} |\theta_{\iota}(h)| \,
    (d \cdot s) \, L \max_{J,(k,\nu)} \left\| \left( {H^{(k,\nu)}}^J \right)^{(q)}
    \right\| .
    \end{split}
\end{equation}
Let $C_8 = C_7 \, d \, s \, L$. Then we get with $\max_{\iota \in
\overline{\mathcal{M}}} |\theta_{\iota}(h)| \leq \tfrac{1}{2 \,
C_8}$ that
\begin{equation} \label{Abschaetzung_qte_Ableitung_H_1}
    \begin{split}
    \max_{J,(k,\nu)} & \left\| \left( {H^{(k,\nu)}}^J \right)^{(q)}
    \right\| \leq 2 \, q \, C_7 \sum_{\textbf{u} \in
    SLTS_{q}^{(\sigma)}} (d \cdot s)^{m(\textbf{u})} \, C_6
    \prod_{i=1}^{m(\textbf{u})} \max_{J,(k,\nu)} \left\| \left( {H^{(k,\nu)}}^J \right)^{(\delta_i)}
    \right\| \\
    & + 2 \, C_7 \, \max_{\iota \in \overline{\mathcal{M}}} |\theta_{\iota}(h)| \sum_{\substack{ \textbf{u} \in
    SLTS_{q+1}^{(\sigma)} \\ m(\textbf{u}) > 1 }} (d \cdot s)^{m(\textbf{u})} \,
    C_6 \prod_{i=1}^{m(\textbf{u})} \max_{J,(k,\nu)} \left\| \left( {H^{(k,\nu)}}^J \right)^{(\delta_i)}
    \right\|
    \end{split}
\end{equation}
holds with $\delta_i=\delta_i(\textbf{u}) \leq q-1$ because
$m(\textbf{u})>1$. Especially for $q=1$ where due to the linear
growth condition $C_6 = C_9 \left( 1 + \| Y_n \| \right)$, we
arrive at
\begin{equation} \label{Abschaetzung_1.Ableitung_H_1}
    \begin{split}
    \max_{J,(k,\nu)} \left\| \sum_{\nu_1 \in
    \overline{\mathcal{M}}} \frac{\partial
    {H^{(k,\nu)}}^J}{\partial \theta_{\nu_1}(h)} \right\| &\leq
    C_{10} \left( 1 + \| Y_n \| \right).
    \end{split}
\end{equation}
Applying formula~(\ref{Abschaetzung_qte_Ableitung_H_1})
recursively and using finally (\ref{Abschaetzung_H_1}) yields an
upper bound $C_q(Y_n)$ of the $q$th derivative of $H^{(k,\nu)}$
only depending on $\| Y_n \|$ for all $q \in \mathbb{N}$. Due to
the definition of $C_2$ and the same structure of
$Y_{n+1}=A(t_n,Y_n,\theta(h))$ as $H^{(k,\nu)}$, the same upper
bound holds also for the $q$th derivative of $A(t_n,Y_n,
\theta(h))$. Since $f \in C_P^{2p+2} (\mathbb{R}^d,\mathbb{R})$,
we obtain for $\xi \in \, ]0,1[$ and $|\theta_{\iota}(h)| \leq
\sqrt{h} \, C_{\vartheta}$ with the Jensen inequality
\begin{equation}
    \begin{split}
    & \left\| {E}^{t_0,x_0} \left( \sum_{\nu_1,
    \ldots, \nu_{2p+2} \in \overline{\mathcal{M}}}
    \Delta \theta_{\nu_1} \cdot \ldots
    \cdot \Delta \theta_{\nu_{2p+2}} \frac{\partial^{2p+2}
    f(A(t_0,Y(t_0), \xi \, \theta(h)))}{\partial \theta_{\nu_1} \ldots
    \partial \theta_{\nu_{2p+2}}} \right) \right\| \\
    & \leq {E}^{t_0,x_0} \left( \left(\max_{\iota \in
    \overline{\mathcal{M}}} |\theta_{\iota}(h)| \right)^{2p+2}
    \times \right. \\
    & \left. \times \sum_{\textbf{u}
    \in SLTS_{2p+3}^{(\gamma)}} \sum_{K_1, \ldots, K_m=1}^d \left\|
    f_{K_1 \ldots K_m} \left( A(t_0,Y(t_0), \xi \, \theta(h)) \right)
    \right\| \prod_{i=1}^{m(\textbf{u})} C_{\delta_i}(Y(t_0))
    \right) \\
    & \leq h^{p+1} C_{\vartheta}^{2p+2} \sum_{\textbf{u}
    \in SLTS_{2p+3}^{(\gamma)}} d^{m(\textbf{u})} C_f \left( 1 + \left( C_4 (1 +
    \| x_0 \|) \right)^{2r(\textbf{u})} \right) \prod_{i=1}^{m(\textbf{u})}
    C_{\delta_i}(x_0)
    \end{split}
\end{equation}
and it follows $E^{t_0,x_0}(\mathcal{R}_{2p+1}(t,t_0))
= O(h^{p+1})$. \hfill $\Box$ \\ \\
The result of Proposition~\ref{St-SRK-expansion-expectation-1:Wm}
can also be proved for general unbounded random variables in the
case of explicit SRK methods (see \cite{Roe03},
Proposition~2.6.1). However, especially for weak approximations it
is usual to use bounded random variables which are often easier to
generate (see, e.g., \cite{KP99}, \cite{Mil95}, \cite{Ta90}). \\ \\
The approximation $Y$ has to be uniformly bounded with respect to
the number $N$ of steps in order to guarantee convergence.
Therefore, sufficient conditions for the random variables and for
some coefficients of the stochastic Runge-Kutta method such that
$Y$ is uniformly bounded are calculated.
\begin{Sat} \label{St-general-order-cond-bound-lemma1}
    Let $a^i, b^{i,j} \in C^1(\mathbb{R}^d, \mathbb{R})$ satisfy
    a Lipschitz and a linear growth condition
    and let for all $1 \leq k \leq m$ and $\nu \in \mathcal{M}$
    \begin{equation} \label{St-general-order-cond-bound-lemma1-eqn1}
        E \left( {z^{(k,\nu)}}^T e \right)=0 .
    \end{equation}
    Further assume that each random variable can be expressed
    as $\theta_{\iota}(h) = \sqrt{h} \cdot \vartheta_{\iota}$ for $\iota \in \mathcal{M}$
    with a bounded random variable $\vartheta_{\iota}$.
    Then the approximation $Y$ by the stochastic
    Runge-Kutta method~(\ref{St-srk-method-m}) has
    uniformly bounded moments, i.e.\ for $r \in \mathbb{N}$
    the expectation $E( \| Y_n \|^{2r} )$ is uniformly bounded w.r.t.\ the number of steps $N$
    for all $n=0,1, \ldots, N$.
\end{Sat}
{\bf{Proof.}} Let $h < 1$. Using the
notation~(\ref{SRK-method-implizit-Kronecker-Darstellung}) we get
with the linear growth condition and with (\ref{Abschaetzung_H_1})
\begin{equation} \label{boundedness-estimate1}
    \begin{split}
    \left\| A(t_n,Y_n,\theta(h)) - Y_n \right\| & \leq \sum_{k=0}^m
    \sum_{\nu,\iota \in \overline{\mathcal{M}}}
    |\theta_{\iota}(h)| \, C_2 \, \left\| G_k \left( H^{(k,\nu)}
    \right) \right\| \\
    & \leq (m+1) \, |\overline{\mathcal{M}}|^2 \,
    \max_{\iota \in \overline{\mathcal{M}}}
    |\theta_{\iota}(h)| \, C_2 \, C_9 \left( 1 + \max_{(k,\nu)}
    \left\| H^{(k,\nu)} \right\| \right) \\
    & \leq C_{11} \left( 1 + \| Y_n \| \right) \, \sqrt{h} .
    \end{split}
\end{equation}
Next, we get with one step of the Taylor-expansion of $G_k$ for
$\xi \in \, ]0,1[$ that
\begin{equation}
    \begin{split}
    & \left\| E \left( A(t_n,Y_n,\theta(h)) - Y_n \right)
    \right\| \leq \left\| \sum_{k=1}^m \sum_{\nu,\iota \in
    \mathcal{M}} E \left( \theta_{\iota}(h)
    \right) \left( {\gamma^{(\iota)}}^{(k,\nu)} \otimes I \right)
    G_k \left( \left( e \otimes I \right) Y_n \right) \right\| \\
    & + h \, \left\| {\gamma^{(0)}}^{(0,0)} \otimes I \right\| \,
    \left\| G_0 \left( \left( e \otimes I \right) Y_n \right) \right\|
    + \left\| E \left( \sum_{k=0}^m \sum_{\nu,\iota \in
    \overline{\mathcal{M}}} \theta_{\iota}(h) \left(
    {\gamma^{(\iota)}}^{(k,\nu)} \otimes I \right) \right. \right. \times \\
    & \times \left. \left. \sum_{\mu \in
    \overline{\mathcal{M}}} \sum_{J=1}^{d \cdot s} \frac{\partial
    G_k \left( H^{(k,\nu)}\left(\xi \theta(h) \right)
    \right)}{\partial x^J} \frac{\partial H^{(k,\nu)}\left(\xi
    \theta(h) \right)^J}{\partial \theta_{\mu}} \, \Delta
    \theta_{\mu}(h) \right) \right\| .
    \end{split}
\end{equation}
The first summand on the right hand side vanishes due to
(\ref{St-general-order-cond-bound-lemma1-eqn1}). With a
Lip\-schitz constant $L$ for $G$ and the linear growth condition,
we get with the Jensen inequality
\begin{equation}
    \begin{split}
    & \left\| E \left( A(t_n,Y_n,\theta(h)) - Y_n \right)
    \right\| \leq h \, C_2 \, C_9 \left( 1 + \left\| ( e \otimes I
    ) Y_n \right\| \right) \\
    & + E \left( \sum_{k=0}^m
    \sum_{\nu,\iota \in \overline{\mathcal{M}}} \left\|
    {\gamma^{(\iota)}}^{(k,\nu)} \otimes I \right\| \sum_{J=1}^{d
    \cdot s} L \left\| \sum_{\mu \in \overline{\mathcal{M}}}
    \frac{\partial H^{(k,\nu)}\left(\xi \theta(h)
    \right)^J}{\partial \theta_{\mu}} \right\| \, \left(
    \max_{\iota \in \overline{\mathcal{M}}} |\theta_{\iota}(h)|
    \right)^2 \right) .
    \end{split}
\end{equation}
Finally, applying (\ref{Abschaetzung_1.Ableitung_H_1}) and the
condition $|\theta_{\iota}(h)| \leq \sqrt{h} \, C_{\vartheta}$, we
get
\begin{equation} \label{boundedness-estimate2}
    \begin{split}
    & \left\| E \left( A(t_n,Y_n,\theta(h)) - Y_n \right)
    \right\| \leq C_{12} \left( 1 + \| Y_n \| \right) \, h
    \end{split}
\end{equation}
Now, Lemma~\ref{St-lg-Lem-bound-1}
can be applied because
(\ref{boundedness-estimate1}) and (\ref{boundedness-estimate2})
are fulfilled. This yields the existence of $E(\|Y_n\|^{2r})$ for
all $r \in \mathbb{N}$ and provides that the moments are uniformly
bounded with respect to
$N$ and $n=1, \ldots, N$. \hfill $\Box$ \\ \\
The next step is to compare the representations of the solution of
the stochastic differential equation in
Theorem~\ref{St-tree-expansion-exact-sol:Wm} with the
representation of the approximation in
Proposition~\ref{St-SRK-expansion-expectation-1:Wm}. According to
Theorem~\ref{St-lg-theorem-main} these representations have to
coincide up to order $p+1$ locally. This leads to different
conditions w.r.t.\ trees in $TS(I)$ and $TS(S)$ on the one hand
and trees in $TS(\Delta) \setminus TS(I)$ and $TS(\Delta)
\setminus TS(S)$ on the other hand. Having in mind that for
$\textbf{t} \in TS(I)$ or $\textbf{t} \in TS(S)$ we have
$s(\textbf{t})/2$ different variable indices while for the same
tree $\textbf{t} \in TS(\Delta)$ we have twice as much, i.e.\
$s(\textbf{t})$ different variable indices, we use the following
helpful definition.
\begin{Def}
    Let $|\textbf{t}|$ denote the tree which is obtained if the nodes
    $\sigma_{j_i}$ of $\textbf{t}$ are replaced by $\sigma$, i.e.\
    by omitting all variable indices.
    Let a tree $\textbf{t} \in TS(*)$ for
    $* \in \{I,S\}$ with variable indices $j_1, \ldots, j_{s(\textbf{t})/2}$ be given
    and let $\textbf{u} \in TS(\Delta)$
    with variable indices $\hat{j}_1, \ldots, \hat{j}_{s(\textbf{u})}$
    denote the tree which is equivalent to $\textbf{t}$ except for the
    variable indices, i.e.\ $|\textbf{t}| \sim |\textbf{u}|$ with
    $s(\textbf{t})=s(\textbf{u})$.
    For a fixed choice of correlations of type $j_k=j_l$ or $j_k \neq
    j_l$, $1 \leq k < l \leq s(\textbf{t})/2$,
    between the indices $j_1,
    \ldots, j_{s(\textbf{t})/2}$, let $\beta(\textbf{t})$ denote the
    number of all possible correlations between the
    indices $\hat{j}_1, \ldots, \hat{j}_{s(\textbf{u})}$ of tree
    $\textbf{u}$ such that $\textbf{t} \sim \textbf{u}$ holds. In
    the case of $s(\textbf{t})=0$ or $\textbf{t} \in TS(\Delta) \setminus TS(*)$,
    $* \in \{I,S\}$, define $\beta(\textbf{t})=1$.
\end{Def}
Note that in case of $m=1$ we have $\beta(\textbf{t})=1$ for all
$\textbf{t} \in TS(*)$, $* \in \{I,S\}$. As an example consider
the trees $\textbf{t} = (\sigma_{j_1}, \sigma_{j_1}, \sigma_{j_2},
\sigma_{j_2}) \in TS(I)$ and $\textbf{u} = (\sigma_{\hat{j}_1},
\sigma_{\hat{j}_2}, \sigma_{\hat{j}_3}, \sigma_{\hat{j}_4}) \in
TS(\Delta)$. For the correlation $j_1 = j_2$ of $\textbf{t}$ we
have exactly one possibility for the choice of a correlation of
$\textbf{u}$: We have to choose $\hat{j}_1 = \hat{j}_2 = \hat{j}_3
= \hat{j}_4$, i.e.\ in this case we have $\beta(\textbf{t}) = 1$.
However, in case of the correlation $j_1 \neq j_2$ for
$\textbf{t}$, there are three different possible correlations for
$\textbf{u}$: We can choose $\hat{j}_1 = \hat{j}_2 \neq \hat{j}_3
= \hat{j}_4$, $\hat{j}_1 = \hat{j}_3 \neq \hat{j}_2 = \hat{j}_4$
or $\hat{j}_1 = \hat{j}_4 \neq \hat{j}_2 = \hat{j}_3$, thus we
have $\beta(\textbf{t}) = 3$. As a second example, for the trees
$\textbf{t} = (\sigma_{j_1}, \sigma_{j_2}, \{ \sigma_{j_2}
\}_{j_1} ) \in TS(I)$ and $\textbf{u} = ( \sigma_{\hat{j}_1},
\sigma_{\hat{j}_2}, \{ \sigma_{\hat{j}_4} \}_{\hat{j}_3} ) \in
TS(\Delta)$, two different correlations are distinguished. On the
one hand we have the correlation $j_1=j_2$ for $\textbf{t}$ where
we get the only possible correlation $\hat{j}_1 = \hat{j}_2 =
\hat{j}_3 = \hat{j}_4$ for $\textbf{u}$, i.e.\
$\beta(\textbf{t})=1$. On the other hand we have $j_1 \neq j_2$ as
a correlation for $\textbf{t}$ allowing us two different
correlations $\hat{j}_1 = \hat{j}_3 \neq \hat{j}_2 = \hat{j}_4$
and $\hat{j}_2 = \hat{j}_3 \neq \hat{j}_1 = \hat{j}_4$ for
$\textbf{u}$. Thus we get $\beta(\textbf{t})=2$
in the latter case. \\ \\
The main theorem for stochastic Runge-Kutta methods of
type~(\ref{St-srk-method-m}) yields general conditions for the
coefficients and the random variables of the method such that
convergence with some order $p$ in the weak sense is assured. Note
that for every tree $\textbf{t} \in TS(*)$ with variable indices
$j_1, \ldots, j_{s(\textbf{t})/2}$ there exists a tree $\textbf{u}
\in TS(\Delta)$ with $|\textbf{u}| \sim |\textbf{t}|$ and variable
indices $\hat{j}_1, \ldots, \hat{j}_{s(\textbf{u})}$ such that for
some suitable correlation of type $\hat{j}_{k} = \hat{j}_{l}$ or
$\hat{j}_{k} \neq \hat{j}_{l}$, $1 \leq k < l \leq s(\textbf{u})$,
we have $\textbf{t} \sim \textbf{u}$ and thus $\textbf{u} \in
TS(*)$ with $\alpha_*(\textbf{u})=\alpha_*(\textbf{t})$ for $* \in
\{I,S\}$. However, we have $\alpha_*(\textbf{u})=0$ for all
$\textbf{u} \in TS(\Delta) \setminus TS(*)$ for $* \in \{I,S\}$.
\begin{The} \label{Ito-St-Theo-conv-cond-tree-main:Wm}
Let $X$ be the solution of either an It{\^o} or a Stratonovich
stochastic differential
equation~(\ref{Intro-Ito-St-SDE1-integralform-Wm}) considered
w.r.t.\ an $m$-dimensional Wiener process, and with $f, a^i,
\tilde{a}^i, b^{i,j} \in C_P^{2(p+1)}(\mathbb{R}^d, \mathbb{R})$
for $i=1, \ldots, d$ and $j=1, \ldots, m$. Then the
approximation $Y$ by the stochastic Runge-Kutta
method~(\ref{St-srk-method-m}) with maximum step size $h$ is of
weak order $p$, if for all $\textbf{t} \in TS(\Delta)$ with
$\rho(\textbf{t}) \leq p + \tfrac{1}{2}$ and all correlations of
type $j_k=j_l$ or $j_k \neq j_l$, $1 \leq k < l \leq
s(\textbf{t})$, between the indices $j_1, \ldots,
j_{s(\textbf{t})} \in \{1, \ldots, m\}$ of $\textbf{t}$ the
equations
\begin{equation} \label{Ito-St-Theo-conv-cond-tree-main-eqn1:Wm}
    \frac{\alpha_{*}(\textbf{t}) \cdot h^{\rho(\textbf{t})}}{2^{s(\textbf{t})/2}
    \cdot \rho(\textbf{t})!} =
    \frac{\alpha_{\Delta}(\textbf{t}) \cdot \beta(\textbf{t})
    \cdot \gamma(\textbf{t})
    \cdot E(\Phi_S(\textbf{t}))}{(l(\textbf{t})-1)!}
\end{equation}
hold for $* = I$ in case of It{\^o} SDEs and $* = S$ in case of
Stratonovich SDEs, provided that (\ref{SRK-weights-condition}) and
(\ref{St-SRK-moment-condition-1:Wm}) hold and that the
approximation $Y$ has uniformly bounded moments w.r.t.\ the number
$N$ of steps.
\end{The}
{\bf Proof.} Apply Theorem~\ref{St-lg-theorem-main} and compare
the coefficients from the representations of the solution in
Theorem~\ref{St-tree-expansion-exact-sol:Wm} with the coefficients
of the SRK method in
Proposition~\ref{St-SRK-expansion-expectation-1:Wm}, where $TS(*)
\subseteq TS(\Delta)$, $* \in \{I, S\}$. Finally, we take into
account the summation w.r.t.\ variable indices. Therefore, the
correlation index $\beta(\textbf{t})$ has to be added and we
obtain the conditions
(\ref{Ito-St-Theo-conv-cond-tree-main-eqn1:Wm}). \hfill $\Box$
\begin{Bem}
    Theorem~\ref{Ito-St-Theo-conv-cond-tree-main:Wm}
    provides uniform weak convergence with order $p$ on the interval $\mathcal{I}=[t_0,T]$
    for the stochastic Runge-Kutta method in the case of a non-random time discretization
    $\mathcal{I}_h$. That is for each $f \in C_P^{2(p+1)}(\mathbb{R}^d, \mathbb{R})$
    there exists a finite constant $C_f$ not depending on the
    maximum step size $h$ such that
    \begin{equation}
        \max_{0 \leq k \leq N} \big| E(f(X_{t_k})) - E(f(Y_k))
        \big| \leq C_f \, h^p
    \end{equation}
    holds. This is a consequence of Theorem~\ref{St-lg-theorem-main}
    (see, e.g., \cite{KP99}, \cite{Mil95}).
\end{Bem}
Table~\ref{tabelle1} contains all S-trees of $TS(I)$ and $TS(S)$
up to order two with the corresponding cardinalities $\alpha_I$
and $\alpha_S$. Table~\ref{tabelle2} contains all S-trees of
$TS(\Delta)$ up to order 2.5 with the values of $\alpha_{\Delta}$.
The cardinalities can be determined very easily as the number of
possibilities to build up the considered tree due to the
corresponding rules of growth. Together with Table~\ref{tabelle1b}
containing the values of $\beta$, we can consider the following
example:
\begin{Bsp} Assume that $m \geq 1$.
\begin{enumerate}[a)]
    \item As a first example, let us have a look at tree $\textbf{t}_{2.5} =
    (\sigma_{j_1},[\sigma_{j_2}]) \in TS(\Delta)$ with
    parameters $l(\textbf{t}_{2.5})=4$, $\gamma(\textbf{t}_{2.5})=2$,
    $s(\textbf{t}_{2.5})=2$, $\alpha_{\Delta}(\textbf{t}_{2.5})=3$
    and $\rho(\textbf{t}_{2.5})=2$. Then the following
    correlations have to be distinguished: For $j_1=j_2$ follows
    that $\textbf{t}_{2.5} \in TS(*)$ with
    $\alpha_I(\textbf{t}_{2.5})=\alpha_S(\textbf{t}_{2.5})=2$ and
    $\beta(\textbf{t}_{2.5})=1$. Then for $j_1 \in \{1, \ldots,
    m\}$ Theorem~\ref{Ito-St-Theo-conv-cond-tree-main:Wm} yields
    the conditions
    \[
        E \left( \left( \sum_{\nu \in \mathcal{M}}
        {z^{(j_1,\nu)}}^T e \right) \left( {z^{(0,0)}}^T \left(
        \sum_{\mu \in \mathcal{M}} Z^{(0,0),(j_1,\mu)} e \right)
        \right) \right) = \frac{2 \cdot 3! \cdot h^2}{2^1 \cdot 2!
        \cdot 3 \cdot 1 \cdot 2}.
    \]
    Here, the conditions for It{\^o} and Stratonovich calculus
    coincide. However, for $j_1 \neq j_2$ follows $\textbf{t}_{2.5} \notin
    TS(*)$, i.e.\
    $\alpha_I(\textbf{t}_{2.5})=\alpha_S(\textbf{t}_{2.5})=0$, and
    thus one gets for $j_1,j_2 \in \{1, \ldots, m\}$ with $j_1 \neq j_2$
    the additional conditions
    \[
        E \left( \left( \sum_{\nu \in \mathcal{M}}
        {z^{(j_1,\nu)}}^T e \right) \left( {z^{(0,0)}}^T \left(
        \sum_{\mu \in \mathcal{M}} Z^{(0,0),(j_2,\mu)} e \right)
        \right) \right) = 0 .
    \]
    \item Consider $\textbf{t}_{2.11}= (\sigma_{j_1}, \sigma_{j_2},
    \sigma_{j_3}, \sigma_{j_4}) \in TS(\Delta)$ with
    $l(\textbf{t}_{2.11})=5$, $\gamma(\textbf{t}_{2.11})=1$, $s(\textbf{t}_{2.11})=4$,
    $\alpha_{\Delta}(\textbf{t}_{2.11})=1$ and $\rho(\textbf{t}_{2.11})=2$. The following
    correlations have to be analyzed:
    For $j_1=j_3 \neq j_2=j_4$ we have $\textbf{t}_{2.11}
    \in TS(*)$ with $\alpha_I(\textbf{t}_{2.11})= \alpha_S(\textbf{t}_{2.11})=1$
    and $\beta(\textbf{t}_{2.11})=3$.
    Then Theorem~\ref{Ito-St-Theo-conv-cond-tree-main:Wm} yields the
    condition
    \begin{equation*}
        \begin{split}
        E \left( \left(\sum_{\nu \in \mathcal{M}} {z^{(j_1,\nu)}}^T e \right)^2
        \left(\sum_{\nu \in \mathcal{M}} {z^{(j_2,\nu)}}^T e \right)^2
        \right) = \frac{4! \cdot h^2}{2^2 \cdot 2! \cdot 3}
        \end{split}
    \end{equation*}
    with $j_1,j_2 \in \{1,\ldots, m\}$, $j_1 \neq j_2$, for both
    It{\^o} and Stratonovich calculus.
    For $j_1=j_2=j_3=j_4$ we have $\textbf{t}_{2.11} \in TS(*)$ with
    $\alpha_I(\textbf{t}_{2.11})=\alpha_S(\textbf{t}_{2.11})=1$ and
    $\beta(\textbf{t}_{2.11})=1$. Again,
    Theorem~\ref{Ito-St-Theo-conv-cond-tree-main:Wm} yields the
    condition
    \begin{equation*}
        \begin{split}
        E \left( \left(\sum_{\nu \in \mathcal{M}} {z^{(j_1,\nu)}}^T e
        \right)^4 \right) = \frac{4! \cdot h^2}{2^2 \cdot 2! \cdot 1}
        \end{split}
    \end{equation*}
    with $j_1 \in \{1,\ldots,m\}$ for both It{\^o} and Stratonovich calculus.
    For all remaining correlations of the indices follows that $\textbf{t}_{2.11} \notin
    TS(*)$ and thus
    $\alpha_I(\textbf{t}_{2.11})=\alpha_S(\textbf{t}_{2.11})=0$.
    Therefore, the condition $E(\Phi_S(\textbf{t}_{2.11}))=0$ has to be
    fulfilled for the remaining correlations.
    \item For $\textbf{t}_{2.12} = (\sigma_{j_1}, \sigma_{j_2},
    \{ \sigma_{j_4}\}_{j_3})$ with $l(\textbf{t}_{2.12})=5$,
    $\gamma(\textbf{t}_{2.12})=2$,
    $s(\textbf{t}_{2.12})=4$, $\alpha_{\Delta}(\textbf{t}_{2.12})=6$ and
    $\rho(\textbf{t}_{2.12})=2$, consider the
    following correlations: For $j_1=j_2 \neq j_3=j_4$
    we have $\textbf{t}_{2.12}=\textbf{t}_{2.12a} \in TS(S)$ with
    $\alpha_S(\textbf{t}_{2.12a})=2$ and
    $\beta(\textbf{t}_{2.12a})=1$. Therefore, we get the condition
    \begin{equation*}
        \begin{split}
        E \left( \left( \sum_{\nu \in \mathcal{M}} {z^{(j_1,\nu)}}^T e
        \right)^2
        \left( \sum_{\nu \in \mathcal{M}} {z^{(j_3,\nu)}}^T \left(
        \sum_{\mu \in \mathcal{M}} Z^{(j_3,\nu)(j_3,\mu)} e \right)
        \right) \right) = \frac{2 \cdot 4! \cdot h^2}{2^2 \cdot 2!
        \cdot 6 \cdot 2}
        \end{split}
    \end{equation*}
    with $j_1,j_3 \in \{1, \ldots, m\}$, $j_1 \neq j_3$, for
    Stratonovich calculus. However, since $\textbf{t}_{2.12a}
    \notin TS(I)$ we get for It{\^o} calculus the condition
    $E(\Phi_S(\textbf{t}_{2.12a}))=0$ since
    $\alpha_I(\textbf{t}_{2.12a})=0$.
    For $j_1=j_3 \neq j_2=j_4$ or
    $j_2=j_3 \neq j_1=j_4$ we have $\textbf{t}_{2.12}=\textbf{t}_{2.12b}
    \in TS(*)$ with $\alpha_I(t_{2.12b})=\alpha_S(\textbf{t}_{2.12b})=4$
    and $\beta(t_{2.12b})=2$. Here we get the condition
    \begin{equation}
        \begin{split}
        & E \left( \left( \sum_{\nu \in \mathcal{M}} {z^{(j_1,\nu)}}^T e \right)
        \left( \sum_{\nu \in \mathcal{M}} {z^{(j_2,\nu)}}^T e \right)
        \left( \sum_{\nu \in \mathcal{M}} {z^{(j_1,\nu)}}^T \left(
        \sum_{\mu \in \mathcal{M}} Z^{(j_1,\nu)(j_2,\mu)} e \right)
        \right) \right) \\
        & = \frac{4 \cdot 4! \cdot h^2}{2^2 \cdot 2! \cdot 6 \cdot
        2 \cdot 2}
        \end{split}
    \end{equation}
    with $j_1,j_2 \in \{1, \ldots, m\}$, $j_1 \neq j_2$, for
    It{\^o} and Stratonovich calculus. Further, for
    $j_1=j_2=j_3=j_4$ we have $\textbf{t}_{2.12} \in TS(*)$ with
    $\alpha_I(\textbf{t}_{2.12})=0+4$,
    $\alpha_S(\textbf{t}_{2.12})=2+4$ and
    $\beta(\textbf{t}_{2.12})=1$. Therefore, we get the conditions
    \begin{equation*}
        \begin{split}
        E \left( \left( \sum_{\nu \in \mathcal{M}} {z^{(j_1,\nu)}}^T e
        \right)^2
        \left( \sum_{\nu \in \mathcal{M}} {z^{(j_1,\nu)}}^T \left(
        \sum_{\mu \in \mathcal{M}} Z^{(j_1,\nu)(j_1,\mu)} e \right)
        \right) \right) = \frac{\alpha_*(\textbf{t}_{2.12}) \, 4!
        \, h^2}{2^2 \cdot 2! \cdot 6 \cdot 2}
        \end{split}
    \end{equation*}
    with $j_1 \in \{1, \ldots, m\}$. Thus we have different
    conditions for It{\^o} and Stra\-to\-no\-vich calculus. Finally, for
    all remaining correlations the conditions
    $E(\Phi_S(\textbf{t}_{2.12}))=0$ have to hold due to $\textbf{t}_{2.12}
    \notin TS(*)$ in these cases.
\end{enumerate}
\end{Bsp}
\section{Conclusions}
\label{Sec:conclusions}
The present paper introduces a very general class of stochastic
Runge-Kutta methods for the approximation of stochastic
differential equations. Explicit as well as implicit SRK methods
for non-autonomous SDE systems w.r.t.\ to a multi-dimensional
Wiener process are considered. A rigorous analysis of the weak
convergence for the SRK method is given. Therefore, colored rooted
trees are introduced and an expansion of the solution and of the
approximation process is given. Finally, a theorem giving directly
the order conditions for arbitrary high order of convergence is
proved. The main advantages of the rooted tree analysis are as
follows: The required colored rooted trees can be easily
determined. So in contrast to the usual direct comparison of the
Taylor expansions, one does not need to calculate the derivatives
of $f$, $a$ and $b$. It has to be pointed out that the calculated
order conditions depend on the coefficients and the random
variables of the SRK method. Therefore, the order conditions can
also be used for the determination of suitable random variables
for the SRK method. In order to get a closed theory, the presented
results cover SRK methods for the approximation of both It{\^o}
and Stratonovich SDE systems. Finally, the presented colored
rooted tree theory and the introduced SRK methods generalize the
well known theory for deterministic Runge-Kutta methods due to
Butcher~\cite{Butcher87}. In the case of $b \equiv 0$ and
$f(x)=x$, i.e.\ an deterministic ordinary differential equation,
the SRK method coincides with a deterministic Runge-Kutta method
and also the order conditions coincide with the deterministic
order conditions. For some examples of SRK methods and the
corresponding analysis of order conditions with colored rooted
trees, we refer to R{\"o}{\ss}ler~\cite{Roe03}.
\begin{ack}
    The author is very grateful to the referee for his comments
    and fruitful suggestions.
\end{ack}
\appendix
\section{Tables}
\label{Sec:Tables}
\begin{longtable}[c]{|c|c|c|c|c||c|c|c|c|c|}
\caption{Trees $\textbf{t} \in TS(*)$, $* \in \{I,S\}$, of order
$\rho(\textbf{t}) \leq 2.5$ with variable indices $j_1, j_2 \in
\{1,\ldots,m\}$.} \label{tabelle1} \\
    \hline
    $\textbf{t}$ & tree & $\alpha_I$ & $\alpha_S$ & $\rho$ &
    $\textbf{t}$ & tree & $\alpha_I$ & $\alpha_S$ & $\rho$ \\
    \hline
    \hline
    \endfirsthead
    \hline
    $\textbf{t}$ & tree & $\alpha_I$ & $\alpha_S$ & $\rho$ &
    $\textbf{t}$ & tree & $\alpha_I$ & $\alpha_S$ & $\rho$ \\
    \hline
    \hline
    \endhead
%
    \hline
    \endfoot
%
    \hline
    \endlastfoot
%
    $\textbf{t}_{0.1}$ & $\gamma$ & 1 & 1 & 0 & $\textbf{t}_{2.11}$ & $(\sigma_{j_1},\sigma_{j_1},\sigma_{j_2},\sigma_{j_2})$ & 1 & 1 & 2 \\
    \cline{1-5}
    $\textbf{t}_{1.1}$ & $(\tau)$ & 1 & 1 & 1 & $\textbf{t}_{2.12a}$ & $(\sigma_{j_1},\sigma_{j_1},\{\sigma_{j_2}\}_{j_2})$ & 0 & 2 & 2\\
    $\textbf{t}_{1.2}$ & $(\sigma_{j_1},\sigma_{j_1})$ & 1 & 1 & 1 & $\textbf{t}_{2.12b}$ & $(\sigma_{j_1},\sigma_{j_2},\{\sigma_{j_2}\}_{j_1})$ & 4 & 4 & 2 \\
    $\textbf{t}_{1.3}$ & $(\{\sigma_{j_1}\}_{j_1})$ & 0 & 1 & 1 & $\textbf{t}_{2.13a}$ & $(\sigma_{j_1},\{\sigma_{j_2},\sigma_{j_2}\}_{j_1})$ & 2 & 2 & 2 \\
    \cline{1-5}
    $\textbf{t}_{2.1}$ & $([\tau])$ & 1 & 1 & 2 & $\textbf{t}_{2.13b}$ & $(\sigma_{j_2},\{\sigma_{j_2},\sigma_{j_1}\}_{j_1})$ & 0 & 2 & 2 \\
    $\textbf{t}_{2.2}$ & $(\tau,\tau)$ & 1 & 1 & 2 & $\textbf{t}_{2.14a}$ & $(\sigma_{j_1},\{\{\sigma_{j_2}\}_{j_2}\}_{j_1})$ & 0 & 2 & 2 \\
    $\textbf{t}_{2.3}$ & $([\{\sigma_{j_1}\}_{j_1}])$ & 0 & 1 & 2 & $\textbf{t}_{2.14b}$ & $(\sigma_{j_2},\{\{\sigma_{j_2}\}_{j_1}\}_{j_1})$ & 0 & 2 & 2 \\
    $\textbf{t}_{2.4}$ & $([\sigma_{j_1},\sigma_{j_1}])$ & 1 & 1 & 2 & $\textbf{t}_{2.15a}$ & $(\{\sigma_{j_1}\}_{j_1},\{\sigma_{j_2}\}_{j_2})$ & 0 & 1 & 2 \\
    $\textbf{t}_{2.5}$ & $(\sigma_{j_1},[\sigma_{j_1}])$ & 2 & 2 & 2 & $\textbf{t}_{2.15b}$ & $(\{\sigma_{j_2}\}_{j_1},\{\sigma_{j_2}\}_{j_1})$ & 2 & 2 & 2 \\
    $\textbf{t}_{2.6}$ & $(\{\sigma_{j_1}\}_{j_1},\tau)$ & 0 & 2 & 2 & $\textbf{t}_{2.16}$ & $(\{\sigma_{j_1},\sigma_{j_2},\sigma_{j_2}\}_{j_1})$ & 0 & 1 & 2 \\
    $\textbf{t}_{2.7}$ & $(\sigma_{j_1},\sigma_{j_1},\tau)$ & 2 & 2 & 2 & $\textbf{t}_{2.17a}$ & $(\{\sigma_{j_1},\{\sigma_{j_2}\}_{j_2}\}_{j_1})$ & 0 & 1 & 2 \\
    $\textbf{t}_{2.8}$ & $(\sigma_{j_1},\{\tau\}_{j_1})$ & 2 & 2 & 2 & $\textbf{t}_{2.17b}$ & $(\{\sigma_{j_2},\{\sigma_{j_2}\}_{j_1}\}_{j_1})$ & 0 & 2 & 2 \\
    $\textbf{t}_{2.9}$ & $(\{\{\tau\}_{j_1}\}_{j_1})$ & 0 & 1 & 2 & $\textbf{t}_{2.18}$ & $(\{\{\sigma_{j_2},\sigma_{j_2}\}_{j_1}\}_{j_1})$ & 0 & 1 & 2 \\
    $\textbf{t}_{2.10}$ & $(\{\sigma_{j_1},\tau\}_{j_1})$ & 0 & 1 & 2 & $\textbf{t}_{2.19}$ & $(\{\{\{\sigma_{j_2}\}_{j_2}\}_{j_1}\}_{j_1})$ & 0 & 1 & 2 \\
\end{longtable}
\begin{Bem}
    If we choose $j_1=j_2$ then some trees of Table~\ref{tabelle1} may
    coincide. In this case $\alpha_*$ has to be taken as the sum
    of the values $\alpha_*$ from the coinciding trees. As an
    example, for $j_1=j_2$ we get $\alpha_I(\textbf{t}_{2.15}) =
    0+2$ and $\alpha_S(\textbf{t}_{2.15})=1+2$.
\end{Bem}
\begin{longtable}[c]{|c|c|c|c|c||c|c|c|c|c|}
\caption{The correlation coefficient $\beta(\textbf{t})$ for some
trees $\textbf{t} \in TS(*)$, $* \in \{I,S\}$, and $j_1,j_2 \in
\{1, \ldots, m\}$. For trees with $\rho(\textbf{t}) \leq 2.5$
which are not listed holds $\beta(\textbf{t})=1$.}
\label{tabelle1b} \\
    \hline
    $\textbf{t}$ & correlation & $\alpha_I$ & $\alpha_S$ & $\beta$ &
    $\textbf{t}$ & correlation & $\alpha_I$ & $\alpha_S$ & $\beta$ \\
    \hline
    \hline
    \endfirsthead
    \hline
    $\textbf{t}$ & correlation & $\alpha_I$ & $\alpha_S$ & $\beta$ &
    $\textbf{t}$ & correlation & $\alpha_I$ & $\alpha_S$ & $\beta$ \\
    \hline
    \hline
    \endhead
%
    \hline
    \endfoot
%
    \hline
    \endlastfoot
%
    $\textbf{t}_{2.11}$ & $j_1 \neq j_2$ & 1 & 1 & 3 & $\textbf{t}_{2.12b}$ & $j_1 \neq j_2$ & 4 & 4 & 2 \\
    $\textbf{t}_{2.11}$ & $j_1 = j_2$ & 1 & 1 & 1 & $\textbf{t}_{2.12}$ & $j_1 = j_2$ & 4 & 6 & 1 \\
    $\textbf{t}_{2.13b}$ & $j_1 \neq j_2$ & 0 & 2 & 2 & $\textbf{t}_{2.16}$ & $j_1 \neq j_2$ & 0 & 1 & 3 \\
    $\textbf{t}_{2.13}$ & $j_1 = j_2$ & 2 & 4 & 1 & $\textbf{t}_{2.16}$ & $j_1 = j_2$ & 0 & 1 & 1
\end{longtable}
\begin{longtable}[c]{|c|c|c||c|c|c|}
\caption{Trees $\textbf{t} \in TS(\Delta)$ of order
$\rho(\textbf{t}) \leq 2.5$ with arbitrary choice of $j_1, j_2,
j_3, j_4, j_5 \in \{1,\ldots,m\}$.} \label{tabelle2}
\\
    \hline
    $\textbf{t}$ & tree & $\alpha_{\Delta}$ &
    $\textbf{t}$ & tree & $\alpha_{\Delta}$ \\
    \hline
    \hline
    \endfirsthead
    \hline
    $\textbf{t}$ & tree & $\alpha_{\Delta}$ &
    $\textbf{t}$ & tree & $\alpha_{\Delta}$ \\
    \hline
    \hline
    \endhead
%
    \hline
    \endfoot
%
    \hline
    \endlastfoot
    $\textbf{t}_{0.1}$ & $\gamma$ & 1 & $\textbf{t}_{0.5.1}$ & $(\sigma_{j_1})$ & 1 \\
    \hline
    $\textbf{t}_{1.1}$ & $(\tau)$ & 1 & $\textbf{t}_{1.2}$ & $(\sigma_{j_1},\sigma_{j_2})$ & 1 \\
    $\textbf{t}_{1.3}$ & $(\{\sigma_{j_2}\}_{j_1})$ & 1 & & & \\
    \hline
    $\textbf{t}_{1.5.1}$ & $([\sigma_{j_1}])$ & 1 & $\textbf{t}_{1.5.2}$ & $(\{\tau\}_{j_1})$ & 1 \\
    $\textbf{t}_{1.5.3}$ & $(\tau,\sigma_{j_1})$ & 2 & $\textbf{t}_{1.5.4}$ & $(\sigma_{j_1},\sigma_{j_2},\sigma_{j_3})$ & 1 \\
    $\textbf{t}_{1.5.5}$ & $(\{\sigma_{j_2}\}_{j_1},\sigma_{j_3})$ & 3 & $\textbf{t}_{1.5.6}$ & $(\{\sigma_{j_2},\sigma_{j_3}\}_{j_1})$ & 1 \\
    $\textbf{t}_{1.5.7}$ & $(\{\{\sigma_{j_3}\}_{j_2}\}_{j_1})$ & 1 & & & \\
    \hline
    $\textbf{t}_{2.1}$ & $([\tau])$ & 1 & $\textbf{t}_{2.2}$ & $(\tau,\tau)$ & 1 \\
    $\textbf{t}_{2.3}$ & $([\{\sigma_{j_2}\}_{j_1}])$ & 1 & $\textbf{t}_{2.4}$ & $([\sigma_{j_1},\sigma_{j_2}])$ & 1 \\
    $\textbf{t}_{2.5}$ & $(\sigma_{j_1},[\sigma_{j_2}])$ & 3 & $\textbf{t}_{2.6}$ & $(\{\sigma_{j_2}\}_{j_1},\tau)$ & 3 \\
    $\textbf{t}_{2.7}$ & $(\sigma_{j_1},\sigma_{j_2},\tau)$ & 3 & $\textbf{t}_{2.8}$ & $(\sigma_{j_1},\{\tau\}_{j_2})$ & 3 \\
    $\textbf{t}_{2.9}$ & $(\{\{\tau\}_{j_2}\}_{j_1})$ & 1 & $\textbf{t}_{2.10}$ & $(\{\sigma_{j_2},\tau\}_{j_1})$ & 2 \\
    $\textbf{t}_{2.11}$ & $(\sigma_{j_1},\sigma_{j_2},\sigma_{j_3},\sigma_{j_4})$ & 1 & $\textbf{t}_{2.12}$ & $(\sigma_{j_1},\sigma_{j_2},\{\sigma_{j_4}\}_{j_3})$ & 6 \\
    $\textbf{t}_{2.13}$ & $(\sigma_{j_1},\{\sigma_{j_3},\sigma_{j_4}\}_{j_2})$ & 4 & $\textbf{t}_{2.14}$ & $(\sigma_{j_1},\{\{\sigma_{j_4}\}_{j_3}\}_{j_2})$ & 4 \\
    $\textbf{t}_{2.15}$ & $(\{\sigma_{j_2}\}_{j_1},\{\sigma_{j_4}\}_{j_3})$ & 3 & $\textbf{t}_{2.16}$ & $(\{\sigma_{j_2},\sigma_{j_3},\sigma_{j_4}\}_{j_1})$ & 1 \\
    $\textbf{t}_{2.17}$ & $(\{\sigma_{j_2},\{\sigma_{j_4}\}_{j_3}\}_{j_1})$ & 3 & $\textbf{t}_{2.18}$ & $(\{\{\sigma_{j_3},\sigma_{j_4}\}_{j_2}\}_{j_1})$ & 1 \\
    $\textbf{t}_{2.19}$ & $(\{\{\{\sigma_{j_4}\}_{j_3}\}_{j_2}\}_{j_1})$ & 1 & $\textbf{t}_{2.20}$ & $(\{[\sigma_{j_2}]\}_{j_1})$ & 1 \\
    \hline
    $\textbf{t}_{2.5.1}$ & $(\tau,\tau,\sigma_{j_1})$ & 3 & $\textbf{t}_{2.5.2}$ & $([\sigma_{j_1}],\tau)$ & 3 \\
    $\textbf{t}_{2.5.3}$ & $([\tau],\sigma_{j_1})$ & 3 & $\textbf{t}_{2.5.4}$ & $([\tau,\sigma_{j_1}])$ & 2 \\
    $\textbf{t}_{2.5.5}$ & $([[\sigma_{j_1}]])$ & 1 & $\textbf{t}_{2.5.6}$ & $(\{\tau\}_{j_1},\tau)$ & 3 \\
    $\textbf{t}_{2.5.7}$ & $(\{\tau,\tau\}_{j_1})$ & 1 & $\textbf{t}_{2.5.8}$ & $(\{[\tau]\}_{j_1})$ & 1 \\
    $\textbf{t}_{2.5.9}$ & $([\{\tau\}_{j_1}])$ & 1 & $\textbf{t}_{2.5.10}$ & $(\tau,\sigma_{j_1},\sigma_{j_2},\sigma_{j_3})$ & 4 \\
    $\textbf{t}_{2.5.11}$ & $([\sigma_{j_1}],\sigma_{j_2},\sigma_{j_3})$ & 6 & $\textbf{t}_{2.5.12}$ & $(\tau,\{\sigma_{j_2}\}_{j_1},\sigma_{j_3})$ & 12 \\
    $\textbf{t}_{2.5.13}$ & $([\sigma_{j_1},\sigma_{j_2}],\sigma_{j_3})$ & 4 & $\textbf{t}_{2.5.14}$ & $([\{\sigma_{j_2}\}_{j_1}],\sigma_{j_3})$ & 4 \\
    $\textbf{t}_{2.5.15}$ & $([\sigma_{j_1}],\{\sigma_{j_3}\}_{j_2})$ & 6 & $\textbf{t}_{2.5.16}$ & $(\tau,\{\sigma_{j_2},\sigma_{j_3}\}_{j_1})$ & 4 \\
    $\textbf{t}_{2.5.17}$ & $(\tau,\{\{\sigma_{j_3}\}_{j_2}\}_{j_1})$ & 4 & $\textbf{t}_{2.5.18}$ & $([\sigma_{j_1},\sigma_{j_2},\sigma_{j_3}])$ & 1 \\
    $\textbf{t}_{2.5.19}$ & $([\sigma_{j_1},\{\sigma_{j_3}\}_{j_2}])$ & 3 & $\textbf{t}_{2.5.20}$ & $([\{\sigma_{j_2},\sigma_{j_3}\}_{j_1}])$ & 1 \\
    $\textbf{t}_{2.5.21}$ & $([\{\{\sigma_{j_3}\}_{j_2}\}_{j_1}])$ & 1 & $\textbf{t}_{2.5.22}$ & $(\{\tau\}_{j_1},\sigma_{j_2},\sigma_{j_3})$ & 6 \\
    $\textbf{t}_{2.5.23}$ & $(\{\tau\}_{j_1},\{\sigma_{j_3}\}_{j_2})$ & 6 & $\textbf{t}_{2.5.24}$ & $(\{\tau,\sigma_{j_2}\}_{j_1},\sigma_{j_3})$ & 8 \\
    $\textbf{t}_{2.5.25}$ & $(\{[\sigma_{j_2}]\}_{j_1},\sigma_{j_3})$ & 4 & $\textbf{t}_{2.5.26}$ & $(\{\tau,\sigma_{j_2},\sigma_{j_3}\}_{j_1})$ & 3 \\
    $\textbf{t}_{2.5.27}$ & $(\{\tau,\{\sigma_{j_3}\}_{j_2}\}_{j_1})$ & 3 & $\textbf{t}_{2.5.28}$ & $(\{[\sigma_{j_2}],\sigma_{j_3}\}_{j_1})$ & 3 \\
    $\textbf{t}_{2.5.29}$ & $(\{[\sigma_{j_2},\sigma_{j_3}]\}_{j_1})$ & 1 & $\textbf{t}_{2.5.30}$ & $(\{[\{\sigma_{j_3}\}_{j_2}]\}_{j_1})$ & 1 \\
    $\textbf{t}_{2.5.31}$ & $(\{\{\tau\}_{j_2}\}_{j_1},\sigma_{j_3})$ & 4 & $\textbf{t}_{2.5.32}$ & $(\{\{\tau\}_{j_2},\sigma_{j_3}\}_{j_1})$ & 3 \\
    $\textbf{t}_{2.5.33}$ & $(\{\{\tau,\sigma_{j_3}\}_{j_2}\}_{j_1})$ & 2 & $\textbf{t}_{2.5.34}$ & $(\{\{[\sigma_{j_3}]\}_{j_2}\}_{j_1})$ & 1 \\
    $\textbf{t}_{2.5.35}$ & $(\{\{\{\tau\}_{j_3}\}_{j_2}\}_{j_1})$ & 1 & $\textbf{t}_{2.5.36}$ & $(\sigma_{j_1},\sigma_{j_2},\sigma_{j_3},\sigma_{j_4},\sigma_{j_5})$ & 1 \\
    $\textbf{t}_{2.5.37}$ & $(\sigma_{j_1},\sigma_{j_2},\sigma_{j_3},\{\sigma_{j_5}\}_{j_4})$ & 10 & $\textbf{t}_{2.5.38}$ & $(\sigma_{j_1},\sigma_{j_2},\{\sigma_{j_4},\sigma_{j_5}\}_{j_3})$ & 10 \\
    $\textbf{t}_{2.5.39}$ & $(\sigma_{j_1},\sigma_{j_2},\{\{\sigma_{j_5}\}_{j_4}\}_{j_3})$ & 10 & $\textbf{t}_{2.5.40}$ & $(\sigma_{j_1},\{\sigma_{j_3}\}_{j_2},\{\sigma_{j_5}\}_{j_4})$ & 15 \\
    $\textbf{t}_{2.5.41}$ & $(\{\sigma_{j_2}\}_{j_1},\{\sigma_{j_4},\sigma_{j_5}\}_{j_3})$ & 10 & $\textbf{t}_{2.5.42}$ & $(\{\sigma_{j_2}\}_{j_1},\{\{\sigma_{j_5}\}_{j_4}\}_{j_3})$ & 10 \\
    $\textbf{t}_{2.5.43}$ & $(\{\sigma_{j_2},\sigma_{j_3},\sigma_{j_4}\}_{j_1},\sigma_{j_5})$ & 5 & $\textbf{t}_{2.5.44}$ & $(\{\sigma_{j_2},\{\sigma_{j_4}\}_{j_3}\}_{j_1},\sigma_{j_5})$ & 15 \\
    $\textbf{t}_{2.5.45}$ & $(\{\{\sigma_{j_3},\sigma_{j_4}\}_{j_2}\}_{j_1},\sigma_{j_5})$ & 5 & $\textbf{t}_{2.5.46}$ & $(\{\{\{\sigma_{j_4}\}_{j_3}\}_{j_2}\}_{j_1},\sigma_{j_5})$ & 5 \\
    $\textbf{t}_{2.5.47}$ & $(\{\sigma_{j_2},\sigma_{j_3},\sigma_{j_4},\sigma_{j_5}\}_{j_1})$ & 1 & $\textbf{t}_{2.5.48}$ & $(\{\sigma_{j_2},\sigma_{j_3},\{\sigma_{j_5}\}_{j_4}\}_{j_1})$ & 6 \\
    $\textbf{t}_{2.5.49}$ & $(\{\sigma_{j_2},\{\sigma_{j_4},\sigma_{j_5}\}_{j_3}\}_{j_1})$ & 4 & $\textbf{t}_{2.5.50}$ & $(\{\sigma_{j_2},\{\{\sigma_{j_5}\}_{j_4}\}_{j_3}\}_{j_1})$ & 4 \\
    $\textbf{t}_{2.5.51}$ & $(\{\{\sigma_{j_3}\}_{j_2},\{\sigma_{j_5}\}_{j_4}\}_{j_1})$ & 3 & $\textbf{t}_{2.5.52}$ & $(\{\{\sigma_{j_3},\sigma_{j_4},\sigma_{j_5}\}_{j_2}\}_{j_1})$ & 1 \\
    $\textbf{t}_{2.5.53}$ & $(\{\{\sigma_{j_3},\{\sigma_{j_5}\}_{j_4}\}_{j_2}\}_{j_1})$ & 3 & $\textbf{t}_{2.5.54}$ & $(\{\{\{\sigma_{j_4},\sigma_{j_5}\}_{j_3}\}_{j_2}\}_{j_1})$ & 1 \\
    $\textbf{t}_{2.5.55}$ & $(\{\{\{\{\sigma_{j_5}\}_{j_4}\}_{j_3}\}_{j_2}\}_{j_1})$ & 1 & & & \\
    \hline
\end{longtable}

\end{document}